\documentclass[11pt]{article}
\usepackage{amsmath}
\usepackage{amsthm}
\input{epsf.tex}

\newtheorem{theorem}{Theorem}[section]

\newtheorem{lemma}[theorem]{Lemma}

\newtheorem{define}[theorem]{Definition}

\def\Empty{}

\catcode`\@=11

\def\section{\@startsection {section}{1}{\z@}{-3.5ex plus -1ex minus 
-.2ex}{2.3ex plus .2ex}{\large\bf}}

\def\fnum@figure{{\small Figure \thefigure}}
\def\fakefigure{\def\@captype{figure}}

\long\def\@makecaption#1#2{
    \vskip 10pt 
    \def\FCap{#2} \def\NoCap{\ignorespaces}
    \ifx \FCap\NoCap
       \setbox\@tempboxa\hbox{#1}  
      \else
       \setbox\@tempboxa\hbox{#1: \small \it #2}
    \fi
    \ifdim \wd\@tempboxa >\hsize   
        \unhbox\@tempboxa\par      
      \else                        
        \hbox to\hsize{\hfil\box\@tempboxa\hfil}  
    \fi}

\def\@oddhead{\hbox{}\rightmark \hfil \rm\thepage}
\def\sectionmark#1{\markright {\sc{\ifnum \c@secnumdepth >\z@
      \S\thesection.\hskip 1em\relax \fi #1}}}

\catcode`\@=12

\def\oplabel#1{
  \def\OpArg{#1} \ifx \OpArg\Empty {} \else
  	\label{#1}
  \fi}

\newlength{\saveu}

\catcode`\@=11

\newcommand{\pa}{\mbox{$<_{\alpha}$}}
\newcommand{\pb}{\mbox{$<_{\beta}$}}
\newcommand{\pu}{\mbox{$<_u$}}
\newcommand{\poa}{\mbox{$>_{\alpha}$}}

\newcommand{\kk}{\mbox{${\cal K}$}}

\newcommand{\aaa}{\mbox{${\cal A}$}}
\newcommand{\bb}{\mbox{${\cal B}$}}
\newcommand{\ak}{\mbox{${\cal A}_{\kappa}$}}
\newcommand{\at}{\mbox{${\cal A}_{\tau}$}}
\newcommand{\ap}{\mbox{${\cal A}_{\alpha}$}}
\newcommand{\ab}{\mbox{${\cal A}_{\beta}$}}

\newcommand{\lat}{\mbox{${\cal L A}_{\tau}$}}

\newcommand{\lal}{\mbox{${\cal L A}_{\alpha}$}}
\newcommand{\lag}{\mbox{${\cal L A}_{\gamma}$}}
\newcommand{\lab}{\mbox{${\cal L A}_{\beta}$}}
\newcommand{\la}{\mbox{${\cal L A}$}}

\newcommand{\uu}{\mbox{${\cal U}$}}
\newcommand{\ro}{\mbox{${\cal R}$}}
\newcommand{\vv}{\mbox{${\cal V}$}}
\newcommand{\we}{\mbox{${\cal W}$}}
\newcommand{\re}{\mbox{${\cal R}$}}

\newcommand{\rrrr}{\mbox{${\bf R}$}}

\newcommand{\mi}{\mbox{$\widetilde M$}}
\newcommand{\wl}{\mbox{$\widetilde \lambda$}}

\newcommand{\gl}{\mbox{${\cal G}$}}

\newcommand{\fol}{\mbox{$\cal F$}}

\newcommand{\fn}{\mbox{$\widetilde {\cal F}$}}

\def\@evenhead{\rm \leftmark \hfil \thepage}
\def\chaptermark#1{\markboth {\sc {\ifnum \c@secnumdepth >\m@ne
      \@chapapp\ \thechapter. \ \fi #1}}{}}%

\catcode`\@=12

\def\centeredepsfbox#1{\centerline{\epsfbox{#1}}}

\begin{document}

\title{Laminar free hyperbolic $3$-manifolds}
\author{S\'{e}rgio R. Fenley
\thanks{Reseach partially supported by NSF grant 
DMS 0296139.}
}
\maketitle

\vskip .2in

\section{Introduction}

We analyse the existence question for essential
laminations in $3$-manifolds.
The purpose of the article is to prove that there
are infinitely many closed hyperbolic $3$-manifolds which
do not admit essential laminations. This
answers in the negative a fundamental question posed
by Gabai and Oertel when they introduced essential
laminations in \cite{Ga-Oe}, see also \cite{Ga4, Ga5}. 
The proof is obtained
by analysing certain group actions on trees and showing
that certain $3$-manifold groups only have trivial
actions on trees. There are corollaries concerning
the existence question for Reebless foliations and
pseudo-Anosov flows.

This article deals with the topological structure of
$3$-manifolds. Two dimensional manifolds are extremely
well behaved in the sense that the universal cover
is always either the plane or the sphere (for closed
manifolds), the fundamental group determines the manifold
and many other important properties. Similarly for
a $3$-manifold one asks: When is the universal cover
$\rrrr^3$? When does the fundamental group determine
the manifold?
Are homotopic homeomorphisms always isotopic? 
An obvious necessary condition is that the manifold
be {\em irreducible}, that is, every embedded sphere bounds
a ball. As for $2$-manifolds, the existence of a compact
codimension one object which is topologically good 
is extremely
useful. 
A properly embedded $2$-sided surface 
not ${\bf S}^2, {\bf D}^2$ is {\em incompressible}
if it injects in the fundamental group level \cite{He}.
A compact, irreducible  manifold with an incompressible
surface is called {\em Haken}. Fundamental work of Haken
\cite{Hak1,Hak2} and Waldhausen \cite{Wa} shows that
Haken manifolds have
fantastic properties, answering in the positive the
3 questions above.

But how common are Haken $3$-manifolds, that is how
common are incompressible surfaces amongst irreducible
$3$-manifolds? In some sense they are not very common.
Recall that {\em Dehn surgery} along an orientation
preserving  simple closed curve $\delta$ is the process
of removing a tubular neighborhood $N(\delta)$ (a solid torus)
and glueing back by a homeomorphism of the boundary - which
is a two dimensional torus $T_1$ \cite{Rol,Bu-Zi}. 
The surgered manifold is completely determined
by which simple closed curve in $T_1$ becomes the new
{\em meridian}, that is, which curve of $T_1$ is glued to
the null homotopic curve in the boundary of $N(\delta)$.
Hence this is parametrized by a pair of relatively
prime integers $(q,p)$, corresponding to the description
of simple closed curves in $T_1$.  
When viewed this way, 
this set of relatively prime $(q,p)$ is the
{\em Dehn surgery space} $-$ a subset of
${\bf Z}^2 \subset {\rrrr}^2$.
The same can be done iterating the process doing Dehn
surgery on links \cite{He,Rol,Bu-Zi}.
Notice that all closed, orientable $3$-manifolds can be obtained
from ${\bf S}^3$ by some Dehn surgery on an appropriate link
in ${\bf S}^3$ \cite{Rol}. So one can interpret
how common a property is by
verifying how many of the Dehn surgered manifolds have
that property. Along these lines some of the many
results on incompressible surfaces are:
If $K$ is a two bridge knot in ${\bf S}^3$ then almost all
Dehn surgeries on $K$ yield manifolds without incompressible
surfaces \cite{Ha-Th}. The same is true for any knot $K$ in
a manifold $M$ so that $M - K$ does not have any
closed incompressible surfaces \cite{Hat1}.
Notice that there are also results on the other direction:
for example Oertel \cite{Oe} proved that for many star links
in ${\bf S}^3$, 
then any non trivial Dehn surgery
yields a manifold with incompressible surfaces.
There are similar results for Montesinos knots \cite{Ha-Oe}.
Basically a lot of it depends on whether the complement
has closed incompressible surfaces or not.
In many cases the complement does not have
such surfaces, yielding the non existence results
for most Dehn surgered manifolds.

This amongst other reasons led to the concept of an
essential lamination as introduced by Gabai and Oertel
in the seminal paper \cite{Ga-Oe} of the late 80's. 
A {\em lamination} is a foliation of a closed subset
of the manifold. Roughly a lamination in a closed
$3$-manifold is {\em essential} if it has no sphere
leaves or tori leaves bounding solid tori, the complement
of the lamination is irreducible and the leaves in
the boundary of the complement are incompressible
and end incompressible in their respective complementary
components \cite{Ga-Oe}.
Gabai and Oertel proved the fundamental result that
essential laminations have deep consequences: 
the manifold $M$ is irreducible, its universal cover
is $\rrrr^3$, leaves of the lamination inject
in the fundamental group level, efficient closed 
transversals are not null homotopic; amongst
other consequences \cite{Ga-Ka3}. In addition
manifolds with {\em genuine} essential laminations
satify the weak hyperbolization conjecture \cite{Ga-Ka4}: either
there is a ${\bf Z} \oplus {\bf Z}$ subgroup of 
the fundamental group or the fundamental group
is Gromov hyperbolic \cite{Gr,Gh-Ha}.
Genuine means that not all complementary
regions are $I$-bundles, or equivalently it is not
just a blow up of a foliation.
Brittenham also proved properties concerning homotopy
equivalences for manifolds with essential laminations
\cite{Br2}.

In addition essential laminations are extremely common:
For example if $K$ is a non trivial knot in ${\bf S}^3$
then off of at most two lines and a couple of points
in Dehn surgery space,
the surgered manifold contains
an essential lamination.
This is obtained as follows: first Gabai
constructed a Reebless foliation $\fol$
in $({\bf S}^3 - N(K))$ which is transverse to the boundary
\cite{Ga1,Ga2,Ga3}.
Reebless means it does not have a Reeb component:
a foliation of the solid torus with the boundary being
a leaf, all other leaves are planes spiralling to the
boundary \cite{Re,No}.
Then results of Mosher, Gabai \cite{Mo2}
show that either there is
an incompressible torus transverse to $\fol$ or there is
an essential lamination in ${\bf S}^3 - N(K)$ with
solid torus complementary regions. This lamination remains
essential off of at most two lines in Dehn surgery space
\cite{Mo2} - see more on solid torus complementary regions 
later.
Also Brittenham produced examples of essential laminations
which remain essential after all non trivial Dehn surgeries
\cite{Br3,Br4}.  Roberts has also obtained many important
existence results concerning alternating knots in the
sphere 
\cite{Ro1,De-Ro} (partly jointly with
Delman) and punctured surface bundles \cite{Ro2,Ro3}.

So succesful was the search for essential laminations that
at first one might wonder whether all manifolds that
can (irreducible, with infinite fundamental group),
in fact do admit essential laminations. 
Given that an incompressible torus is an essential lamination,
the Geometrization conjecture \cite{Th2} suggests that one should
only have to analyse Seifert fibered spaces and hyperbolic
manifolds \cite{Sc,Th2}. 

The situation for Seifert fibered spaces has
been resolved: Brittenham produced examples of Seifert
fibered spaces  which are irreducible,
have infinite fundamental group, universal cover $\rrrr^3$,
but which do not have essential laminations \cite{Br1}.
Naimi \cite{Na}, using work of Bieri, Neumann and Strebel
\cite{BNS},
completely determined which Seifert fibered manifolds
admit essential laminations.

For hyperbolic $3$-manifolds there were two fundamental 
open questions: 1) (Thurston) Does every closed hyperbolic
$3$-manifold admit a Reebless foliation? 2) (Gabai-Oertel \cite{Ga-Oe},
see also \cite{Ga4,Ga5})
Does every closed hyperbolic $3$-manifold admit
an essential lamination?
In the past year question 1) was answered in the negative
by Roberts, Shareshian and Stein \cite{RSS} who produced
infinitely many counterexamples. The goal of this article
is to answer question 2) in the negative. We now proceed
to describe the examples.

Basically one starts with a torus bundle $M$ over the circle and
do Dehn surgery on a particular closed curve. Let
$\phi$ be the monodromy of the fibration associated to
a $2$ by $2$ integer matrix $A$, so that $A$ is hyperbolic.
Let $R$ be a fiber which is a torus.
There are two foliations in $R$ which are invariant 
under the monodromy $\phi$, the stable and unstable foliations.
The suspension flow in $M$ induces two foliations in $M$
with leaves being planes, annuli and M\"{o}ebius bands.
Suppose there is a M\"{o}ebius band leaf. 
Blow up that leaf, producing a lamination with
a solid torus complementary component with closure a solid
torus with core $\delta$ 
and with some curves $\eta$ removed from the boundary.
The curves $\eta$ are called the {\em degeneracy locus} of
the complementary region of the lamination \cite{Ga-Ka1}.
One can think of $\eta$ as lying in the boundary of $N(\delta)$,
which is a two dimensional torus.
Let $(1,0)$ be the curve in $\partial N(\delta)$ which bounds
the fiber in $M - N(\delta)$. Under an appropriate choice for
the curve $(0,1)$ of $\partial N(\delta)$ then $\eta$ is
represented by $(1,2)$. Do Dehn surgery along $\delta$.
If $\xi$ is the new meridian (the Dehn surgery slope),
then results of essential laminations \cite{Ga-Oe,Ga-Ka1}
show that $\lambda$ remains essential in the
Dehn surgery manifold $M_{\xi}$
if the intersection number of $\xi$ and $\eta$ is at least
$2$ in absolute value. If $\xi$ is described as $(q,p)$ then
this is equivalent to $|p - 2q| \geq 2$.
Therefore the open cases for essential laminations are
$|p - 2q| \leq 1$.

For simplicity of notation we omit the explicit dependence
of $M$ on $\phi$. It is always understood that $M$ depends
on the particular $\phi$.

In a beautiful and fundamental result, Hatcher \cite{Hat2}, showed
that if $p < q$ then then Dehn surgery
manifold $M_{\xi} = M_{p/q}$ has a Reebless foliation. This
is done via an explicit construction involving train tracks
and branched surfaces. In the last year Roberts, Shareshian
and Stein considered a particular type of monodromy,
namely generated by  the matrix

$$A \ \ = \ \
\left[
\begin{array}{rr}
 m & -1 \\ 1 & 0 \\ 
\end{array}
\right] \ \ \ \ m \leq -3$$

The eigenvalues of $A$ are negative. Consider the point
$(0,0)$ in $\rrrr^2$ and its projection $O$ to 
the fibering torus $R$. Let $\delta$ be the closed
orbit of the suspension flow through $O$. Because
the eigenvalues are negative, the leaf of 
the stable foliation through $O$ is a M\"{o}ebius
band. When it is blown open into an annulus the
degeneracy locus is $(1,2)$ as described above.
In a groundbreaking work, Roberts, Shareshian and
Stein \cite{RSS} considered Dehn surgery on these
manifolds and proved a wonderful result: 
if $p$ is odd, $m$ is odd and $p \geq q$ then
$M_{p/q}$ does not admit Reebless foliations.
In this article we consider a subclass of these
manifolds and show they do not admit essential 
laminations:

\vskip .2in
\noindent
{\bf Main Theorem:} \ Let $M$ be a torus bundle over the
circle with monodromy induced by the matrix $A$ above.
Let $\delta$ be the orbit of the suspension flow coming
from the origin and $M_{(q,p)} = M_{p/q}$ be the manifold obtained
by $(q,p)$ Dehn surgery on $\delta$. Here $(1,0)$ 
bounds the fiber in \ $M - N(\delta)$ \ and $(1,2)$ is
the degeneracy locus. 
Then if $m \leq -4$ and $|p - 2q| = 1$, the
manifold $M_{p/q}$ does not admit essential laminations.
\vskip .18in

The manifold $M - \delta$ is atoroidal \cite{Th4,Bl-Ca}
and fibers over the circle
with fiber a punctured torus.
By Thurston's hyperbolization
theorem in the fibering case $M - \delta$ has
a complete hyperbolic structure of finite volume \cite{Th3}.
By Thurston's Dehn surgery theorem $M_{p/q}$
is hyperbolic for almost all $p/q$ \cite{Th1}. Therefore:

\vskip .15in
\noindent
{\bf Corollary:} \ There are infinitely many closed, hyperbolic
$3$-manifolds which do not admit essential laminations.
\vskip .13in

We mention that Calegari and Dunfield \cite{Ca-Du}
approached the existence problem from a different point
of view. Following ideas
and results of Thurston \cite{Th5,Th6} concerning
the universal circle for foliations they showed
that 
a wide class of essential laminations also possess
a universal circle. 
One consequence is that tight essential laminations
with torus guts (see \cite{Ca-Du} for definitions)
have universal circles.
Hence the fundamental groups
act on the circle. Under certain conditions and if
the fundamental group is orderable then the 
action lifts to a non trivial action in $\rrrr$ and
they obtain nonexistence results for these types of laminations.
See more below.

Another immediate corollary is:

\vskip .15in
\noindent
{\bf Corollary:} \ If $m \leq -4$ and $|p-2q| = 1$, then
the manifolds $M_{p/q}$ above do not admit Reebless foliations.
\vskip .13in

About half of this result has already been established by
Roberts, Stein and Shareshian \cite{RSS}, namely the situation when
$m$ is odd. See more on $m$ odd below.
Another consequence is:

\vskip .2in
\noindent
{\bf Corollary:} \ If $m \leq -4$ and $|p-2q| = 1$ then
$M_{p/q}$ does not admit pseudo-Anosov flows.
\vskip .2in

For basic definitions and properties of pseudo-Anosov flows
consult \cite{Mo1,Mo2}. This result provides infinitely
many hyperbolic manifolds without pseudo-Anosov flows.
We stress that 
Calegari and Dunfield \cite{Ca-Du}
previously obtained conditions implying manifolds
do not admit pseudo-Anosov flows and showed for example that
the Weeks manifold does not admit pseudo-Anosov flows.

We remark that Dehn surgery on torus bundles over the circle
has been widely studied, for example: a) Which surgered 
manifolds have incompressible surfaces \cite{Fl-Ha,CJR},
b) Virtual homology \cite{Bk1,Bk2}, c) geometrization
\cite{Jo,Th1,Th2,Th3,Th4}.

\vskip .2in
We now describe the key ideas of the proof of the
main theorem.
The proof is done by looking at group actions on trees.
For simplicity first consider the case of a
Reebless foliation $\fol$. 
Novikov proved that leaves are incompressible
and transversals are never null homotopic \cite{No}.
Hence the lift to the universal
cover $\fn$ is a foliation by planes and its
leaf space is a simply connected $1$-dimensional
manifold, which may not be Hausdorff.
This $1$-manifold can be collapsed to a tree.
The fundamental group acts on this tree.
Roberts et al analysed group actions 
on simply connected non Hausdorff $1$-manifolds and also on trees
$-$ under the conditions
$p \geq q$ and $p, m$ odd, they ruled out the existence
of Reebless foliations \cite{RSS}.
Notice that the leaf space of the lifted foliation
$\fn$ is an orientable 
object and it makes sense to talk about orientation
preserving homeomorphisms. In order to stay in
the orientation preserving world they restricted
to $p, m$ odd.

Now consider an essential lamination $\lambda$.
The results of Gabai and Oertel \cite{Ga-Oe} 
imply that the lift to the universal cover
$\wl$ is a lamination by planes
in $\mi$.
To get the leaf space
blow down closures of complementary regions
to points and also non isolated leaves (on both
sides) to points. This produces an order 
tree as defined by Gabai-Kazez \cite{Ga-Ka2}
also called a non Hausdorff tree in this
situation \cite{Fe}.
A further appropriate collapsing of the (possible)
non Hausdorff points yields an actual tree
where the group acts non trivially.
The strategy is to show there are no nontrivial
actions of the group on trees.
An action is {\em trivial} if it has a global
fixed point.
A crucial difference from the case of foliations
is that in the case of laminations the
tree does not have an orientation in general.
Hence orientation dependent arguments cannot be used.
This was very important and widely used in \cite{RSS}.
Since we do not have an orientation here, the condition
$m$ odd does not play a role, which allows us
to consider $m$ even as well. 
Notice that the
condition $|p - 2q| =1$ obviously implies that
$p$ is odd. On the other hand there are many examples
with $p$ even so that $M_{p/q}$ has a Reebless foliation
- for example $p = 4, q = 1$ or
 $p = 8, q = 3$ (this has $p > q$!).
So to rule out Reebless foliations,
some further condition on $p, q$ should be necessary
when $p$ is even.
Except for ruling out trivial actions,
the proof here is done entirely in the tree $-$
we never go back to the original non Hausdorff tree.
For the sake of completeness we state this result
from which the main theorem is an easy corollary:

\vskip .2in
\noindent
{\bf {Theorem:}}  \ Let $M_{p/q}$ be the manifold
described above. If $m \leq -4$ and $|p - 2q| = 1$,
then every action of $\pi_1(M_{p/q})$ on a tree
is trivial.
\vskip .2in

The fundamental group of $M_{p/q}$ denoted by
$\gl$ can be generated by two elements
$\alpha$ and $\tau$. 
Actions of a homeomorphism on a tree are easy
to understand: either there is a fixed point
or in the free case there is an invariant axis.
An axis is a properly embedded copy of the reals
where the homeomorphism acts by translation.
The proof breaks down as to whether the generators
above act freely or not yielding 3 main cases
to consider (when $\tau$ acts freely it does
not matter the behavior of $\alpha$).
The proof subdivides into various subcases.
Invariably the analysis goes like this:
apply a certain relation in the group to a well
chosen point. One side of the relation
implies the image of the point is in a certain
region of the tree while the other side of the
relation implies it is in a different region - contradiction!
An important idea is that of a {\em local axis},
which has all the properties of axis except perhaps
being properly embedded. Homeomorphisms with fixed points
may have local axes. This is extremely useful
in a variety of cases.

We note that ${\bf Z}$ actions on non Hausdorff  trees
had been previously analysed in \cite{Fe} and \cite{Ro-St1,Ro-St2},
with consequences for pseudo-Anosov flows \cite{Fe}
and Seifert fibered spaces \cite{Ro-St1,Ro-St2}.

There is a large literature of group actions
on trees which were brought to the
forefront by Serre's fundamental monograph 
\cite{Se}. The analyis usually involve a metric
which is  invariant under the actions
\cite{Mo-Sh1,Mo-Sh2,Mo-Sh3} or actions on simplicial
trees \cite{Se}.
We stress that the tree involved in here is 
not simplicial and it is not presented in
general with a group invariant metric $-$ unless
there is a holonomy invariant transverse measure
of full support in the lamination, e.g when
there is an incompressible surface.
So the proof is entirely topological and
in that sense elementary.
The topology of the manifold, particularly
the condition $|p-2q| = 1$ plays a crucial role.
Notice that in the foliations case there
is a pseudo-metric lying in the background which is
used from time to time in the proof by Roberts
et al \cite{RSS}. The pseudometric distance
between two points measures
how many jumps between non separated points
are necessary to go from one point to the other.
This pseudometric was analysed and used previously by
Barbot in \cite{Ba1,Ba2} with consequences for
foliations.

The results of this article mean that the search for
structures more general than essential laminations,
but still useful takes an added relevance.
One idea previously proposed by Gabai \cite{Ga5} is
that of a {\em loosesse lamination}. We will have more comments
on that in the final remarks section.

We are very thankful to Rachel Roberts who introduced
the idea of considering group actions in
the foliations case and other ideas.

\section{The group}
\label{grou}

Here we compute the fundamental group of $M_{p/q}$.
Start with $M$ the torus bundle over the circle
with monodromy induced by 

$$A \ \ = \ \ \left[ \begin{array}{rr} 
m & -1 \\ 1 & 0 \end{array} \right] \ \ \ \ \ 
{\rm where} \ \ m \leq -3$$

For notational simplicity the dependence of $M$ on $A$ is omitted.

The eigenvalues of $A$ are 

$$\frac{m \pm \sqrt{m^2 - 4}}{4}$$

\noindent
which are both negative  and the matrix is
hyperbolic. The eigenvector directions produce two
linear 
foliations in $\rrrr^2$ with irrational slope
and invariant under $A$. They induce two foliations
in the torus $T^2$. Since $A$ is integral it induces
a homeomorphism $\phi$ of $T^2$, which leaves
the foliations invariant. Let $O$ in $T^2$ be the
image of the origin. Let $M$ be the suspension
of $\phi$ and let $\fol$ be (say) the suspension
of the stable foliation of $T^2$. Then $\fol$ has
leaves which are planes, annuli and M\"{o}ebius bands.
Identify $T^2$ with a fiber in $M$ and let $\delta$
be the orbit through $O$, which is a closed orbit
intersecting $T^2$ once.
Since the eigenvalues of $A$ are negative, the stable
leaf containing $\delta$ is a M\"{o}ebius band. 
We do Dehn surgery on $\delta$. 
We first determine the fundamental group of $M - N(\delta)$.
To do that 
let 

$$D \ = \ N(\delta) \cap T^2 \ \ {\rm (a \ disk)}, \ \ \ 
V \ = \ T^2 - D \ \ {\rm (a \ punctured \ torus)}.$$

\noindent
Choose a basis for the homology of $\partial N(\delta) 
= T_1$, a torus. 
Let $(1,0)$ be the curve in $T_1$ bounding the fiber $V$ of
$M - N(\delta)$.
Blow up the leaf of $\fol$ through $\delta$. It blows
to a single annulus and the complementary region
is a solid torus with core $\delta$. 
The closure of the complementary region 
is a solid torus with a  closed curve in the boundary removed.
The removed curve is the degeneracy locus of the
complementary component \cite{Ga-Ka1}.
Since the leaf of $\fol$ was a M\"{o}ebius band,
the degeneracy locus intersects the curve
$(1,0)$ twice. Choose the curve $(0,1)$ so that
the degeneracy locus is the curve $(1,2)$ in this basis.
Let $M_{p/q}$ be the manifold obtained from
$M$ by doing $(q,p)$ Dehn surgery on $\delta$.
By results about essential laminations, the lamination
$\lambda$ remains essential in $M_{p/q}$ if
$|p - 2q| \geq 2$.
Let $\gamma$ be the curve $(0,1)$ in $T_1$ and
$\tau$ be the curve $(1,0)$.
The degeneracy locus is the curve $\gamma \tau^2$.
Notice there are two tori here: one
which is a fiber of the original fibration (here denoted
by $T^2$), another which
is the boundary of $N(\delta)$ (here denoted by $T_1$).
The Dehn surgery coefficients refer to $T_1$.

\begin{figure}
\centeredepsfbox{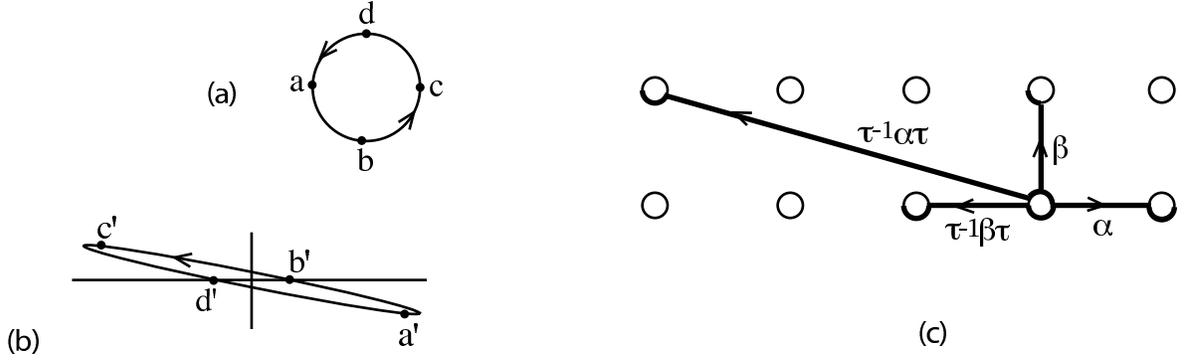}
\caption{
Computing the fundamental group of $M - N(\delta)$.}
\label{01}
\end{figure}

Suppose the disk $D$ above
is a round disk of radius $\epsilon$ sufficiently
small.
The universal abelian cover of $T^2 - D$ is the plane
with disks of radius $\epsilon$ around integer lattice
points removed. Let $E$ be the one around the origin.
We pick $4$ points in $\partial E$:
$a = (-\epsilon, 0),
b = (0, - \epsilon),
c = (\epsilon, 0)$ and $d = (0, \epsilon)$,
see fig. \ref{01}, a.
Let $a'$ be the image of $a$ under $A$, etc.., see fig.
\ref{01}, b.

The image of $\partial E$ under $A$ is an ellipse
which can be deformed back to $\partial E$, see fig.
\ref{01}, b.
Notice $b', d'$ are in the $x$ axis
and $d' = a$.

Let the image of $a$ in $T^2 - D$ be the basepoint
of the fundamental group of $M - N(\delta)$
for simplicity still denoted by $a$ and likewise
for $b, c, d$.
Let $l$ be an arc along the image of $\partial E$
under $A$, going counterclockwise from $d'$ to $a'$.


We pick a basis for $\pi_1(T^2 - D)$:
Let \ $\alpha  \ =  \ \overline{ac} * l_1$ \ (see fig. \ref{01}, c)
where the arc $\overline{ac} \subset \partial E$
is traversed in the counteclockwise direction
and $l_1$ is parametrized as \ $\{ (t,0) \ | \
\epsilon \leq t \leq 1 - \epsilon\}$.
Here $*$ denotes concatenation of arcs.
Let also 

$$\beta \ = \ \overline{ad}_{clo} * l_2 * \overline{ba}_{clo},$$

\noindent
where $l_2$ is parametrized as \ $\{ (0,t) \ | \
\epsilon \leq t \leq 1 - \epsilon \}$,
and the ``clo" subscript means the arcs are traversed
clockwise in $\partial E$.
We identify $\alpha$ and $\beta$ with their images
in $T_2 - D$, so they generate the fundamental group
of $T_2 - D$.
It is easy to see that 
the curve $\gamma = [\alpha,\beta]$ is just 
an counterclockwise turn around $\partial E$.
Then 

$$ \tau^{-1} \alpha \tau \ \ = \ \ l * \overline{a'c'} * l'_1 
* l^{-1}.$$

\noindent
The composition $l * \overline{a'c'}$ is roughly one 
counterclockwise turn around
$\partial E$ so it is the curve $\gamma$.
The straight arc $l'_1$ goes from
$c' = (m \epsilon, \epsilon)$ to
$(m(1-\epsilon), 1 - \epsilon)$ - roughly going 
one step up and $|m|$ steps to the left.
This together with $l^{-1}$ can be isotoped to
$\beta \alpha^m$  \ (where we are identifying $\alpha,
\beta$ with the appropriate covering translates).
We conclude that $\tau^{-1} \alpha \tau = \gamma \beta \alpha^m$.
Similarly 

$$\tau^{-1} \beta \tau \ \ = \ \ l *  \overline{a'd'}_{clo} *  l'_2 
*  \overline{b'a'}_{clo} *  l^{-1}$$

\noindent
So in the same way it is easy to see that 
$\tau^{-1} \beta \tau = \alpha^{-1}$.
Notice that $\alpha, \tau$ generate $\pi_1(M - N(\delta))$.
Hence

$$\pi(M - N(\delta)) \ \ = \ \ 
\{ \alpha, \tau \ \ | \ \ \tau^{-1} \alpha \tau = \gamma \beta \alpha^m,
\ \ \tau^{-1} \beta \tau = \alpha^{-1}, \ \gamma = [\alpha, \beta] \ \}$$

\noindent
After $(q,p)$ Dehn surgery on $\delta$ we obtain
$q \gamma + p \tau$ is the new meridian or
$\tau^p \gamma^q = 1$. Hence

$$\gl \ \ = \ \ \pi_1(M_{p/q}) \ \ =   \ \ 
\{ \alpha, \tau \ | \ \tau^{-1} \alpha \tau = \gamma \beta \alpha^m,
\ \ \tau^{-1} \beta \tau = \alpha^{-1}, \ \gamma = [\alpha, \beta],
\ \tau^p \gamma^q = 1 \}$$

In the proof we will use these and
the following variations of these
relations extensively:

$$\tau^{-1} \beta \tau = \alpha^{-1}, \ \ \tau \alpha \tau^{-1} = \beta,$$

$$ \tau^{-1} \alpha \tau \ = \ 
\gamma \beta \alpha^m = \alpha \beta \alpha^{m-1}
\ = \ \alpha \tau \alpha^{-1} \tau^{-1} \alpha^{m-1}$$

$$
\alpha \tau  = \tau \gamma \beta \alpha^m = \tau \alpha \beta \alpha^{m-1}$$

$$\alpha \beta = \gamma \beta \alpha, 
\ \ {\rm or} \ \ \alpha \tau \alpha^{-1} \tau^{-1}= \tau \gamma 
\alpha^{-1} \tau^{-1} \alpha$$

\noindent
A little manipulation with the relations also yields

$$ \tau \beta \tau^{-1} \ = \ \beta \alpha \beta^{m-1}
\ = \ \gamma^{-1} \alpha \beta^m$$

These and circular variations of these will be used throughtout
the article.

Since $q, p$ are relatively prime there are $e,f$ in ${\bf Z}$ with
$ep + fq = 1$. Let $\kappa = \tau ^f \gamma^{-e}$.
Then $\kappa$ is a generator of the ${\bf Z}$ subgroup
of $\gl$ generated by $\tau, \gamma$ and
$\tau = \kappa^q, \gamma = \kappa^{-p}$.

\vskip .1in
\noindent
{\bf {NOTATION:}} $-$ In the arguments group elements act
on the {\underline {right}}.
\vskip .1in

\section{Outline of the proof}

Given the presentation of $\gl$ above the proof
of the main theorem
is broken into 4 cases:

\begin{itemize}

\item 
Case R - $\rrrr$-covered case

\item 
Case A - $\tau$ acts freely

\item 
Case B - $\alpha$ acts freely $\tau$ has a fixed point.

\item 
Case C - $\alpha$ and $\tau$ have fixed points.

\end{itemize}

If $\mu$ acts freely on a tree, let $\aaa_{\mu}$ be
its axis.
If $\mu$ has a local axis, we denote it by ${\cal L A} _{\mu}$.
Unlike a true axis, a homeomorphism may have more than
one local axis.
The context will make it clear which one is being considered.

\vskip .15in
\noindent
{\underline {Case R}}  $-$ $\rrrr$ covered case.

The $\rrrr$-covered case is simple.
Given that $p$ is odd, this implies that
$\tau$ is orientation preserving in $\rrrr$.
The case $\alpha$ orientation preserving 
is simple. The other case (which implies $m$ is even)
leads to $p > 3q$ which for our purposes is
enough. It also leads us to move away from
orientation preserving arguments, which is more
like the laminations case.
We note that there is an easy linear non trivial
action on $\rrrr$ when $p = 4, q = 1$. Notice that
in this case $p$ is even.

\vskip .15in
\noindent
{\underline {Case A}} $-$ $\tau$ acts freely.

This implies that $\kappa$ also acts freely
and $\aaa_{\kappa} = \aaa_{\tau}$.
We analyse how $\aaa_{\kappa}$ intersects
$\ak \alpha$ and other translates (here $\ak \alpha$
is the image of $\ak$ under $\alpha$).
Let $u = \alpha \beta$. One uses the
relation $\alpha \beta = \gamma \beta \alpha$ to
analyse how $\ak$ intersects $\ak u$ which
breaks down into various cases as to whether
this intersection is empty, a single point
or a segment.
One particularly tricky case needs the 
condition $m \not = -3$.

\vskip .15in
\noindent
{\underline {Case B}} $-$ $\alpha$ acts freely, $\tau$ has
a fixed point.

Let $z$ be a fixed point of $\tau$.
First suppose that $z$ is not in the axis 
$\ap$ 
of $\alpha$.
Suppose there is no fixed point of
$\tau$ between $z$ and $\ap$.
Here let $\uu$ be the component
of $T - \{ z \}$ containing $\ap$.
The case $\uu \tau \not = \uu$ is easy to
deal with.
It follows that $\uu \tau = \uu$
producing a local axis $\lat$ of $\tau$ which
is contained in $\uu$ and has one limit point
in $z$. The proof breaks down as to whether
$\lat$ intersects $\ap$ or not.
Empty intersections are easy to deal with,
the other case being trickier.

Then suppose $z$ is in $\ap$.
We remark this is a crucial case, because this
is likely what happens for the essential laminations
we know to exist when $|p - 2q| \geq 2$.
These come from the original stable lamination
on the fibering manifold. In that manifold,
$\alpha$ acted freely and $\tau$ had a fixed
point in $\ap$. After the surgery $\alpha$ would
still have at least a local axis, which contains
a fixed point of $\tau$. So one knows the exact
condition $|p - 2q| =  1$ will have to be used here!

In this case consider $\uu_1$ be the component
of $T - \{ z \}$ containing $z \alpha$ and
$\uu_2$ the one containing $z \alpha^{-1}$.
It is easy to show that $\uu_1 \tau$ is not
$\uu_1$ and that $\uu_1 \tau$ is in fact
equal to $\uu_2$.
When $\uu_1 \tau^{-1} = \uu_2$ then one
produces a contradiction just using that
$p$ is odd.
The case $\uu_1 \tau^{-1} \not = \uu_2$ or
$\uu_2 \tau \not = \uu_1$ is much more
interesting. Here the  exact condition
$|p-2q| = 1$ is used to show it would imply
$\uu_1 \tau = \uu_1$ which was disallowed
at the beginning.
This actually has connections with the topology
of the situation, see detailed explanation
in section \ref{caseb}.
This is a crucial part of the proof.
One very tricky issue is that a priori
$z$ is only a fixed point of $\tau$ and not
of $\gamma$ $-$ part of the proof is ruling this out.

\vskip .15in
\noindent
{\underline {Case C}} \ 

Generally an axis is good because it gives
information about where points go. The
case of fixed points is trickier
and one many times searches for local axis.

Given two points $a, b$ in a tree let
$[a,b]$ be the unique embedded segment
connecting them. Let $(a,b) = [a,b] - \{ a, b \}$.
Notice $(a,b)$ is exactly the set of points in the
tree separating $a$ from $b$.

Here let $s$ be a fixed point of $\kappa$ and
$w$ a fixed point of $\alpha$ so that there
is no fixed point of either in $(s,w)$.
Notice there may be fixed points of
$\tau$ in $(s,w)$!
Let $\we$ be the component of $T - \{ s \}$ containing
$w$ and $\vv$ the component of $T - \{ w \}$ containing
$s$.
The first part of the proof shows that $\we \tau = \we$ and
$\vv \alpha = \vv$.
These are moderately involved cases.
This immediately produces a local axis 
$\lal$ of $\alpha$ contained in $\vv$ and with
one ideal point $w$.
One does not have yet a local axis for $\tau$ 
because we do not know a priori that $\tau$ has
no fixed points in $(s,w)$.
Some technical complications ensue.

One then shows that 
$s \alpha, s \alpha^{-1}$ are in $\we$.
Let $z$ be the fixed point of $\tau$ in $[s,w)$ which
is closest to $w$ $-$ $z$ could be $s$.
Using the previous results show that
the component $\uu$ of $T - \{ z \}$ containing
$w$ is invariant under $\tau$.
Now this produces a local axis $\lat$ of $\tau$ 
in $\uu$ with ideal point $z$ and some further
properties.
One then shows that $w$ is not 
in $\lat$ and $z$ not in $\lal$.

We are now in familiar ground.
If $\lal \cap \lat$ has at most one
point, then it is easy.
When
$\lal \cap \lat$ has more than one point
we use arguments done in case B $-$ 
this part of the arguments in case B is done
in more generality using local axis (rather than
axis as needed in case B) and can be used
in case C as well. This finishes the proof of case C.
This finally yields the proof of the main theorem.

The arguments in this article are very involved.
One possibility to read the article and get
a quick grasp of the proof is to first analyse
the $\rrrr$-covered proof. Then go to the proof
of case B.2 - $\alpha$ acts freely and $\tau$ has
 a fixed point in the axis of $\alpha$ $-$ this
case admits essential laminations if $|p - 2q| \geq 2$
and the topology can be detected. Then read the
proof of $\tau$ acts freely and the other proofs.

\section{Preliminaries}

Let $\lambda$ be an essential lamination on
a $3$-manifold $N$. We'll modify $\lambda$ if
necessary to eventually obtain a group action
on a tree which is essentially the leaf space
of the lifted lamination $\widetilde \lambda$
to the universal cover $\widetilde N$.
First if there are any leaves
of $\lambda$ which are isolated on both sides,
then blow each of them into an $I$-bundle of
leaves $-$ needs to be done at most countably
many times. Now 
$\widetilde \lambda$ is a lamination by planes
with no leaves isolated on both sides \cite{Ga-Oe}.

Suppose $L$ is a leaf of $\wl$
which is non separated from another leaf $F$
$-$ that is, there are $L_i$ leaves of $\wl$
with $L_i$ converging to both $L$ and $F$.
We do not want that $L$ is not separated from
some other leaf in the other side (the one
not containing $F$). If that happens, blow up
$L$ into an $I$-bundle of leaves.
This can also be achieved by a blow up in $\lambda$.
Since there are at most countably many leaves
non separated from some other leaf we can get
rid of leaves non separated from leaves on
both sides. If needed use blow ups so that
non separated
leaves of $\wl$ are not boundary leaves
of a complementary region of $\wl$ (on the opposite side).
After all these possible modifications assume
this is the original lamination $\lambda$.

Now define a set $T_*$ whose elements are:
closures of complementary components of $\wl$ and
also leaves of $\wl$ which are non isolated on
both sides. Then $T_*$ is an {\em order tree}
\cite{Ga-Ka2,Ro-St2}, also called non
Hausdorff tree \cite{Fe}.
The fundamental group $\pi_1(N)$ naturally acts
on $T_*$. If $e$ is any point of $T_*$ which is
non separated another point $e'$, collapse all
points non separated from $e$ together with
$e$. This is OK since no such $e$ is
non separated on more than one side and $e$ also
does not come from a complementary region of $\wl$.
The collapsed object is now an actual tree $T$ and
the action of $\pi_1(N)$ on $T_*$ induces
a natural action of $\pi_1(N)$ on $T$.
In our proof $N = M_{p/q}$ and we will analyse
group actions of $\gl = \pi_1(M_{p/q})$ on the tree $T$.

\begin{define}{}
A group action on a tree $T$ is nontrivial if no
point of $T$ is fixed by all elements of the group.
\end{define}

A lot of results on group actions on trees are
to rule out non trivial group actions \cite{Cu-Vo}.

Given point $a, b$ on a tree $T$ let

$$(a,b) \ = \ \{ c \in T \ \ | \ \ c \ \ {\rm separates} \ \
a \ \ {\rm from} \ \ b \}.$$

\noindent
If $a = b$, then $(a,b)$ is empty, otherwise it is
an open segment.
Let $[a,b]$ be the union of $(a,b)$ and $\{ a, b \}$.
Then $[a,b]$ is always a closed segment.

One fundamental concept here is the following:

\begin{define}{(bridge)}
If $x$ is a point of a tree $T$ not contained in a 
connected set $B$, then there is a unique 
embedded path $[x,y]$ from $x$ to $B$.
This path has $(x,y) \cap B = \emptyset$ and
either  $y$ is in $B$ or $y$ is an accumulation
point of $B$. We say that $[x,y]$ is the bridge
from $x$ to $B$ and if $y$ is in $B$ we say that
$x$ bridges to $B$ in $y$ or that $x$ bridges
to $y$ in $B$.
\end{define}

For example if $T$ is the reals and $B = (0,1)$, $x = 2$, then
the bridge from $x$ to $B$ is $[2,1]$. 
One common use of bridges will be: if $x$ is not 
in a properly embedded line $l$ (as an axis defined below)
let $[x,y]$ be the bridge
from $x$ to $l$.
The crucial property of the bridge is that given $x$ and
$B$, the bridge is {\underline {unique}}.
In various situations this will force some useful 
equalities of points.
Another fundamental concept is:

\begin{define}{(axis)}
Suppose that $g$ is a homeomorphism acting freely
on a tree $T$. Then $g$ has an axis $\aaa_g$, a properly
embedded line in $T$, invariant under $g$ and $g$ acts
by translations on $\aaa_g$.
\end{define}

This is classical. Here $y$ is in $\aaa_g$
if and only if $yg$ is in $(y,yg^2)$, that is
$yg$ separates $y$ from $yg^2$.
Then it is easy to see that the axis must be
the union of $[y g^i, y g^{i+1}]$ where $i \in {\bf Z}$
\cite{Ba1,Fe}.
To obtain an element in $\aaa_g$ consider any $x \in T$.
If $x g \in (x, x g^2)$ done. Else 
there is a unique

$$y \ \in \ \ [x,xg] \cap [x, x g^2] \cap [x g, x g^2].$$

\noindent
$y$ is the basis of the tripod with corners $x, xg, xg^2$
\cite{Gr,Gh-Ha}.
A simple analysis of cases using free action yields
$y$ is in the axis.

Another simple but fundamental concept for us is:

\begin{define}{(local axis)}
Suppose $l$ is a line in  a tree $T$ where a homeomorphism
$g$ acts by translation. Then $l$ is a local axis for $g$
and is denoted by $\la_g$. The local axis may not be unique,
the context specifies which one we refer to.
\end{define}

For example if $g$ acts in $\rrrr$ by 
$x g = 2x$, then $\rrrr_+, \rrrr_-$ are
both local axes of $g$ with accumulation
point $x = 0$.
Another characterization of local axis: $x$ is in a 
local axis of $g$ if and only if $x g$ separates
$x$ from $x g^2$ (same definition as for axis
except requiring that $g$ acts freely).
Another characterization: suppose $x g$ is not $x$ 
and let $\uu$ be  the component of
$T - \{ x \}$ containing $x g$. Then $x$
is in a local axis of $g$ if and only
if $\uu g \subset \uu$.

Let $x$ be a point in a tree $T$. A {\em prong} at
$x$ is a non degenerate segment $I$ of $T$ so that $x$ is
one of the endpoints of $I$.
Two prongs at $x$ are equivalent if they
share a subprong at $x$. Associated to a subprong $I$
at $x$ there is a unique component $\uu$ of $T - \{ x \}$
containing $I - \{ x \}$.

\vskip .1in
\noindent
{\bf {Notation}} $-$ If $x, y, z$ are elements
in  a tree we will write $x \prec y \prec z$ if
$y$ separates $x$ from $z$,
or $y$ is in $(x,z)$.
We say that $x, y, z$ (in this order)
are {\em aligned}.
Also $x \prec y \preceq z$ if one also allows
$y = z$ and so on.
Notice that this
is invariant under homeomorphims of the tree.

The following simple results will be very useful:

\begin{lemma}{}{}
Let $x$ be a point in a tree $T$. Then two
prongs $I_1, I_2$ at $x$ are equivalent
if and only if the associated complementary
components $\uu_1, \uu_2$ are the same.
\label{inva}
\end{lemma}

\begin{proof}{}
If $I_1, I_2$ are equivalent, there is $y$ in $I_1 - \{ x \}$
also in $I_2$. Then clearly $y \in \uu_1$ and $y \in \uu_2$,
so $\uu_1 = \uu_2$.
Conversely suppose $\uu_1 = \uu_2$.
If $I_1$ is not equivalent to $I_2$, then $I_1 \cap I_2 = \{ x \}$
because $T$ is a tree and it also follows that $x$ separates
$I_1$ from $I_2$. This would imply $\uu_1$, $\uu_2$ disjoint,
contradiction.
\end{proof}

\begin{lemma}{}{}
Let $T$ be a tree and $\eta$ a homeomorphism so that
there are two points $x, y$ of $T$ so that
$x \prec x \eta \prec y \prec y \eta$ \ \ 
or \ \ $x \prec y \prec x \eta \prec y \eta$.
Then $x$ and $y$ are in a local axis of $\eta$.
\label{chax}
\end{lemma}

\begin{proof}{}
We do the proof for the first situation, the other
being very similar. Let $\uu$ be the component of $T - \{ 
x \}$ containing $x \eta$. Using $x \prec x \eta \prec y$
this is also the component of $T - \{ x \}$ containing
$y$. Apply $\eta$, then $\uu$
is taken to the component
of $T - \{ x \eta \}$ containing $y \eta$. Then
$\uu \eta$ is contained in $\uu$ and $x$ is in
a local axis. Apply $\eta^{-1}$ to $y$ to get 
$y$ is in a local axis as well.
We stress the two local axes produced in
this way a priori may not be the same: there may
be a fixed point of $\eta$ in $(x,y)$.
\end{proof}

\vskip .1in
\noindent
{\bf Global fixed points}

Here we consider the case that an essential
lamination $\lambda$ on $N$ would produce
a trivial group action on a tree $T$. 

Recall the notion of {\em efficient} transversal
to a lamination: let $\eta$ be a transversal
to a lamination $\lambda$.
Then $\eta$ is efficient \cite{Ga-Oe} if
for any subarc $\eta_0$ with both endpoints 
in leaves of $\lambda$ and interior disjoint
from $\lambda$, then $\eta_0$ is not homotopic
rel endpoints into a leaf of $\lambda$.
Gabai and Oertel showed  that if $\lambda$ is
essential then any efficient transversal
cannot be homotoped rel endpoints into a leaf
of $\lambda$. Also closed efficient transversals
are not null homotopic.

\begin{lemma}{}{}
If $\lambda$ is an essential lamination in
$N$ then the associated group action of
$\pi_1(N)$ on a tree $T$ as described above
has no global fixed point.
\end{lemma}

\begin{proof}{}
Suppose on the contrary that a point $x$ of
$T$ is left invariant by the whole group.
Look at the preimage of $x$ in the possibly
non Hausdorff tree $T_*$.
There are 3 options:

$1 -$ $x$ comes from a non singular, Hausdorff leaf
$E$ of $\wl$. Then $E$ is left invariant by the whole
group $\pi_1(N)$,

$2 -$   $x$ comes
from the closure $R$ of a complementary region 
of $\wl$ in the universal
cover. Then $R$ is left invariant by the whole group.
In this case let $E$ be a boundary leaf of $R$.

$3 - $ Finally $x$ may come from a non Hausdorff 
leaf $E$.  
Then the orbit of $E$ under $\pi_1(N)$
consists only of the non separated leaves from $E$.

By construction of the tree $T$ above
these 3 cases are mutually 
exclusive.
It follows that in any of the 3 options there is
at least one component $B$ of $\widetilde N - E$ which does
not contain any translate of $E$. In option 1)
any component will do, in option 2) choose the
component not containing $R - E$ and in option
3) choose the component not containing leaves
non separated from $E$. 

Let $A = \pi(E)$ where $\pi: \widetilde N \rightarrow N$ is
the universal covering map. 
Suppose first that $A$ is not compact.
Then it limits on some leaves of $\lambda$ and
there is a laminated box where $A$ intersects it
in at least $3$ leaves and the box intersects
an efficient transversal to $\lambda$.
Lifting to $\widetilde N$ so that the middle leaf
is $E$ then the other 2 leaves are not $E$ (efficient
transversal) and one of them is contained in $B$
producing a covering translate of $E$ in $B$, contradiction.
The same is of course true if $A$ intersects
an efficient closed transversal.

Now $A$ is compact. If $A$ is non separating, then
it intersects a closed transversal associated to $g$ in
$\pi_1(N)$ only once. Same proof yields either
$Eg$ or $Eg^{-1}$ in $B$, done.

Finally suppose that $A$ is separating. Then $C = \pi(B \cup E)$
is a compact submanifold of $N$ which has $A$ as its
unique boundary component.
For any $g$ in $\pi_1(C)$ then $Eg$ is contained in $B \cup E$,
so by hypothesis must be $E$, therefore $\pi_1(A)$ surjects
in $\pi_1(C)$. As $\lambda$ is essential then $\pi_1(A)$
also injects \cite{Ga-Oe}, so $\pi_1(A)$ is isomorphic
to $\pi_1(C)$. As $C$ is irreducible \cite{Ga-Oe}, then
theorem 10.5 of Hempel \cite{He} implies that $C$ is
homemorphic to $A \times I$ with $A$ corresponding
to $A \times \{ 0 \}$. This contradicts the
fact that $A$ is the only boundary component
of $C$. 
This finishes the proof of the lemma.
\end{proof}

\vskip .1in
\noindent
{\bf {Remark:}} $-$ Notice that leaves of essential
laminations may not intersect a closed transversal.
For example this occurs for separating incompressible
surfaces. It also occurs for leaves
of Reebless foliations which have a separating
leaf (which necessarily must be a torus or Klein bottle)
$-$ there are many examples of these.
So Reebless foliations which are also essential
laminations need not be taut foliations!

\section{Case R $-$ the $\rrrr$-covered case}

For the remainder of the article we consider
the manifold $M_{p/q}$ as described in
section \ref{grou} with fundamental group
$\gl$. The goal is to show it does not
admit an essential lamination. Suppose then on
the contrary that there is an essential
lamination $\lambda$ on $M_{p/q}$.
Let $T$ be the associated tree with non trivial
action of $\gl$ on it.
Notice that since $\alpha, \tau$ generate $\gl$ then
no point of $T$ is fixed by both $\alpha$ and $\tau$.

The conditions on the parameters are $|p - 2q| = 1$ and
$m \leq -4$. They will not be used in full force
for all the arguments. Many times all we need
is $p \geq q$ or $p$ odd or $m$ negative or none
of these. 
The proof is done by subdiving into subcases and
showing each subcase is impossible leading to 
various contradictions.

In this section we assume that $T$ is homeomorphic to
the real numbers and study non trivial actions of $\gl$
in $\rrrr$.
Notice that $\gamma$ being a commutator is an orientation
preserving homeomorphism of $\rrrr$. 
Since $\tau^p \gamma^q = id$, then $\tau^p$ is also
orientation preserving. 

We use the relations from the group presentation
of $\gl$ or variations thereof.

Suppose first the action is orientation preserving on $\rrrr$:

\vskip .1in
\noindent
{\bf {Case R.1}} $-$ $\alpha$, $\tau$ are orientation
preserving.

As $\beta = \tau \alpha \tau^{-1}$ then $\beta$ also is
orientation preserving and so is the whole group $\gl$.
We subdivide into subcases:

\vskip .1in
\noindent
{\bf {Case R.1.1}} $-$ $\tau$ has a fixed point $x$.

Then $x \alpha$ is not $x$. Orient $\rrrr$ so that 
$x \alpha > x$.
As $\gamma$ is orientation preserving then $x \gamma = x$.
Then  applying $\gamma \tau \beta \alpha^m = \alpha \tau$
to $x$:
 
$$x \gamma \tau \beta \alpha^m = x \alpha \tau > x \tau = x$$

\noindent
which uses $\tau$ orientation preserving. Hence $x \beta \alpha^m > x$ \
or \ $x \beta > x \alpha^{-m} > x$ \ (as $-m > 0$).
\ Hence $x \beta^{-1} < x$.
But also 

$$x \beta^{-1} \ = \ x \tau \alpha \tau^{-1} \ = \ 
x \alpha \tau^{-1} \ > \ x \tau^{-1} \ = \ x.$$

\noindent
This is a contradiction, ruling out this case.

\vskip .1in
\noindent
{\bf {Case R.1.2}} $-$ $\tau$ acts freely, $\alpha$ has
a fixed point $x$.

Assume $\tau$ is increasing in $\rrrr$.
As $\tau = \kappa^q$ and $q$ is positive
then $\kappa$ is increasing.
Here use $x \alpha \tau = x \tau = x \gamma \tau \beta \alpha^m$.
Hence $x \tau \alpha^{-m} = x \gamma \tau \beta$.
As $x \tau > x$ then $x \tau \alpha^{-m} > x$.
Hence $x \gamma \tau > x \beta^{-1}$.
Here $\gamma = \kappa^{-p}$ and $\gamma \tau = \kappa^{q - p}$.
As $q \leq p$ then $q - p \leq 0$ and $\gamma \tau$ is
monotone decreasing or constant.
Hence

$$x \beta^{-1} \ < x \gamma \tau \ \leq \ x.$$

\noindent
One fact that will be used in a lot of arguments is that
under the condition $p \geq q$ when $\gamma, \tau$ act
freely and $x \tau > x$ then
$x \gamma \leq x \tau^{-1}$.
Notice that $x \tau^{-1} \beta = x \alpha \tau^{-1} = x \tau^{-1}$.
On the other hand

$$x \beta \ = \ x \alpha \beta \ = \ 
x \gamma \beta \alpha \ \leq \ x \tau^{-1} \beta \alpha
\ = \ x \tau^{-1} \alpha \ < \ x \alpha \ = \ x.$$

\noindent
leading to the contradiction that both $x \beta$ and
$x \beta^{-1}$ are $< x$.

Notice a lot of these arguments are using orientation
preserving homeomorphims.

\vskip .1in
\noindent
{\bf {Case R.1.3}} $-$ $\tau$ acts freely increasing in $\rrrr$
and $\alpha$ acts freely, also increasing in $\rrrr$.

Take any $x$ in $\rrrr$. Then $x \alpha \tau > x$
so $x \gamma \tau \beta \alpha^m > x$.
So $x \gamma \tau \beta > x \alpha^{-m} > x$.
Since $x \gamma \tau \leq x$ this implies 
$x \beta > x$.
On the other hand,

$$x \beta \ = \ x \tau \alpha^{-1} \tau^{-1} \ <
\ x \tau \tau^{-1} \ =  \ x,$$

\noindent
contradiction.

\vskip .1in
\noindent
{\bf {Case R.1.4}} $-$ $\tau$ acts freely and increasing
in $\rrrr$, $\alpha$ acts freely and decreasing in $\rrrr$.

This implies $z \alpha^{-1} > z$ for all $z$ in $\rrrr$.
For any $x$ in $\rrrr$,
$x \beta = x \tau \alpha^{-1} \tau^{-1} > 
x \tau \tau^{-1} = x$.
Also $x \tau^{-1} \alpha \tau < x$ for all $x$. 
Hence

$$x \alpha \beta \alpha^{m-1} \ = \ x \tau^{-1} \alpha \tau 
\ < \ x,$$

\noindent
for all $x$. Hence $x \alpha \beta < x \alpha \alpha^{-m}
< x \alpha$ for all $x$ \ ($-m > 0$). But this contradicts
$(x \alpha) \beta > x \alpha$ because $\beta$ is increasing
everywhere as proved above.

This finishes the analysis of $T$ homeomorphic to
$\rrrr$ and orientation preserving action.

%
%
%
%
%
%
%

\vskip .2in
We now deal with orientation reversing cases.
The general case of $\tau$ orientation reversing is
hard, so we use one of the hypothesis to discard it
as follows: $\tau^p = \gamma^{-q}$ is orientation
preserving as $\gamma$ always is.
We are mainly interested in $|p-2q| = 1$, which
implies $p$ odd and if $p$ is odd and $\tau^p$ orientation
preserving then $\tau$ is also orientation preserving.
We now deal with the case $\alpha$ orientation reversing.

\vskip .1in
\noindent
{\bf {Case R.2}} $-$ $\alpha$ orientation reversing,
$\tau$ orientation preserving.

Let $x$ be the unique fixed point of $\alpha$. As
$x \tau \not = x$, assume $x \tau > x$.
As $\tau = \kappa^q$ and $q > 0$, this implies $k$ is
increasing in $x$.
Notice that $x \tau^{-1}$ is the unique fixed point
of $\beta$.
The subcases depend on the relative position
of $x \tau \alpha$ and $x \tau^{-1}$. Notice that $x \tau > x$,
so $x \tau \alpha < x \alpha = x$.

\vskip .1in
\noindent
{\bf Case R.2.1} $-$ $x \tau \alpha < x \tau^{-1}$

Then $x \tau \alpha \tau^{-1} = x \beta^{-1} < x \tau^{-2}$.
Notice 

$$x \tau \gamma \beta \alpha^m = x \alpha \tau = x \tau > x$$

\noindent
so $x \tau \gamma \beta > x \alpha^{-m} = x$ and so

$$x \tau \gamma < x \beta^{-1} < x \tau^{-2}$$

\noindent
or $x \tau^3 \gamma < x$. 
As $\tau^3 = \kappa^{3q}$ and $\gamma = \kappa^{-q}$, then
$x \kappa^{3q - p} < x$. As $\kappa$ is increasing
in $x$ then $3q - p < 0$ or $p > 3q$.
Arguments such as this will be used in various
parts of the proof.
Since in the end we want $p = 2q \pm 1$ we can discard
this case.

\vskip .1in
\noindent
{\bf {Remark}} $-$ What we really wanted was to rule out
this case without using $p = 2q \pm 1$, but we were
unable to do that. Our partial results (without
using $p = 2q \pm 1$) show that
$x \tau \alpha^3 > x \tau \alpha$
so $x < x \tau \alpha^2 < x \tau$.
Also there is a fixed point of $\alpha^2$ between
$x \tau$ and $x \tau^2$ and 
$\alpha^2$ acts expandingly (away from $x$) in some point.
Something similar is also true in the following case.

\vskip .2in
\noindent
{\bf Case R.2.2} $-$ $x \tau \alpha > x \tau^{-1}$

%
%
%
%
%
%
%
%
%
%
%
%

First notice that $x \beta^{-1} < x \tau^{-1}$.
Use

$$(x \tau) \tau \gamma \beta \alpha^m = (x \tau) \alpha \tau > 
x \tau^{-1} \tau = x$$

\noindent
so $x \tau^2 \gamma \beta > x \alpha^{-m} = x$ and 

$$ x \tau^2 \gamma < x \beta^{-1} < x \tau^{-1}.$$

\noindent
We conclude as in the previous case that
$x \tau^3 \gamma < x$ or $p > 3q$, also disallowed.

The reader may think we just got lucky to get $p > 3q$ as 
we have the hypothesis 
$p = 2q \pm 1$. The remaining case explains why this
has happened.

\vskip .2in
\noindent
{\bf {Case R.2.3}} $-$ $x \tau \alpha = x \tau^{-1}$.

This case is much more interesting.
First

$$x \alpha \tau = x \tau \alpha \beta \alpha^{m-1}$$

\noindent
Since $x \tau \alpha = x \tau^{-1}$ this is left
invariant by $\beta$, so the right side
is $x \tau \alpha \alpha^{m-1} = x \tau \alpha^m$ 
equal to $x \tau$.
Since $m$ is even,  $\alpha^m$ preserves orientation,
therefore $x \tau \alpha^2 = x \tau$.
Also $x \tau \alpha = x \tau \alpha^{-1} = x \tau^{-1}$.
Now notice that 

$$x \tau \gamma \beta \alpha^m = x \alpha \tau
= x \tau, \ \ \ {\rm so} \ \ \ 
x \tau \gamma = x \tau \alpha^{-m} \beta^{-1},$$

\noindent
or $x \tau \gamma = x \tau \beta^{-1}$.
Now we show that $x \tau^2 \alpha = x \tau^{-2}$.
To show this use $x \beta^{-1} \tau = x \tau \alpha
= x \tau^{-1}$, hence $x \beta^{-1} = x \tau^{-2}$.
Use 

$$ \tau^{-2} \beta \tau^2 \ = \ \tau^{-1} \alpha^{-1} \tau
\ = \ \alpha^{1-m} \beta^{-1} \alpha^{-1}$$

\noindent
applied to $x$:

$$x \tau^{-2} \beta \tau^2 = x \alpha^{1-m} \beta^{-1} \alpha^{-1}$$

\noindent
or 
$x \beta^{-1} \beta \tau^2 = x \beta^{-1} \alpha^{-1}$ so

$$x \tau^2 = x \tau^{-2} \alpha^{-1}.$$

\noindent
Then

$$x \tau^{-2} = x \tau^2 \alpha =  (x \tau) \tau \alpha
= x \tau \beta^{-1} \tau
= x \tau \gamma \tau$$

\noindent
or 

$$x  \gamma \tau^4 = x$$

\noindent
As seen before this implies $p = 4q$ or $p = 4$,
$q = 1$. 
This is disallowed by $p$ being odd.

We remark that in this case the group in fact acts non trivially
in $\rrrr$. For instance let

$$x \alpha = -x, \ \ \ x \tau = x + 1$$

\noindent
It is easy to check they satisfy the equations if $m$ is
even!

It may be true that this is the only possibility and when
$x \tau \alpha \not = x \tau^{-1}$ we get a perturbation
of this, namely that $p$ is close to $4q$ and in fact
$p > 3q$.

\section{Case A $-$ $\tau$ acts freely}

In this section we consider the case that $\tau$ acts
freely in $T$.
This implies that $\kappa^q$ acts freely in the tree,
and therefore $\kappa$ itself acts freely.
In addition the axes are the same $\aaa_{\kappa} = \at$.
Here we will use the relation $\alpha \beta = \gamma \beta \alpha$
in the following form, defining an element
$u$ of $\gl$:

$$u \ = \ \alpha \tau \alpha^{-1} \tau^{-1}
\ = \ \gamma \tau \alpha^{-1} \tau \alpha$$

\noindent
We will consider the intersections $\ak \cap \ak \alpha$ and
$\ak \cap \ak u$.
The axis $\ak$ is homeomorphic to the real numbers.
Put an order $<$ in $\ak$ so that 
$x <  x \tau$ for any $x$ in $\ak$.
This induces an order $\pa$ in $\ak \alpha$
so that $x < y$ in $\ak$ if and only if \
$x  \alpha \ \pa \ y \alpha$ \ in $\ak \alpha$ and similarly
put order $\pu$  in $\ak u$ so that $x < y$
in $\ak$ if and only if \ $x u \ \pu \ y u$ \ in $\ak u$.

\vskip .15in
\noindent
{\bf {Case A.1}} $-$ $\ak \alpha \cap \ak$ has at
most one point.

If the intersection is a single point $x$,
let $y = x$ as well.

If they are disjoint, there is a single point
$x$ in $\ak$ bridging to $\ak \alpha$. 
For intance $x$ is the unique point so that
there is a path from $x$ to $\ak \alpha$ intersecting
$\ak$ only in $x$. Another way to characterize $x$,
it is the only point so that $x$ separates the
rest of $\ak$ from $\ak \alpha$. In other words
the components of $T - \{ x \}$ containing $\ak
\alpha$ and the rest of $\ak$ are all disjoint.
In the same way there is a single $y$ in $\ak \alpha$ which
is the closest to $\ak$. Then $[x,y]$ is a path
from $\ak$ to $\ak \alpha$ so that $(x,y)$ does 
not intersect either $\ak$ or $\ak \alpha$ $-$
this is an equivalent way to get the segment $[x,y]$.
This path $[x,y]$ is called the {\em bridge}
from $\ak$ to $\ak \alpha$. This extended notion
of bridges will also be used in the article.
It is invariant by homemorphisms of the tree.
The bridge between connected sets is also unique.

We now use the relation above. 
The proof is very similar to ping pong lemma
arguments. Since $\ak$ is invariant
under $\gamma$ and $\tau$, the right side
says that $\ak u = \ak \alpha^{-1} \tau^{-1} \alpha$.


\begin{figure}
\centeredepsfbox{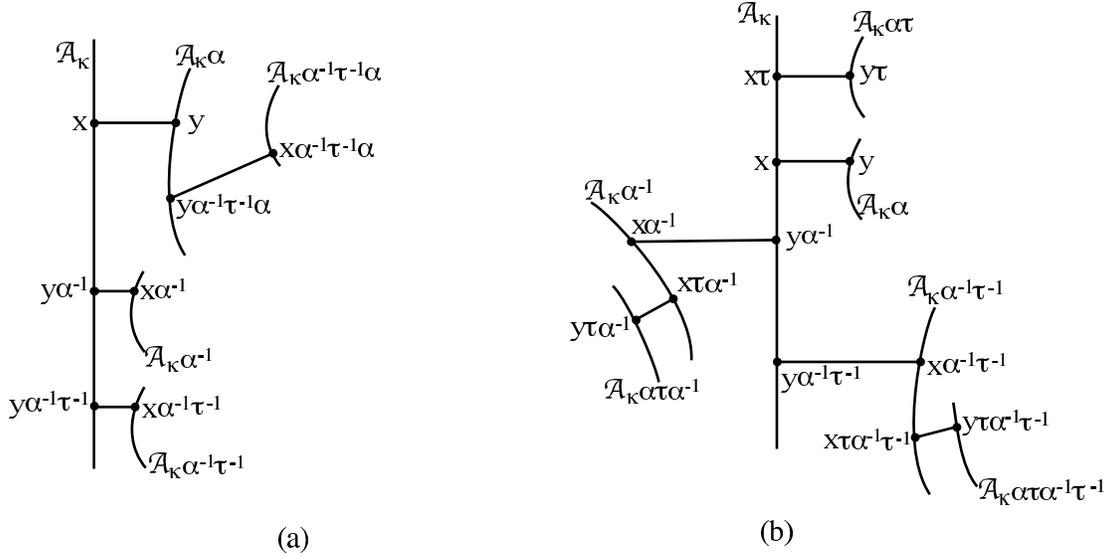}
\caption{
The case $\ak \cap \ak \alpha = \emptyset$.
The same arguments can be used for intersection
a single point. a. Using $\ak u = \ak \alpha^{-1} \tau^{-1} \alpha$,
b. Using $\ak u = \ak \alpha \tau \alpha^{-1} \tau^{-1}$.}
\label{02}
\end{figure}

The bridge from $\ak$ to $\ak \alpha$ is $[x,y]$ -
degenerate $[x,x]$ when they intersect in a point.
Therefore the bridge from $\ak \alpha^{-1}$ to
$\ak$ \ is $[x \alpha^{-1}, y \alpha^{-1}]$, see fig. 
\ref{02}, a.
Then the bridge $\ak \alpha^{-1} \tau^{-1}$ to
$\ak$ is $[x \alpha^{-1} \tau^{-1},
y \alpha^{-1} \tau^{-1}]$.
This implies that 
the bridge from

$$\ak \alpha^{-1} \tau^{-1} \alpha \ \ \ {\rm  to} \ \ \
\ak \alpha \ \ \ {\rm  is } \ \ \
[x \alpha^{-1} \tau^{-1} \alpha, \  y \alpha^{-1} \tau^{-1} 
\alpha].$$

\noindent
Notice that $y \alpha^{-1} \tau^{-1}$ is not $y \alpha^{-1}$.
Therefore $y \alpha^{-1} \tau^{-1} \alpha$ is not $y$. 
It now follows that

$${\rm the \ bridge \ from} \ \ \ 
\ak u = \ak \alpha^{-1} \tau^{-1} \alpha \ \ \ {\rm  to}
\ \ \ \ak \ \ \ 
{\rm is} 
\ \ \ [x \alpha^{-1} \tau^{-1} \alpha, x].$$

\noindent

On the other hand use that $\ak u = \ak \alpha \tau
\alpha^{-1} \tau^{-1}$.
The bridge from $\ak \alpha \tau$ to $\ak$ is
$[y \tau, x \tau]$, see fig. \ref{02}, b.
The bridge from $\ak \alpha \tau \alpha^{-1}$
to $\ak \alpha^{-1}$ is 
$[y \tau \alpha^{-1}, x \tau \alpha^{-1}]$ and
the bridge from $\ak \alpha^{-1}$ to $\ak$ is
$[x \alpha^{-1}, y \alpha^{-1}]$.
Since $x \alpha^{-1}$ is not equal 
$x \tau \alpha^{-1}$ then the bridge from
$\ak \alpha \tau \alpha^{-1}$ to $\ak$ is
$[y \tau \alpha^{-1}, y \alpha^{-1}]$. 
Finally 

$${\rm the \ bridge \ from} \ \ \ \ak u
\ \ {\rm  to} \  \ \ak  \ \ {\rm is} \ \ \
[y \tau \alpha^{-1} \tau^{-1}, \ y \alpha^{-1} \tau^{-1}].$$

\noindent
Since the bridge from $\ak u$ to $\ak$ is uniquely defined
this implies

$$ y \alpha^{-1} \tau^{-1} \ = \ x, 
\ \ \ \ y \tau \alpha^{-1} \tau^{-1} \ = \ 
 x \alpha^{-1} \tau^{-1} \alpha.$$

\noindent
So  $y = x \tau \alpha$ and

$$x \alpha^{-1} \tau^{-1} \alpha \ = \ 
x \tau \alpha \tau \alpha^{-1} \tau^{-1},
\ \ \ \ {\rm or} \ \ \ \
x \alpha^{-1} \tau^{-1} \alpha \tau \alpha \ = \ 
x \tau \alpha \tau.$$

\noindent
Use $\tau^{-1} \alpha \tau = \alpha \beta \alpha^{m-1}$,
so 

$$\alpha^{-1} \tau^{-1} \alpha \tau \alpha \ = \
\alpha^{-1} \alpha \beta \alpha^{m-1} \alpha \ = \ 
\beta \alpha^m \ = \ 
\gamma^{-1} \tau^{-1} \alpha \tau,$$

\noindent
so $x \gamma^{-1} \tau^{-1} \alpha \tau = x \tau \alpha \tau$,
or $x \gamma^{-1} \tau^{-1} = x \tau$. This implies
$x \gamma \tau^2 = x$ and as seen before implies
$p = 2q$. This is disallowed by $p$ odd.

\vskip .1in
We now consider intersections with more than one point.

\vskip .15in
\noindent
{\bf {Case A.2}} $-$ $\ak \cap \ak \alpha \ = \ [x,y]$.

Here $x$ is not equal to $y$ and $x < y$ in $\ak$.
We include some ideal point cases: $x$ could $-\infty$
and $y$ could be $+\infty$, in which case
the intersection is a ray in $\ak$.
On the other hand we can never have $\ak = \ak \alpha$.
Otherwise $\alpha, \tau$ leave $\ak$ invariant, so
the whole group does. But $\ak$ is homemorphic to
$\rrrr$ $-$ this was disallowed by no actions
on $\rrrr$.

Since the intersection is a non trivial interval one
considers separately whether the orders $<$, $\pa$ agree
on the intersection.

\vskip .15in
\noindent
{\bf {Case A.2.1}} $-$ The orders $<$ and $\pa$ agree on
$\ak \cap \ak \alpha$.

It is easy to check that this is equivalent to
\ $x \alpha^{-1} < \ y \alpha^{-1}$ \ in $\ak$,
by applying $\alpha$ to the pair $x \alpha^{-1},
y \alpha^{-1}$ both of which are in $\ak$.

We now consider $\ak u$. We first use $\ak u = \ak \alpha^{-1}
\tau^{-1} \alpha$. Notice that

$$\ak \cap \ak \alpha^{-1} \ = \ 
[x \alpha^{-1}, y \alpha^{-1}] \ \ \ \ \
{\rm so}
\ \ \ \ \ 
\ak \alpha^{-1} \tau^{-1} \cap \ak \ = \
[x \alpha^{-1} \tau^{-1}, \ y \alpha^{-1} \tau^{-1}],$$

\noindent
in the correct order. Hence

$$\ak u \cap \ak \ = \
[x \alpha^{-1} \tau^{-1} \alpha, \
y \alpha^{-1} \tau^{-1} \alpha].$$

\noindent
In addition \
$x \alpha^{-1} \tau^{-1} \alpha \ \pa \ 
y \alpha^{-1} \tau^{-1} \alpha$.


\begin{figure}
\centeredepsfbox{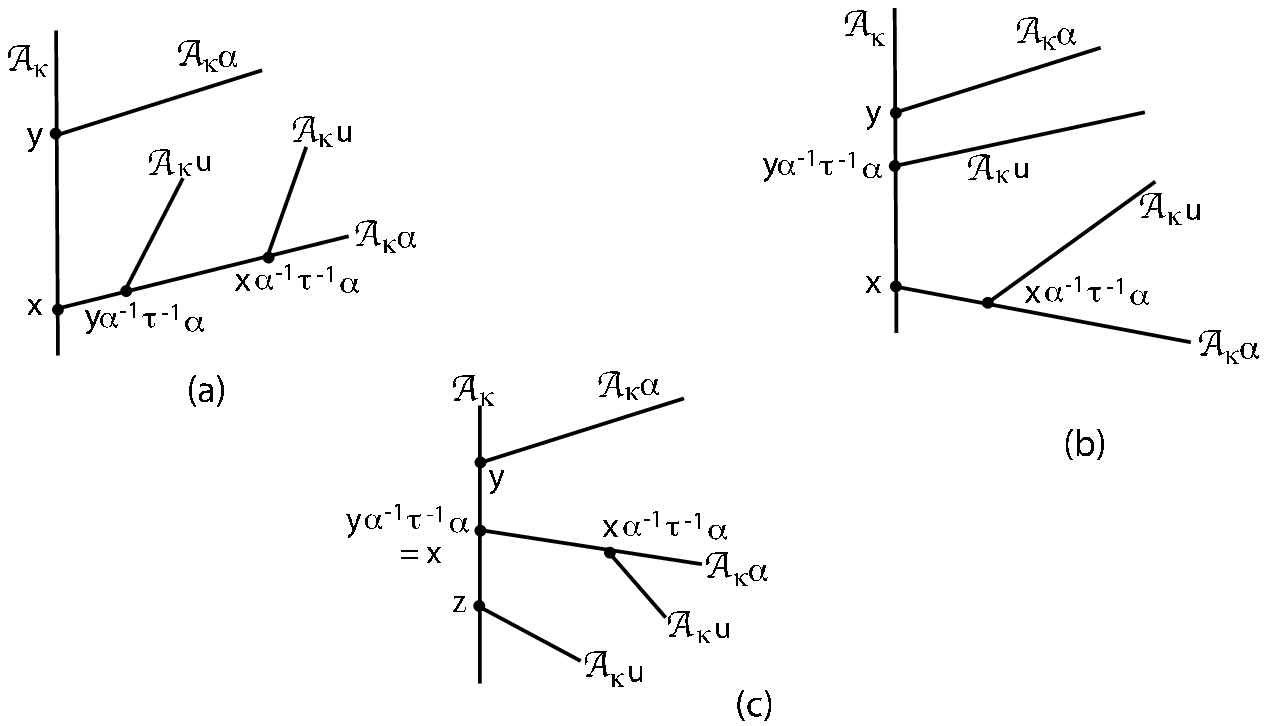}
\caption{
Evaluating $\ak u \cap \ak$,
using $\ak u = \ak \alpha^{-1} \tau^{-1} \alpha$,
\ a. $y \alpha^{-1} \tau^{-1} \alpha \pa x$,
\ b. $y \alpha^{-1} \tau^{-1} \alpha >_{\alpha}  x$,
\ c. $y \alpha^{-1} \tau^{-1} \alpha =  x$.}
\label{03}
\end{figure}

Notice that $x \alpha^{-1} \tau^{-1} < x \alpha^{-1}$
in $\ak$, hence \ $x \alpha^{-1} \tau^{-1} \alpha \
\pa \ x$ \ in $\ak \alpha$.
Also \ $y \alpha^{-1} \tau^{-1} \alpha \ \pa \ y$ 
\ in $\ak \alpha$.
Given this there are 3 options:


1) \ If \ $y \alpha^{-1} \tau^{-1} \alpha \ \pa \ x$ \ in
$\ak \alpha$ then 
$\ak u \cap \ak = \emptyset$ and the
bridge from $\ak$ to $\ak u$ is $[x, y \alpha^{-1}
\tau^{-1} \alpha]$, fig. \ref{03}, a.

2) \ If \ $y \alpha^{-1} \tau^{-1} \alpha \ >_{\alpha} \ x$ \ in
$\ak \alpha$ then $y \alpha^{-1} \tau^{-1} \alpha$
is in $(x,y)$ and $\ak u \cap \ak = [x, y \alpha^{-1}
\tau^{-1} \alpha]$. In addition the orders
$<$ and $\pu$ agree on $\ak \cap \ak \alpha$,
see fig. \ref{03}, b.

3) \ If $y \alpha^{-1} \tau^{-1} \alpha = x$, then
$\ak \alpha  \cap \ak = [z,x]$. In addition if $z$ is not
$x$ then 
the orders $<$ and $\pu$ disagree
on $\ak \cap \ak u$, see fig. \ref{03}, c. In this case both
$x$ and $y$ are finite.
The last option can occur because $\ak u$ can enter
$\ak$ in $x$ but rather than going up, going in
the opposite direction $-$ the one containing
$x \tau^{-1}$.


\vskip .1in

Notice that the 3 options are mutually exclusive.
We now consider $\ak u = \ak \alpha \tau \alpha^{-1} \tau^{-1}$.
Use 

$$\ak u \cap \ak \ = \ 
(\ak \alpha \tau \cap \ak \alpha) \alpha^{-1}
\tau^{-1}.$$

\noindent
Here $\ak \alpha \tau \cap \ak = [x \tau, y \tau]$.
So whether $\ak \alpha \tau \alpha^{-1}$
and $\ak$ intersect, depends on the relative positions
of $x \tau$ and $y$.
Notice that $x \tau > x$ in $\ak$.


\begin{figure}
\centeredepsfbox{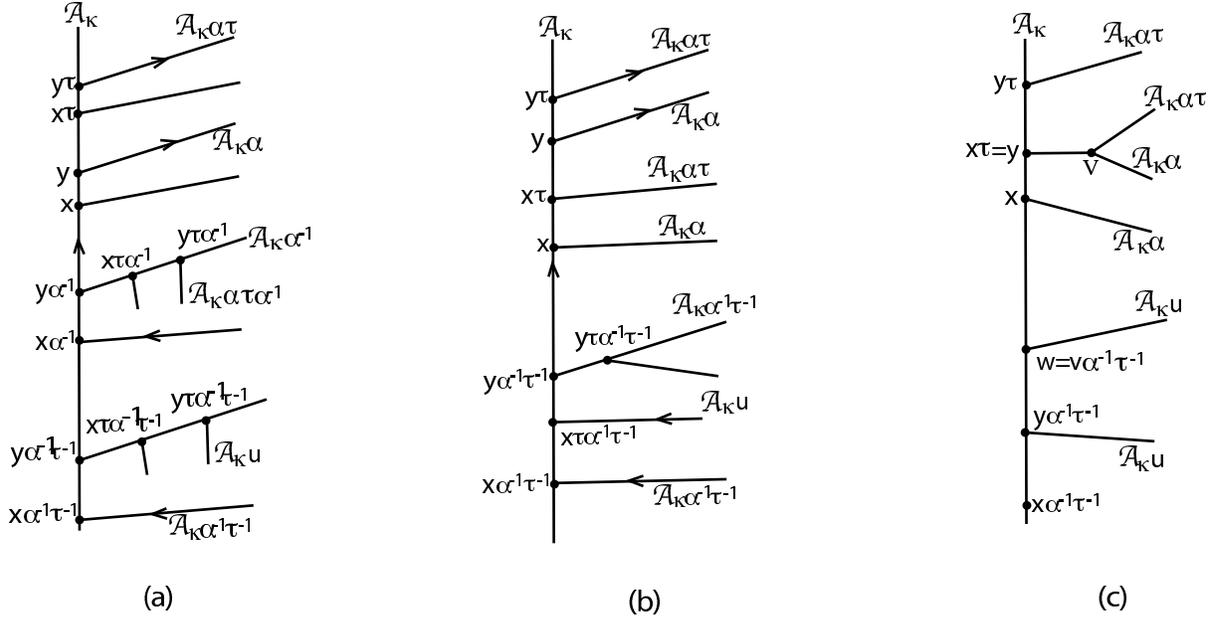}
\caption{
Using $\ak u = \ak \alpha \tau \alpha^{-1}
\tau^{-1}$,
\ a. $x \tau > y$, \ b. $x \tau < y$, \ c. $x \tau = y$.}
\label{04}
\end{figure}


1')  \ If $x \tau > y$ in $\ak$ then \ $\ak \alpha \tau \cap 
\ak  \alpha = \emptyset$, so
\ $\ak \alpha \tau \alpha^{-1} \cap \ak = \emptyset$.
Therefore $\ak u \cap \ak = \emptyset$ and
the bridge from $\ak$ to $\ak u$ is
\ $[y \alpha^{-1} \tau^{-1}, \  x \tau \alpha^{-1} \tau^{-1}]$,
see fig. \ref{04}, a.
Here $x, y$ finite.

2') \  If $x \tau < y$ in $\ak$ then 
$\ak \alpha \tau \cap \ak \alpha =
[x \tau, y]$, then $\ak \cap \ak u$ is
\ $[x \tau \alpha^{-1} \tau^{-1},
y \alpha^{-1} \tau^{-1}]$
(the first term smaller in $\ak$),
and the orders $<$ and $\pu$ agree on
$\ak \cap \ak u$, see fig. \ref{04}, b.

3') \ If $x \tau = y$, then  $\ak \alpha \tau \cap \ak \alpha
= [y, v]$.
Notice we may have $v \not = y$. So 
$\ak u \cap \ak = [y \alpha^{-1} \tau^{-1}, w]$,
where $w = v \alpha^{-1} \tau^{-1}$.
Here $x$ and $y$ are finite and if
$w$ is not equal to $x \tau \alpha^{-1} \tau^{-1}$, then
the orders $<$ and $\pu$ disagree on $\ak \cap \ak u$.
Notice that order in $\ak \alpha \tau$ goes from
$v$ to $y$, so the increasing order
$\pu$ in $\ak u$
from $w = v \alpha^{-1} \tau^{-1}$
to $y \alpha^{-1} \tau^{-1}$, see fig. \ref{04}, c.


\vskip .1in
Notice that again all 3 cases are mutually exclusive.
Therefore we can match the 2 pairs of 3 possibilities to
get 3 mutually exclusive cases:

\vskip .1in

{\bf I} $-$ $y \alpha^{-1} \tau^{-1} \alpha \ \pa \ x$ \ in
$\ak \alpha$ or $x \tau > y$ in $\ak$ and
$\ak \cap \ak u = \emptyset$. In this case

$$[x, \  y \alpha^{-1} \tau^{-1} \alpha] \ = \
[y \alpha^{-1} \tau^{-1}, \ x \tau \alpha^{-1} \tau^{-1}]$$

{\bf II} $-$ $y \alpha^{-1} \tau^{-1} \alpha \ \poa \ x$ \ in
$\ak \alpha$ or $x \tau < y$ in $\ak$ and

$$\ak \cap \ak u \ = \ 
[x, \ y \alpha^{-1} \tau^{-1} \alpha] \ = \ 
[x \tau \alpha^{-1} \tau^{-1}, \ y \alpha^{-1} \tau^{-1}]$$

{\bf III} $-$ $y \alpha^{-1} \tau^{-1} \alpha = x$
or $x \tau = y$. 
Then 

$$\ak \cap \ak u \ = \ 
[z, x] \ = \ [y \alpha^{-1} \tau^{-1}, w]$$

\noindent
If $z$ is not $x$ then the orders
$<$ and $\pu$ disagree on $\ak \cap \ak u$.

\vskip .1in

We now deal with each situation separately.

\vskip .1in
\noindent
{\bf {Situation II}} $-$ 

Here $x \tau \alpha = x \tau$ and 
$x \tau$ is in $(x,y)$.
Let  $\uu_1$, (respectively $\uu_2$)
be the component of $T - \{ x \tau \}$ 
containing $y$ (respectively $x$).
Here $[x,y] = \ak \cap \ak \alpha$,
$x \tau$ is in the interior of $[x,y]$ 
and then the orders
$<$, $\pa$ agree on $[x,y]$. It
follows that the prongs $[x \tau, y]$,
$[x \tau, y \alpha]$ are equivalent. 
By lemma \ref{inva}, $\uu_1 \alpha = \uu_1$. 
In the same way
$\uu_2 \alpha = \uu_2$. 
This situation is disallowed by the following lemma.

\begin{lemma}{}{}
Suppose that ${\cal L}$ is a local axis for $\tau$
and $r$ is a point in ${\cal L}$ with
$r \alpha = r$.
Suppose that $\uu_1$ ($\uu_2$ respectively) is
the component of $T - \{ r \}$ containing $r \tau$
($r \tau^{-1}$ respectively). Then at least one
of $\uu_1$ or $\uu_2$ is not invariant under $\alpha$.
\end{lemma}

\begin{proof}{}
On the contrary suppose that $\uu_i \alpha = \uu_i$ for
$i = 1,2$.
We will arrive at a contradiction.
Let $\vv_i = \uu_i \tau^{-1}$. Then the conjugation of
$\beta$ with $\alpha^{-1}$ by $\tau$ implies that
$\vv_i \beta = \vv_i$, $i = 1,2$.
Use

$$r \tau^{-1} \alpha \tau \ = \ r \gamma \beta \alpha^m$$

\noindent
Since $p \geq q$, then $r \gamma \leq r \tau^{-1}$ in
${\cal L}$ (with $\tau$ increasing in ${\cal L}$
and so $r \gamma \beta$ is in $\vv_2 \cup \{ r \tau^{-1} \}$
contained  in $\uu_2$. Therefore $r \gamma \beta \alpha^m$
is in $\uu_2$.
Consequently $r \tau^{-1} \alpha \tau$ is in $\uu_2$ and
$r \tau^{-1} \alpha$ is in 
$\uu_2 \tau^{-1} = \vv_2$ \ \ \ (*).

On the other hand $r \gamma \in \vv_2 \cup \{ r \tau^{-1} \}$,
so 

$$r \beta \alpha^{-1} \ = \ 
r \gamma \beta \ \in \ \vv_2 \cup \{ r \tau^{-1} \},$$

\noindent
so $r \tau^{-1}$ is in $[r \beta \alpha^{-1},r)$. Apply 
$\alpha$ to obtain 

$$r \tau^{-1} \alpha \  \in \
[r \beta, r) \ \ \ \ \ (**).$$

\noindent
Now

$$r \beta \ = \ r \tau \alpha^{-1} \tau^{-1}
\ \ \ {\rm and} \ \ \ 
r \tau \in \uu_1 \ \Rightarrow  \  r \tau \alpha^{-1} \in \uu_1
\ \Rightarrow \ r \beta = r \tau \alpha^{-1} \tau^{-1}
\in \vv_1.$$

\noindent
As $r$ is also in $\vv_1$, it follows
from (**)
that $r \tau^{-1} \alpha$ is also in $\vv_1$.
This contradicts (*) above and finishes the proof.
\end{proof}

\vskip .15in
\noindent
{\bf {Situation III}} $-$  Here 
$\ak u \cap \ak = [z, x]$ with $z \leq x$ in $\ak$.
Then 

$$\ak u \tau \cap \ak \ = \ 
\ak \alpha \tau \alpha^{-1} \cap \ak \ = \ 
[z \tau, x \tau] \ = \ [z \tau, y].$$


\begin{figure}
\centeredepsfbox{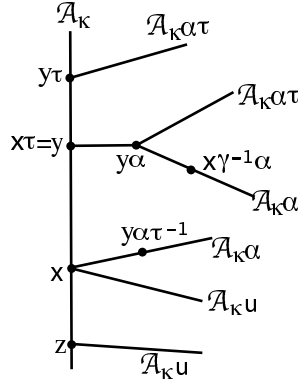}
\caption{
Situation III leading to a contradiction.}
\label{05}
\end{figure}

\noindent
Hence $\ak \alpha \tau \cap \ak \alpha 
= [z \tau \alpha, y \alpha]$
and \ $z \tau \alpha  = y \ \leq_{\alpha} \ y \alpha$
in $\ak \alpha$ $-$ this is the crucial fact.
Now

$$x \gamma^{-1} \alpha \ = \ 
x \beta \alpha \beta^{-1} \ = \ 
x \tau \alpha^{-1} \tau^{-1} \alpha \beta^{-1}
\ = \ y \alpha^{-1} \tau^{-1} \alpha \beta^{-1} $$

$$ \ = \ x \beta^{-1} \ = \ x \tau \alpha \tau^{-1}
\ = \ y \alpha \tau^{-1}.$$

\noindent
Here the bridge of $y \alpha$ to $\ak$ is \
$[y \alpha, y]$ \ (which a priori could be the single point
$y$). So the bridge from 
$y \alpha \tau^{-1}$ to
$\ak$ is $[y \alpha \tau^{-1}, \ y \tau^{-1}]
\ = \ [y \alpha \tau^{-1}, x]$.
On the other hand 
$y \leq  x \gamma^{-1}$ in $\ak$ so
$y \alpha \leq_{\alpha} x \gamma^{-1} \alpha$
in $\ak \alpha$.
It follows that the bridge
from 
$x \gamma^{-1} \alpha$ to $\ak$ 
is $[x \gamma^{-1} \alpha, y]$.
This would imply $x = y$, contradiction.

\vskip .2in
\noindent
{\bf {Situation I}} 
$-$ Surprisingly this is the most difficult
case.
Here 

$$y \alpha^{-1} \tau^{-1} \alpha \pa x
\ {\rm in} \
\ak \alpha, \  
\ \ x \tau > y \ {\rm in } \ \ak,
\ 
\ \ x = y \alpha^{-1} \tau^{-1}, \ 
\ \ y \alpha^{-1} \tau^{-1} \alpha = 
x \tau \alpha^{-1} \tau^{-1}.$$

As \ $y \alpha^{-1} \tau^{-1} \alpha \ \pa \ x$ in
$\ak \alpha$ then $y \alpha^{-1} \tau^{-1} \alpha$ is 
not in $\ak$. Also 

$$x \alpha \ = \ 
(y \alpha^{-1} \tau^{-1}) \alpha \ = \ 
x \tau \alpha^{-1} \tau^{-1} \ = \ x \beta,$$

\noindent
so $x \alpha = x \beta$ $-$ this is a crucial
fact in this proof. 
The bridge from $x \alpha$ to $\ak$ is
$[x \alpha, x]$. Notice also  that

$$x \alpha^{-1} \tau^{-1} \alpha \ \ \pa  \ \ 
y \alpha^{-1} \tau^{-1} \alpha \ \ \pa \ \ x
\ \ \ {\rm in} \ \ \ \ak \alpha,$$

\noindent
so the bridge
from $x \alpha^{-1} \tau^{-1} \alpha$ to
$\ak$ is
$[x \alpha^{-1} \tau^{-1} \alpha, \ x]$.
It follows that 

$${\rm the \ bridge \ from} \ \ 
x \alpha \tau^{-1} \alpha \tau
\ {\rm to} \ \ak \ \ {\rm is }
\ \ [x \alpha^{-1} \tau^{-1} \alpha \tau, \ x \tau] \ = \ 
[x \alpha^{-1} \tau^{-1} \alpha \tau, \ y \alpha^{-1}].$$

\noindent
Now

$$x \alpha^{-1} \tau^{-1} \alpha \tau \ = \ 
(x \alpha^{-1}) \alpha \beta \alpha^{m-1} \ = \ 
x \beta \alpha^{m-1} \ = \ x \alpha \alpha^{m-1}
\ = \ x \alpha^{m}$$

\noindent
Here \ $x \alpha   \ \prec \ x \ \prec \ y \
\prec \ y \alpha^{-1}$ $-$ they are aligned.
It follows from lemma \ref{chax} that
$x, x \alpha$ are in a local axis $\lal$ for
$\alpha$, similarly $y$ is also in a local axis.
Since $y$ is in $[x \alpha^m, x]$, then
also $y, y \alpha^{-1}$ are in $\lal$.
In the same
way 
$(\lal) \tau^{-1} = \lab$ is a local axis
for $\beta$ and 
$x \beta, x, x \tau^{-1}$ are in $\lab$.
Now 

$$x \beta \ = \ x  \tau \alpha^{-1} \tau^{-1} \ = \ 
x \alpha, \ \ \ \ {\rm so} \ \ 
\ \ x \alpha \tau \ = \ x \tau \alpha^{-1} \ = \
y \alpha^{-2}$$

\noindent
Apply $\alpha \beta \alpha^{m-1}
= \tau^{-1} \alpha \tau$ to $y \alpha^{-1}$:

$$(y \alpha^{-1}) \alpha \beta \alpha^{m-1} \ = \
y \beta \alpha^{m-1} \ =  \
(y \alpha^{-1}) \tau^{-1} \alpha \tau \ = \ 
(x \tau) \tau^{-1} \alpha \tau
\ = \ x \alpha \tau \ = \ y \alpha^{-2}.$$


\begin{figure}
\centeredepsfbox{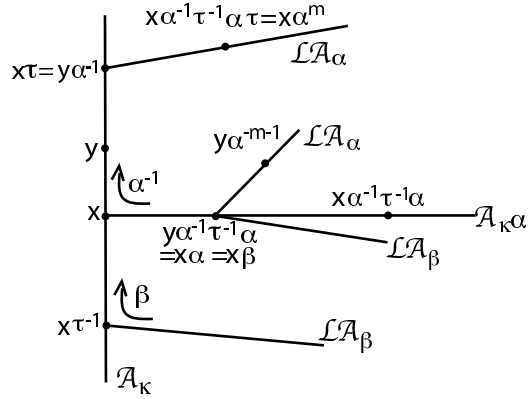}
\caption{
Situation I, the hard case.}
\label{06}
\end{figure}

\noindent
The conclusion is \ \ $y \beta = y \alpha^{-m -1}$ \ \ and
it is in $\lal$.
Now $y$ is not in $\lab$ and the bridge from
$y$ to $\lab$ is $[y,x]$, so the
bridge from $y \beta$ to $\lab$
is $[y \beta, x \beta] \ = \ [y \alpha^{-m-1}, 
\ x \alpha]$.
Therefore $\lal$ and $\lab$ split away from each
other in $x \alpha = x \beta$,
or

$$ \lal \cap \lab \ = \ [x, x \alpha] \ = \ [x, x \beta].$$

\noindent
The homeomorphism $\tau$ conjugates 
the action of $\alpha^{-1}$ in $\lal$ to
the action of $\beta$ in $\lab$ (see fig. \ref{06}).
Now apply $\alpha \tau \alpha^{-m} = \tau \gamma \beta$ to
$x$:

$$ (x \alpha \tau) \alpha^{-m} \ = \
(y \alpha^{-2}) \alpha^{-m} \ = \
y \alpha^{-2-m} \ = \ x \tau \gamma \beta.$$

\noindent
As $x \alpha$ is in $\lab$, then
$x \alpha \tau$ is in $\lal$ and it follows that
$x \tau \gamma \beta$ is in $\lal$.
If $x \tau \gamma \leq x \tau^{-1}$ in $\ak$,
then the bridge from $x \tau \gamma$ to
$\lab$ is $[x \tau \gamma, x \tau^{-1}]$
and so the bridge from
$x \tau \gamma \beta$ to $\lab$ is
$[x \tau \gamma \beta, x \tau^{-1} \beta]$.
But $x \tau^{-1} \beta = x \alpha^{-1} \tau^{-1}$
and 

$$x \alpha^{-1} \tau^{-1} \ <  \
y \alpha^{-1} \tau^{-1}  \ < \ x \ \ \ 
{\rm in } \ \ \ \ak.$$

\noindent
This would imply $x \tau \gamma \beta$ is
not in $\lal$, contradiction.
Notice 

$$x \beta^{-1} \ = \ x \tau \alpha \tau^{-1}
\ = \ y \tau^{-1} \ \ \in \ (x \tau^{-1}, x).$$

\noindent
If $x \tau \gamma$ is in
$[x \tau^{-1}, \ x  \beta^{-1})$
then $x \tau \gamma \beta$ is
in $[x \tau^{-1} \beta, x)$ and
not in $\lal$ either, contradiction again.
Therefore $x \tau \gamma$ is
in $[x \beta^{-1}, x]$.
The case $x \tau \gamma = x$ can only occur
when $p = q = 1$. This case can also be ruled out by
a further argument, but as we are mainly interested
in $|p-2q| = 1$ we assume here that $p > q$.
Then $x \tau \gamma$ is in $[x \beta^{-1}, x)$ and
$x \tau \gamma \beta$ is in $[x, x \beta)$.
We conclude that

$$y \alpha^{-2-m} \ \in \ [x, x \alpha).$$

\vskip .1in
\noindent
{\underline {Claim}} $-$ $y \tau \gamma \beta$ is
in $\lal$.

If $y \tau \gamma \geq x$ in $\ak$, then
$x \ \leq \ y \tau \gamma \ \leq \ y$
in $\ak$.
So $y \tau \gamma \beta$ is in $[x \beta, y \beta]$
or

$$ y \tau \gamma \beta \ \in \ [x \alpha, y \alpha^{-m-1}]
\  \subset \ \lal.$$

\noindent
Notice $x \tau \gamma \beta \in \lal$.
If on the other hand
$y \tau \gamma < x$ in $\ak$, then
$x \tau \gamma < y \tau \gamma < x$ in $\ak$,
and

$$y \tau \gamma \beta \ \in \
(x \tau \gamma \beta, x \beta) \ = \ (x \tau \gamma \beta, 
x \alpha) \ \subset \lal$$ 

\noindent
and again $y \tau \gamma \beta$
is in $\lal$. 

Therefore the claim is proved.

\vskip .1in
It now follows that
$y \tau \gamma \beta \alpha^m \ = \ y \alpha \tau$ is
in $\lal$.
If $y \alpha \ \poa \ x$ in 
$\lal$, then $y \alpha > x$ in $\ak$ as well.
Then
$y \alpha \tau > x \tau = y \alpha^{-1}$ in
$\ak$ and $y \alpha \tau$ is not in $\lal$
contradiction.

Therefore $y \alpha \ \leq_{\alpha} \ x$ in
$\lal$ and so
$y \alpha$ is in $[x, x \alpha)$.
But $y \alpha^{-2-m} \in [x, x \alpha)$.
Since $y$ is in a local axis for $\alpha$ it follows that

$$y \alpha \ = \ y \alpha^{-2-m},
\ \ \ {\rm or} \ \ \ m = -3.$$

\noindent
Since we are assuming $m < -3$ this rules out
this case as well.

This finishes the analysis of situation I and
completes the analysis of the situation 
orders $<$ and $\pa$ agree on $\ak \cap \ak \alpha$.


\begin{figure}
\centeredepsfbox{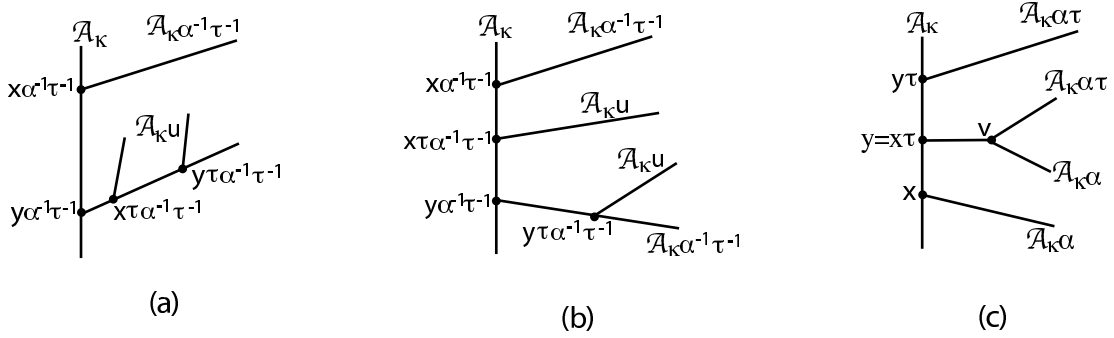}
\caption{
The orientation reversing situation,
\ a. $x \tau > y$, \ b. $x \tau < y$, \ c. $x \tau = y$.}
\label{07}
\end{figure}

\vskip .2in
\noindent
{\bf {Case A.2.2}} $-$ The orders $<$ and $\pa$ 
disagree on $\ak \cap \ak \alpha$.

Notice this is equivalent to $y \alpha^{-1} < 
x \alpha^{-1}$ in $\ak$.
Again use $u = \alpha \tau \alpha^{-1} \tau^{-1} 
= \gamma \tau \alpha^{-1} \tau^{-1} \alpha$.
Then 

$$\ak u \cap \ak \ = \ (\ak \alpha \tau \alpha^{-1}
\cap \ak) \tau^{-1} \ = \
(\ak \alpha \tau \cap \ak \alpha) \alpha^{-1} \tau^{-1}$$

There are the following possibilities:


1) \ If $x \tau > y$ in $\ak$, then
$\ak \alpha \tau \cap \ak \alpha$ is empty
and the bridge from $\ak \alpha \tau$ to
$\ak \alpha$ is $[x \tau, y]$.
Therefore $\ak u \cap \ak = \emptyset$ and
the bridge from $\ak u$ to $\ak$ is
$[x \tau \alpha^{-1} \tau^{-1}, \ y \alpha^{-1} \tau^{-1}]$,
see fig. \ref{07}, a.

2) \ If $x \tau < y$ in $\ak$, then
$\ak \alpha \tau \cap \ak = [x \tau, y]$.
Hence
$\ak \alpha \tau \alpha^{-1} \cap \ak 
= [y \alpha^{-1}, \ x \tau \alpha^{-1}]$,
where the first endpoint is smaller than
the second in $\ak$. Finally

$$\ak u \cap \ak \ = \ [y \alpha^{-1} \tau^{-1}, 
\ x \tau \alpha^{-1} \tau^{-1}]$$

\noindent
and the orders $<$, $\pu$ agree on 
$\ak u \cap \ak$, see fig. \ref{07}, b $-$ because
$y \alpha^{-1} < x \alpha^{-1}$ in
$\ak$ and their images under $u$ satisfy
\ $y \tau \alpha^{-1} \tau^{-1} \ \pu \
x \tau \alpha^{-1} \tau^{-1}$ \ in $\ak u$.

3) \ Finally if $x \tau = y$, then
$\ak \alpha \tau \cap \ak \alpha
= [y, \ v]$, where \
$v \ \leq_{\alpha} \ y$ \ in $\ak \alpha$.
It follows that the intersection
$\ak \alpha \tau \alpha^{-1} \cap \ak 
= [v \alpha^{-1}, y \alpha^{-1}]$,
the first point precedes in $\ak$.
And then

$$ \ak u \cap \ak \ = \ 
[ v \alpha^{-1} \tau^{-1}, \ y \alpha^{-1} \tau^{-1}]
\ = \ 
[t, \ y \alpha^{-1} \tau^{-1}].$$

\noindent
Here if $t$ is not $y \alpha^{-1} \tau^{-1}$ then
$<$ and $\pu$ disagree on $\ak u \cap \ak$
$-$ because $y \tau^{-1} \alpha^{-1} \leq
v \tau^{-1} \alpha^{-1}$ in $\ak$.


\vskip .1in

Now use $\ak u \cap \ak
\ = \ (\ak \alpha^{-1} \tau^{-1} \cap \ak \alpha^{-1})
\alpha$.
Here $\ak \alpha^{-1} \cap \ak = 
[y \alpha^{-1}, \ x \alpha^{-1}]$ the first term
precedes in $\ak$.
Again there are 3 possibilities


\begin{figure}
\centeredepsfbox{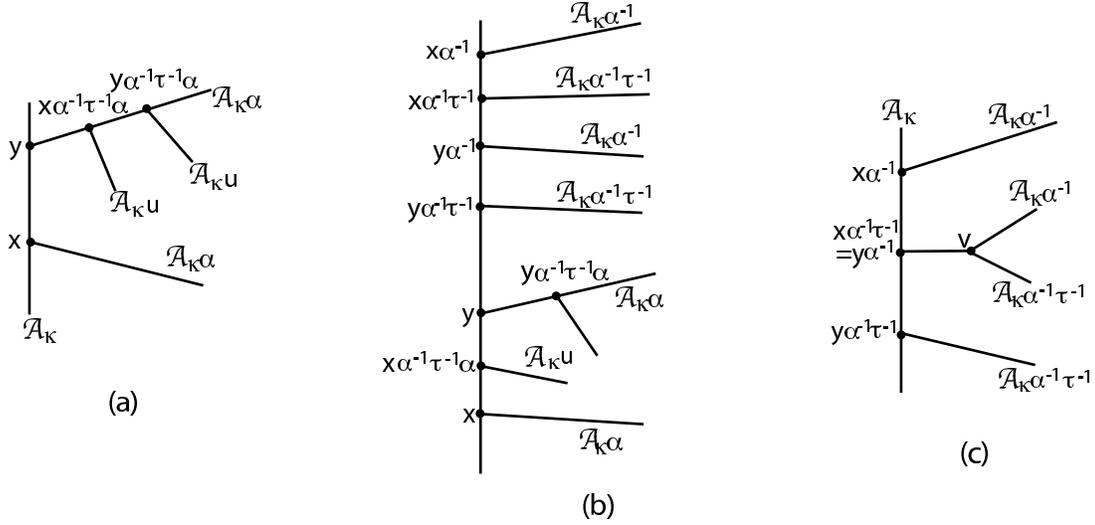}
\caption{
Using 
$\ak u = \ak \alpha \tau \alpha^{-1} \tau^{-1}$:
\ a. $x \alpha^{-1} \tau^{-1} < y \alpha^{-1}$,
\ b. $x \alpha^{-1} \tau^{-1} > y \alpha^{-1}$,
\ a. $x \alpha^{-1} \tau^{-1} = y \alpha^{-1}$.}
\label{08}
\end{figure}


1') \ If $x \alpha^{-1} \tau^{-1} < y \alpha^{-1}$ in $\ak$
then
$\ak \alpha^{-1} \tau^{-1} \cap \ak \alpha^{-1} =
\emptyset$ and the bridge from
$\ak \alpha^{-1} \tau^{-1}$ to $\ak \alpha^{-1}$ 
is $[x \alpha^{-1} \tau^{-1}, \ y \alpha^{-1}]$.
Hence $\ak u \cap \ak = \emptyset$ and the 
bridge from $\ak u$ to $\ak$ is
$[x \alpha^{-1} \tau^{-1} \alpha, y]$, see fig.
\ref{08}, a.

2') \ If $x \alpha^{-1} \tau^{-1} > y \alpha^{-1}$ in
$\ak$, then
$\ak \alpha^{-1} \tau^{-1} \cap \ak \alpha^{-1}
= [y \alpha^{-1}, \ x \alpha^{-1} \tau^{-1}]$
and hence 
$\ak u \cap \ak = [x \alpha^{-1} \tau^{-1} \alpha, y]$
and the orders $<$ and $\pu$ agree on $\ak \cap \ak u$,
see fig. \ref{08}, b $-$ because
$x < y$ in $\ak$ and
\ $x \alpha^{-1} \tau^{-1} \alpha  \   \pu \
y \alpha^{-1} \tau^{-1} \alpha$ \ in $\ak u$.

3') \ If $x \alpha^{-1} \tau^{-1} = y$, then
$\ak \alpha^{-1} \tau^{-1} \cap \ak \alpha^{-1}
= [c, \ y \alpha^{-1}]$ and
$\ak u \cap \ak \ = \ [y, z]$ where $z = c \alpha^{-1}$.
If $z$ is not equal to $y$, then the
orders $<$ and $\pu$ disagree on $\ak u \cap \ak$.


\vskip .1in
Notice both pairs of 3 alternatives are all
mutually exclusive. We match them and obtain 3 possible
situations:


{\bf I} $-$ 
$x \tau > y$ in $\ak$,
$x \alpha^{-1} \tau^{-1} < y \alpha^{-1}$ in $\ak$
and

$$\ak u \cap \ak  \ = \ \emptyset, \ \ \ \
[y \alpha^{-1} \tau^{-1}, \ x \tau \alpha^{-1} \tau^{-1}] \ = \ 
[y, \ x \alpha^{-1} \tau^{-1} \alpha].$$

{\bf II} $-$ 
$x \tau < y$ in $\ak$, $x \alpha^{-1} \tau^{-1} > y \alpha^{-1}$
in $\ak$,

$$ \ak u \cap \ak \ = \ 
[y \alpha^{-1} \tau^{-1}, \ x \tau \alpha^{-1} \tau^{-1}]
\ = \ [x \alpha^{-1} \tau^{-1} \alpha, \ y]$$

\noindent
and the orders $<$, $\pu$ agree on $\ak u \cap \ak$.

{\bf III} $-$
$x \tau = y$, \ $x \alpha^{-1} \tau^{-1} = y \alpha^{-1}$
and

$$ \ak u \cap \ak \ = \
[y, z] \ = \ [t, y \alpha^{-1} \tau^{-1}].$$

\noindent
If $z$ is not $y$ then the orders
$<$, $\pu$ disagree on $\ak u \cap \ak$.


\vskip .1in
We analyse each case in turn:

\vskip .1in
\noindent
{\bf {Situation II}} $-$

Here $x \tau < y$, \ $x \alpha^{-1} \tau^{-1} > y \alpha^{-1}$
and 

$$ y \ = \ x \tau \alpha^{-1} \tau^{-1}, \ \ \ \ 
y \alpha^{-1} \tau^{-1} \ = \ x \alpha^{-1} \tau^{-1} \alpha.$$

\noindent
Suppose first that $[y \alpha^{-1}, x \alpha^{-1}]
\cap [x, y] = \emptyset$.
Since $y \tau  = x \tau \alpha^{-1}$ then
$[y \alpha^{-1}, x \alpha^{-1}]$ is contained
in the set of points $> y$ in $\ak$.

In addition
$y \alpha$ is in $\ak \alpha - \ak$ and
$y \pa y \alpha$.
Hence $y$ is in $(y \alpha^{-1}, y \alpha)$, 
producing a local axis $\lal$ of $\alpha$ which
contains $y$.
Now use $\tau^{-1} \alpha \tau = \alpha \tau \alpha^{-1}
\tau^{-1} \alpha^{m-1}$ applied to $x \alpha^{-1}$:

$$x \alpha^{-1} \tau^{-1} \alpha \tau
\ = \ x \alpha^{-1} \alpha \tau \alpha^{-1}
\tau^{-1} \alpha^{m-1} \ = \
x \tau \alpha^{-1} \tau^{-1} \alpha^{m-1}$$

\noindent
Substitute 
$x \tau \alpha^{-1} \tau^{-1} = y$ in the last
term and $x \alpha^{-1} \tau^{-1} \alpha = y
\alpha^{-1} \tau^{-1}$ in the first term to get

$$(y \alpha^{-1} \tau^{-1}) \tau \ = \ y \alpha^{-1} 
\ = \  y \alpha^{m-1}$$

\noindent
or $y = y \alpha^{m}$. This is impossible because
$y$ is in a local axis of $\alpha$ and $m$ is not zero.

\vskip .1in
From now on in situation II suppose that
$[y \alpha^{-1}, x \alpha^{-1}] \cap [x, y]$ is not
empty.
Since $x \tau \alpha^{-1} = y \tau > y$ in $\ak$,
then $x \alpha^{-1} > y$ in $\ak$.
It follows that $y \alpha^{-1} \leq y$ in $\ak$.


Suppose first that $y \alpha^{-1} < y$ in $\ak$.
There is $r$ in $[y \alpha^{-1}, y]$ which is
fixed by $\alpha$.
Either $r$ is equal to $y$ or $r < y$ in $\ak$.
Let $\uu_1$ (respectively $\uu_2$)
be the component of $T - \{ r \}$ containing
$r \tau$  (respectively $r \tau^{-1}$).
Since

$$x \alpha^{-1} \in \uu_1, \ \ 
x \in \uu_2 \ \ \  \ {\rm then} \ \ \ 
\uu_1 \alpha \ = \ \uu_2.$$

\noindent
If $r < y$ in $\ak$ then also we have $\uu_2 \alpha = \uu_1$.
Otherwise $\uu_2 \alpha = \uu_3$ which is
another component of $T - \{ r \}$ which is
not $\uu_1, \uu_2$.
We will rule out this case, but the result will be
used later on as well, so we state it in more generality:

\begin{lemma}{}{}
Let $\lat$ be a local axis for $\tau$. Let $r$ in $\lat$ which
is fixed by $\alpha$. Let $\uu_1$ (respectively $\uu_2$
be the component of $T - \{ r \}$ containing
$r \tau$ (respectively $r \tau^{-1}$). 
Then $\uu_1 \alpha$ is not $\uu_2$ and
$\uu_2 \alpha$ is not $\uu_1$.
\label{turn}
\end{lemma}

\begin{proof}{}
The proof is as follows:
suppose that either
$\uu_1 \alpha = \uu_2$ or $\uu_2 \alpha = \uu_1$ and
arrive at a contradiction.

First assume that $\uu_1 \alpha = \uu_2$.
Either $\uu_2 \alpha = \uu_1$ or $\uu_2 \alpha$ is another
component $\uu_3$ of $T - \{ u \}$.

Let $\vv_i = \uu_i \tau^{-1}$.
Since 
$\vv_1 \beta \ = \ \vv_1 \tau \alpha^{-1} \tau^{-1} \ = \
\uu_1 \alpha^{-1} \tau^{-1} \ \not = \ \vv_1$,
 we have
that $\vv_1 \beta$ is contained in $\uu_2$.
Therefore $r \beta $ is in $\uu_2$ and
$r \beta \alpha^{m-1}$ is in 
$\uu_2 \alpha^{m-1}$.
Also

$$r \tau^{-1} \alpha \tau \ = \ 
r \alpha \beta \alpha^{m-1} \ = \ 
r \beta \alpha^{m-1}$$

\noindent
As $r \tau^{-1} \in \uu_2$ then $r \tau^{-1} \alpha$ is
in $\uu_2 \alpha$, which is either $\uu_1$ or $\uu_3$.
Therefore $r \tau^{-1} \alpha \tau$ is either
in $\uu_1 \tau \subset \uu_1$ or in $\uu_3 \tau$ again
a subset of $\uu_1$. So $r \tau^{-1} \alpha \tau \in \uu_1$.
Therefore $\uu_2 \alpha^{m-1} \cap \uu_1 
\not = \emptyset$. But both are components of 
$T - \{ r \}$, because $r \alpha = r$,
 so it follows that they are equal. 
As $\uu_2 = \uu_1 \alpha$ then

$$\uu_1 \alpha \alpha^{m-1} \ = \ \uu_1,
\ \ \ {\rm or} \ \ \ \uu_1 \alpha^{m} = \uu_1,
\ \ \uu_2 \alpha^m = \uu_2, 
\ \ \uu_3 \alpha^m = \uu_3 \ \ \ {\rm if \ needed}.$$

\noindent
In case $r \not = y$ this immediately implies
$m$ even.

Now use $r \tau \gamma \beta \alpha^m = r \alpha \tau = r \tau 
\in \uu_1$.
Therefore $r \tau \gamma \beta \ \in \ \uu_1 \alpha^m = \uu_1$.
It follows that

$$r \tau^{-1} \prec r \prec r \tau \gamma \beta$$

\noindent
$-$ recall this means $r$ separates $r \tau^{-1}$
from $r \tau \gamma \beta$.
Applying $\beta^{-1}$ one gets

$$r \tau^{-1} \prec r \beta^{-1} \prec r \tau \gamma 
\ \ \ \ (*)$$

\noindent
Use $r \beta^{-1} = r \tau \alpha \tau^{-1}$:

$$r \tau \in \uu_1 \ \Rightarrow \ r \tau \alpha \in \uu_2, 
\ \ \ r \beta^{-1} = r \tau \alpha^{-1} \tau^{-1} \in \vv_2.$$

\noindent
As $r \tau^{-1}$ is an accumulation point of $\vv_2$,
equation (*) above implies that
$r \tau \gamma$ is in $\vv_2$ or
$r \tau \gamma < r \tau^{-1}$ in $\ak$, which immediately
implies $p > 2q$.

As in the $\rrrr$-covered
case, look at $r \tau \alpha$. If $r \tau \alpha$ is
not in $\vv_2$ then
$r \tau \alpha \tau \not \in \uu_2$ so

$$r \tau \alpha \tau \ = \ 
(r \tau^2) \tau^{-1} \alpha \tau \ = \ 
(r \tau^2) \gamma \beta \alpha^m \ \ \not \in \uu_2
\ \ \ {\rm and} \ \ \ r \tau \gamma \beta \ \not \in \ \uu_2.$$

\noindent
So 
$r \tau^{-1} \prec r \preceq r \tau^2 \gamma \beta$ \ and \
$r \tau^2 \gamma \preceq r \beta^{-1} \prec r \tau^{-1}$. \ 
As $r \beta^{-1}  = r \tau \alpha \tau^{-1} \in \vv_2$,
then 

$$r \tau^2 \gamma \ \in \ \vv_2, \ \ {\rm so} \ \ 
r \tau^2 \gamma < r \tau^{-1} \ \ \ {\rm in } \ \ak.$$

\noindent
As seen before this implies $p > 3q$, which is disallowed
and finishes this case.

If $r \tau \alpha \in \vv_2$ then
$r \beta^{-1} \in \vv_2 \tau^{-1}$.
By $(*)$ 
$r \tau^{-1} \prec r \beta^{-1} \prec r \tau \gamma$, so

$$r \tau \gamma \in \vv_2 \tau^{-1} \ \ 
\Rightarrow \ \ r \tau \gamma < r \tau^{-2} 
 \ \ {\rm in } \ \ \ak.$$

\noindent
As seen before this also implies $p > 3q$ contradiction.

This finishes the analysis of the case $\uu_1 \alpha = \uu_2$.

\vskip .1in
Now suppose that $\uu_2 \alpha = \uu_1$.
If $\uu_1 \alpha = \uu_2$, then this is taken care by the
previous situation. So now assume $\uu_2 \alpha^{-1} = 
\uu_3$ which is not $\uu_1$ or $\uu_2$.

Here use $r \tau^{-1} \alpha \tau = r \alpha \beta \alpha^{m-1}
= r \tau \alpha^{-1} \tau^{-1} \alpha^{m-1}$.
First 

$$r \tau^{-1} \ \in \ \uu_2 \ \ \ \Rightarrow \ \ \ 
r \tau^{-1} \alpha \ \in \ \uu_2 \alpha \ = \ \uu_1 
\ \ \ \Rightarrow \ \ \ r \tau^{-1} \alpha \tau \ \in \ \uu_1.$$

\noindent
On the other hand

$$r \tau \ \in \ \uu_1 \ \ \Rightarrow \ \ 
r \tau \alpha^{-1} \ \in \ \uu_1 \alpha^{-1}  = \uu_2
\ \ \Rightarrow \ \ 
r \tau \alpha^{-1} \tau^{-1} \ \in \ \uu_2 \tau^{-1} \subset \uu_2
\ \ \Rightarrow \ \ 
r \tau \alpha^{-1} \tau^{-1} \alpha^{m-1} \ \in \ 
\uu_2 \alpha^{m-1}.$$

\noindent
From which we conclude that \ $\uu_2 \alpha^{m-1} = \uu_1
= \uu_2 \alpha$.

Now use $r \tau^{-1} \alpha \tau = r \gamma \beta \alpha^m$.
The right side is in $\uu_1 = \uu_2 \alpha$.
The fact that $\uu_2 \alpha^{-1}$ is not $\uu_1$ implies
that $\vv_2 \beta$ is not $\vv_1$, hence 
$\vv_2 \beta$ is contained in $\uu_2$.
We know that $r \gamma$ is $\leq r \tau^{-1}$ in $\lat$ so
it is either in $\vv_2$ or is equal to $r \tau^{-1}$.
Hence $r \gamma \beta$ is either $r \tau^{-1}$ or is
in $\vv_2 \beta$ $-$ in either case it is in $\uu_2$.
Finally $r \gamma \beta \alpha^m$ is in $\uu_2 \alpha^m$
which must be $\uu_1$. But then $\uu_2 \alpha^m = \uu_2
\alpha^{m-1}$ contradiction.

This finishes the analysis of the case $\uu_2 \alpha = \uu_1$
and so 
finishes the
proof of lemma \ref{turn}.
\end{proof}

This finishes the analysis of situation II.

\vskip .2in
\noindent
{\bf {Situation I}} $-$

In this case 
$x \alpha^{-1} \tau^{-1} < y \alpha^{-1}$ in 
$\ak$ and $y < x \tau$ in $\ak$. In addition

$$y \tau \ = \ y \alpha^{-1}, \ \ \ \ \ \ \ \
x \alpha^{-1} \tau^{-1} \alpha \ = \ 
x \tau \alpha^{-1} \tau^{-1} \ \ \ \ \ (*) $$
%
%

\noindent
Here
$x \alpha^{-1} \ > \ y \alpha^{-1} = y \tau$ in $\ak$
(orientation
reversing case) so $x \alpha^{-1} \tau^{-1}  \ > \ y$
in $\ak$.
Therefore $x \alpha^{-1} \tau^{-1} \ \in \ 
(y, y \alpha^{-1})$.
Also $x \tau \ < \ y \tau = y \alpha^{-1}$ in $\ak$,
so one concludes

$$x \alpha^{-1} \tau^{-1}, \ \ \ 
x \tau \ \in \ (y, y \alpha^{-1}) $$

On the other hand \ $y \prec y \alpha^{-1} \prec x \alpha^{-1}$,
so $y \alpha \prec y \prec x$ and $y \alpha$ is
in $\ak \alpha - \ak$. It follows that
$y \alpha^{-1} \prec y \prec y \alpha$ and
$y$ is in a local axis $\lal$ for $\alpha$.
This implies that the translates $[y \alpha^i, y \alpha^{i+1})$
are all disjoint (as $i$ varies in ${\bf Z}$).
%
Use the relation
$\tau^{-1} \alpha \tau \ = \ 
\alpha \tau \alpha^{-1} \tau^{-1} \alpha^{m-1}$
in the form

$$ \alpha^{-1} \tau^{-1} \alpha \tau \alpha^{1-m} 
\ = \ \tau \alpha^{-1} \tau^{-1} $$

\noindent 
applied to $x$ to get

$$ (x \alpha^{-1} \tau^{-1} \alpha) \tau \alpha^{1-m} 
\ = \  x \tau \alpha^{-1} \tau^{-1} \ \ \ \ \ \ (**)$$

\noindent
Now apply the second equality of $(*) $
{\underline {both}} sides of $(**)$ to get

$$ (x \tau \alpha^{-1} \tau^{-1}) \tau \alpha^{1-m} \ = \ 
x \alpha^{-1} \tau^{-1} \alpha \ \ \ \ \ \ \ 
{\rm or }  \ \ \ \ \ \ \ 
(x \tau) \alpha^{-m} \ = \ (x \alpha^{-1} \tau^{-1}) \alpha.$$

\noindent
But \ $x \tau \in (y, y \alpha^{-1})$,  \ so
$x \tau \alpha^{-m} \in
(y, y \alpha^{-1}) \alpha^{-m}$. 
Similarly 
$x \alpha^{-1} \tau^{-1} \alpha$ is in $(y, y \alpha^{-1}) \alpha$.
Since they are equal then $-m = 1$ or
$m = -1$, impossible.

\vskip .2in
\noindent
{\bf {Situation III}} $-$ 

Here $x \tau = y$, $x \alpha^{-1} \tau^{-1} = y \alpha^{-1}$
and 

$$\ak u \cap \ak \ = \ [y,z] \ = \ [t, y \alpha^{-1} \tau^{-1}]$$

\noindent
and if $t \not = y$, then
$<$, $\pu$ disagree on $\ak u \cap \ak$.

Notice that $y \leq z = y \alpha^{-1} \tau^{-1}$
so $y < y \alpha^{-1}$ in $\ak$,
and $y \alpha^{-1}$ is in $\ak - \ak \alpha$.
Also $y \tau \leq y \alpha^{-1}$ in $\ak$.
Now 

$$y \ \prec \ y \alpha^{-1} \ \prec \ x \alpha^{-1}
\ \ \ \Rightarrow \ \ \ 
x \ \prec \ y \ \prec \ y \alpha \ \ \ {\rm all \ in }
\ \ \ak \alpha.$$

\noindent
Hence 
$y \alpha \ \pu \ y$ in $\ak$ and
$y \alpha$ is in $\ak \alpha - \ak$.
Hence $y$ is in $(y \alpha^{-1}, y \alpha)$ and
there is a local axis $\lal$ of $\alpha$ with
$y$ in $\lal$.
Consider the relation $\tau^{-1} \alpha \tau = \alpha \beta
\alpha^{m-1}$. Substitute $\beta = \tau \alpha^{-1} \tau^{-1}$
and rearrange the terms to get
$\alpha^{-1} \tau^{-1} \alpha =
\tau \alpha^{-1} \tau^{-1} \alpha^{m-1} \tau^{-1}$.
Now apply it to
$x$:

$$y = x \alpha^{-1} \tau^{-1} \alpha \ = \ 
x \tau \alpha^{-1} \tau^{-1} \alpha^{m-1} \tau^{-1},$$

\noindent
or $y \tau \alpha^{1-m} = y \alpha^{-1} \tau^{-1}$.
Now $y \tau \in [y, \ y \alpha^{-1}]$, so
$y \tau$ is in $\lal$ and

$$y \tau \alpha^{1-m} \ \in \ [y \alpha^{1-m}, \ y \alpha^{-m}],$$

\noindent
so $y \tau \alpha^{1-m}$ is not in $\ak$.
But $y \alpha^{-1} \tau^{-1}$ is in $\ak$, contradiction.

This finishes the analysis of $\ak u \cap \ak = [x,y]$
with $x$ not equal $y$. Consequently this finishes
the analysis of Case A, $\tau$ acts freely,
which we now proved cannot happen.

\vskip .2in

\section{Case B $-$ $\tau$ has a fixed point, 
$\alpha$ acts freely}
\label{caseb}

Here $\alpha$ has an (actual) axis $\ap$
and so does $\beta$ with axis $\ab = \ap \tau^{-1}$.
Let $Fix(\tau)$ be the set of fixed points
of $\tau$.
As usual there are various possibilities.
This case is very interesting because the topology
of the manifold $M_{p/q}$ will play a key role.

Recall that if $t$ is a point not in a connected set
$B$ of the tree $T$, then the segment $[t,u]$ is
the bridge from $t$ to $B$ if the subsegment
$[x,u)$ does not intersect $B$ and
if $u$ is either in $B$ or is an accumulation point
of $B$. Again the important fact is that
the bridge from $x$ to $B$ is unique:
it is the only embedded path from $x$ to $B$ because
$T$ is a tree. As in case $A$ this will be
explored here. If $u$ is in $B$ we say that
$t$ bridges to $u$ in $B$.

We say that a point $a$ is an ideal point of a local
axis $l$ if $a$ is not in $l$ but is an accumulation
point of $l$. Obviously this implies that $l$ is
not properly embedded in $T$ in the side accumulating
to $a$.

There are two main cases depending on whether $Fix(\tau)$
intersects $\ap$ or not.

\vskip .2in
\noindent
{\bf {Case B.1}} $-$ $Fix(\tau) \cap \ap = \emptyset$.

Then $\kappa$ also has a fixed point $s$. 
Choose $s$ closest to $\ap$, that is,
the bridge $[s,c]$ from
$s$ to $\ap$ has no other fixed point of $\kappa$.
Let $z$ in $[s,c]$ fixed by $\tau$
and closest to $\ap$, that is,
the bridge $[z,c]$ from $z$ to $\ap$ has no other
fixed point of $\tau$ besides $z$.
A priori we do not know whether $z$ is equal to $s$ or
not.
Let $\uu$ be the component of $T - \{ z \}$ containing
$\ap$.

Then $\ab$ is a subset of $\uu \tau \not = \uu$ and
$z$ bridges to $c \tau^{-1}$ in $\ab$.


\begin{figure}
\centeredepsfbox{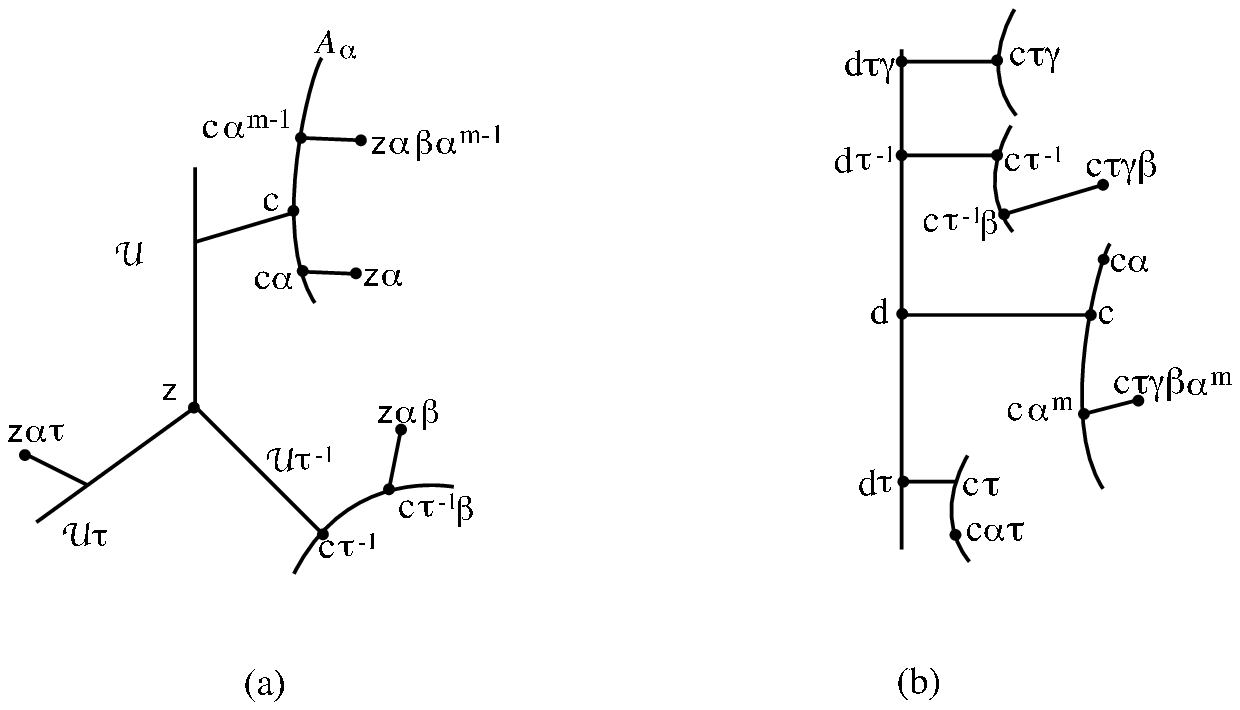}
\caption{
a. The case $\uu \tau \not = \uu$, \
b. The case $\lat \cap \ap = \emptyset$.}
\label{09}
\end{figure}


\vskip .2in
\noindent
{\bf {Case B.1.1}} $-$ Suppose $\uu \tau \not = \uu$.

Then $\uu \tau^{-1} \not = \uu$ as well.
Apply \ $\alpha \tau = \tau \alpha \beta \alpha^{m-1}$ \ to
$z$: the point 
$z$ bridges to $c$ in $\ap$, so $z \alpha$ bridges
to $c \alpha$ in $\ap$. As $c \alpha$ is not $c$ then
$z \alpha$ is in $\uu$, so $z \alpha \tau$ is in
$\uu \tau$, see fig. \ref{09}, a.
On the other hand $z \tau \alpha = z \alpha$ is in
$\uu$ and hence $z$ separates
it from $\ab$. It follows that $z \alpha$ also
bridges to $c \tau^{-1}$ in $\ab$.
Then 

$$z \alpha \tau \beta \ \ \ {\rm bridges \ to} \ \ \ 
c \tau^{-1} \beta \ \ {\rm in} \ \ 
\ab \ \ \ {\rm and} \ \ \  c \tau^{-1} \beta \not = 
c \tau^{-1}, \ \ \ {\rm so} \ \ \  z \tau \alpha \beta \ \in \
\uu \tau^{-1}.$$

\noindent
Therefore $z \tau \alpha \beta$ bridges
to $c$ in $\ap$, so  $z \tau \alpha \beta \alpha^{m-1}$
bridges to $c \alpha^{m-1}$ in $\ap$. This implies 
$z \tau \alpha \beta \alpha^{m-1}$ is in $\uu$,
impossible since it is equal to $z \alpha \tau \in \uu \tau$.

\vskip .1in
We conclude that $\uu \tau = \uu$, which will be assumed
from now on in this proof.

Choose a prong $\eta$ at $z$ which is a subset
of $[z, c]$.
This prong is associated to the component $\uu$ of
$T - \{ z \}$, hence the prong $\eta \tau$ also is 
associated to the component $\uu = \uu \tau$ and
$\eta \cap \eta \tau$ is not just $z$. Let $e$ be
another point in the intersection. Then $e \tau^{-1}, e$
are both in $\eta$ and $e \tau^{-1}$ is not equal $e$ $-$ 
by choice of $z$ as the fixed point of $\tau$ 
in $[z,c]$ closest
to $\ap$.
So either $e$ is in $[z, \ e \tau)$ or
$e \tau$ is in $[z, e)$. In the first case (say) apply
$\tau$ to get $e \tau$ is in $[z, \ e \tau^2)$ and
it now follows that $e \prec e \tau \prec e \tau^2$.
The same alignment of points
 happens in the second case. We conclude that
there is a local axis $\lat$ for $\tau$, with $e$ in 
the local axis.

This construction of a local axis is crucial in case 
B and also in case C of the proof.

\vskip .1in
\noindent
{\bf {Conclusion}} $-$ If $\uu \tau = \uu$ and there is
no fixed point of $\tau$ in $(z,w]$,
then there is a
 local axis of $\tau$ contained in $\uu$ with one
ideal point $z$.

\vskip .2in
\noindent
{\bf {Case B.1.2}} $-$ Suppose that $\lat \cap \ap$
is at most one point.

Let $[d,c]$ be the bridge from 
$\lat$ to $\ap$ $-$ here $d = c$ if
$\lat \cap \ap$ is a single point.
We do the proof for $\lat \cap \ap = \emptyset$, the
case of single point intersection being entirely
similar.
The bridge from $c \alpha \tau$ to $\lat$ is
$[c \alpha \tau, d \tau]$, see fig. \ref{09}, b.
Now the bridge from $c \tau \gamma$ to $\lat$ is
$[c \tau \gamma, d \tau \gamma]$. Here use 
$p$ odd to get
$d \tau \gamma \not = d \tau^{-1}$, so the
bridge from $c \tau \gamma$ to $\ab$ is
$[c \tau \gamma, \ c \tau^{-1}]$.
Therefore 

$$c \tau \gamma \beta \ \  \ {\rm bridges \ to} \ \
c \tau^{-1} \beta \ \ {\rm in } \ \ \ab,
\ \ {\rm hence \ bridges \ to} \ \ 
 d \tau^{-1} \ \ \ {\rm in} \ \  \lat
\ \ \ {\rm and \ to}  \ \ c \ \ \ {\rm in} \ \
\ap.$$

\noindent
Finally 
$c \tau \gamma \beta \alpha^m$ bridges to 
$\ap$ in $c \alpha^{m}
\not = c$ and so bridges to $\lat$ in $c$.

As $c \alpha \tau = c \tau \gamma \beta \alpha^m$,
this implies $c = c \tau$, impossible.
This rules out this case.

\vskip .1in
We conclude that $\lat \cap \ap$ is more than one point.
If $\lat \cap \ap$ is $(z,d]$, then either
$z \alpha = z$ or $\alpha$ has a fixed point
in $\lat$, both impossible.
Therefore from now on in 
case B.1 let 
$\lat \cap \ap = [a, b]$, with $a \not = z$ and
$a$ closest to $z$.
By an abuse of notation $b$ can be $+\infty$, meaning
the intersection is a ray in $\lat$.
Put an order $<$ in $\lat$ so that $a < b$ in $\lat$.
Also let $\pa$ be the order in $\ap$ with 
\ $a  \ \pa  \ b$.

From now on in case B.1 the proof will depend on
whether $\uu \gamma$ is equal to $\uu$ or not.
The arguments here are also very similar to
what will be needed for case C, therefore we will make
the arguments in more generality so that they can
be used in case C,
namely when $\alpha$ has
a fixed point but has a local axis with
certain properties. 
We first specify the conditions under which 
the analysis works.

\vskip .1in
\noindent
{\bf {Conditions}} $-$ Consider two conditions:

\vskip .05in
\noindent
{\bf {Condition F}} $-$ $\tau$ has a fixed point $z$,
$\alpha$
acts freely and $z$ is not in the axis $\ap$.
Let $\ap$ be in the component $\uu$ of $T - \{ z \}$.
There is a fixed point $s$ of $\kappa$ so that
$s$ is either $z$ or $z$ separates $s$ from $\ap$.
Let $(s,c]$ be the bridge from $s$ to $\ap$.
Then $(s,c]$ has no fixed point of $\kappa$
and $(z,c]$ has no fixed point of $\tau$.
Also $\uu \tau = \uu$ and there is a local axis
$\lat$ of $\tau$ in $\uu$ with ideal point $z$.
Finally $\lat \cap \ap = [a,b]$ where $a \not = z$ 
and $a$ is in $(z,b)$.

\vskip .05in
\noindent
{\bf {Condition N}} $-$ $\tau$ has a fixed point $z$; \
$\kappa$ has a fixed point $s$ and $\alpha$ has
 a fixed point $w$ so that $(s,w)$ has no fixed
point of either $\kappa$ or $\alpha$. In addition
either $z = s$ or $z \in (s,w)$ and
 $(z,w)$
has no fixed point of $\tau$.
In addition let $\uu$ be $T_z(w)$ and $\vv$ be
$T_w(z)$. Then $\uu \tau = \uu$ and $\vv \alpha = \vv$.
There is a local axis $\lat$ of $\tau$ in
$\uu$ with one ideal point $z$ and a local
axis $\lal$ of $\alpha$ in $\vv$ with ideal
point $w$. The intersection of $\lal$ and $\lat$ is
$[a,b]$ where $a$ is the closest point
to $z$ and $b$ can be $+\infty$ in $\lat$.

\vskip .05in
Here condition F is for free action of $\alpha$ 
(which is used here) and
condition N is for non free action of $\alpha$
(which is used in Case C).
In either case the order $\pa$ in $\lal$ corresponds
to $a \ \pa \ b$.
This implies the orders $<, \pa$ coincide in the
intersection.
Beware that the order $\pa$ here is in $\lal$ and
not in $({\cal A}_{\tau}) \alpha$ as in case A.

\vskip .05in
\noindent
{\bf {Caution}} $-$ 
An axis is also a local axis.
For the sake of simplicity and to use it
for case C,
we will use the notation $\lal$ even in the
case of $\alpha$ acting freely for the rest of the
proof of case B.1.
In case B.2,
we will return to use the notation $\ap$ for the
axis of $\alpha$.

\vskip .2in
\noindent
{\bf {Case B.1.3}} $-$ $\uu \gamma \not = \uu$.

We first claim that this implies that $\uu \gamma \cap \uu$
is empty. 
Recall that $\partial \uu = z$ and $z \tau = z$.
Notice we do not know a priori that
$z \gamma = z$. If $z \gamma = z$ then $\gamma$ permutes
the components of $T - \{ z \}$ so one has $\uu \gamma
\cap \uu = \emptyset$.
Suppose then that $z \gamma$ is not $z$.
Recall that there is a fixed point $s$ of $\kappa$ with
$z \in [s,w]$ $-$ maybe $s = z$.
If $z \gamma \not = z$, then

$$[s, z] \ \cap \ [s, z \gamma] \ = \ 
[s, t] \ \ \ {\rm with} \ \ t \in [s,z),
\ \ \ {\rm hence} \ \ t \in (z, z \gamma).$$

\noindent
In particular $z$ is not equal to $s$.
Notice $t$ may be equal to $s$.
Here $z$ separates $\uu$ from $s$, hence
$z$ separates $\uu$ from $t$.
Also $z \gamma$ separates $\uu \gamma$ from
$s$, hence $z \gamma$ separates $\uu \gamma$ from
$t$. It follows that $t$ separates 
$\uu$ from $\uu \gamma$ and $\uu \cap \uu \gamma = \emptyset$.
This proves the claim.



\vskip .2in
\noindent
{\bf {Situation I}} $-$ Suppose \ $a   \alpha  \ \pa  \ a$ \ in
$\lal$.

\begin{figure}
\centeredepsfbox{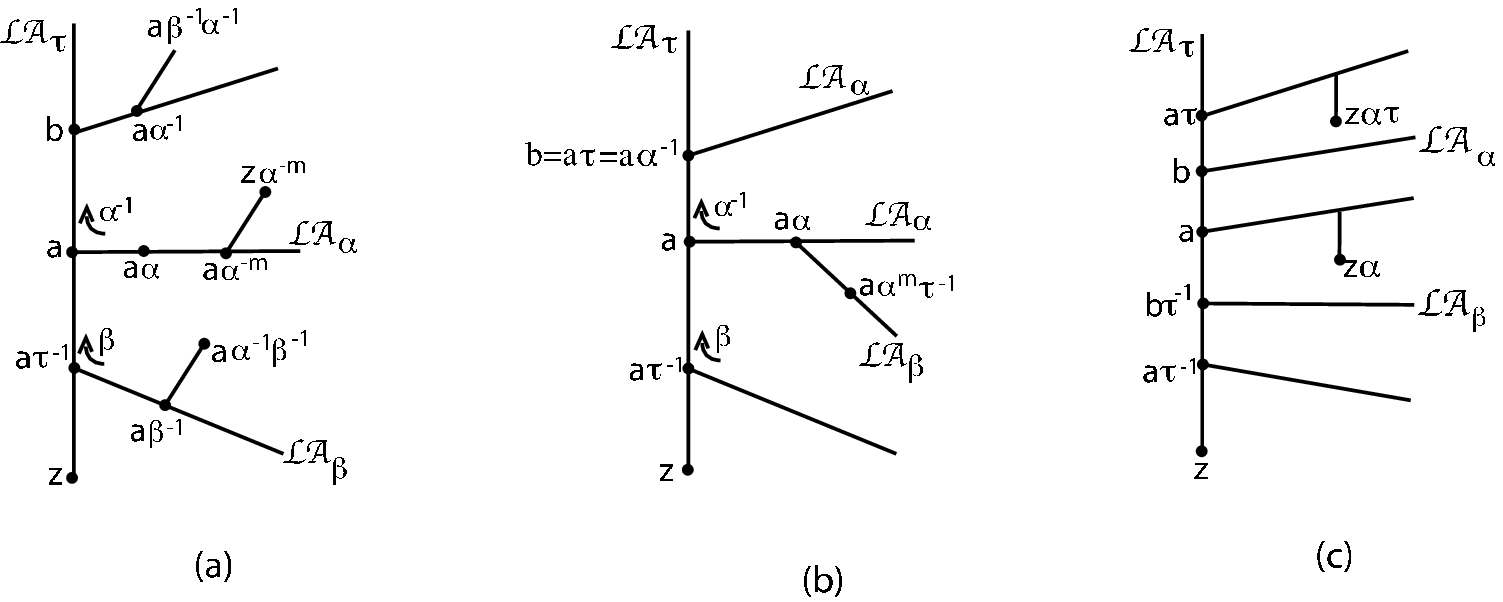}
\caption{
The case 
$\lal \cap \lat = [a,b]$: \  a. Case \ $a \alpha  \ \pa \ a$,
\ $b \ \pa  \ a \alpha^{-1}$, \
b. Case $b = a \tau = a \alpha^{-1}$,
\ c. Case \ $a \tau > b$.}
\label{10}
\end{figure}

\vskip .2in
\noindent
{\bf {Situation I.1}} $-$ Suppose \  $a \alpha^{-1} 
>_{\alpha} \  b$ \ in $\lal$, see fig. \ref{10}, a.

This implies that $a \alpha$ is not in
$\lat$, see fig. \ref{10}, a.
Also this implies $b$ is finite. Notice that

$$z \tau^{-1} \alpha^{-1} \tau \ = \
z \alpha^{-m} \beta^{-1} \gamma^{-1} \ = \
z \alpha^{-m} \tau \alpha \tau^{-1} \gamma^{-1}$$

\noindent
The point $z$ bridges to $\lal$ in $a$.
Hence $z \tau^{-1} \alpha^{-1} = z \alpha^{-1}$
bridges to $\lal$ in $a \alpha^{-1}$, 
so $z \alpha^{-1}$ is in $\uu$ and $z \alpha^{-1} \tau$
is also in $\uu$, which is invariant under $\tau$.
Since $\uu \gamma \cap \uu = \emptyset$, then

$$z \alpha^{-m} \tau \alpha \ \not \in \uu
\ \ \ {\rm and \ it \ bridges \ to} \ \ 
\lal \ \ {\rm in} \ a \ \ \Rightarrow \ \
z \alpha^{-m} \tau \ \ \ {\rm bridges \ to} \ \
\lal \ \ {\rm in} \ 
a \alpha^{-1}$$

\noindent
 and hence bridges to $\lat$ in
$b$.
But $z \alpha^{-m}$ bridges to $\lal$ in $a \alpha^{-m}$
so bridges to $\lat$ in $a$. So $z \alpha^{m} \tau$ 
bridges to $\lat$ in $a \tau$. This implies
$a \tau = b$ and also 
that $\tau$ is increasing in
$(\lat, <)$.


In addition 

$$\lab \ = \ (\lal) \tau^{-1} \ \ \ {\rm so} \ \ \
\lab \cap \lat \ = \ [a \tau^{-1}, a] \ = \ 
[a \tau^{-1}, b \tau^{-1}]$$

\noindent
and $a \beta^{-1}$ is not in $\lat$ and
bridges to $\lat$ in $a \tau^{-1}$. So this point
bridges to $\lal$ in $a$ and $a \beta^{-1} \alpha^{-1}$
bridges to $\lal$ in $a \alpha^{-1}$. As
a result $a \beta^{-1} \alpha^{-1}$ is in $\uu$.

Also $a \alpha^{-1}$ bridges to $\lat$ in $b = a \tau$.
Hence it bridges to $\lab$ in $a$. This implies
that $a \alpha^{-1} \beta^{-1}$ bridges
to $\lab$ in $a \beta^{-1}$ so again
$a \alpha^{-1} \beta^{-1}$ is in $\uu$.
Now  \ 
$(a \beta^{-1} \alpha^{-1}) \gamma \ = \ 
a \alpha^{-1} \beta^{-1}$.
\ Which implies $\uu \gamma \cap \uu$ is not empty.
This contradicts the first claim in Case B.1.3.

Situation I.1 cannot happen.

\vskip .2in
\noindent
{\bf {Situation I.2}} $-$ Suppose \ $a \alpha^{-1} \leq_{\alpha}
b$ \ in $\lal$.

Similarly to the arguments  in situation I.1,
$z \alpha^{-1} \tau$ is in
$\uu$,
so $z \alpha^{-m} \tau \alpha$ is not
in $\uu$ so 

$$z \alpha^{-m} \tau \alpha \ \ \ {\rm bridges \ to} \ \ \ 
\lal \ \ {\rm in} \  a,  \ \ z \alpha^{-m} \tau
\ \ \ {\rm  bridges \   to} \ \ 
\lal \ \ {\rm in} \  \ a \alpha^{-1}.$$

\noindent
Also
$a \alpha^{-1} \leq_{\alpha}  b$ 
in $\lal$, hence $a \alpha^{-1}$ is
in $\lat$ and $a \alpha^{-1} \leq b$ in $\lat$ as well.
On the other hand $z \alpha^{-m}$ bridges to 
$\lat$ in $a$ so $z \alpha^{-m} \tau$ bridges
to $\lat$ in $a \tau$.
From this it follows that $a \tau \geq a \alpha^{-1}$
in $\lat$.
There are two possibilities:

\vskip .05in
The first possibility is that 
$a \alpha^{-1} \not = b$.
In this case $z \alpha^{-m} \tau$ 
bridges to $\lal$ in $a \alpha^{-1}$
which is in the interior of $[a,b]$, hence this point
also bridges to
$\lat$ in $a \alpha^{-1}$.
It follows that

$$a \tau \ = \ a \alpha^{-1} \ \ \ \Rightarrow \ \ \ 
a \beta^{-1} \ = \  a \tau^{-1} \ \ \ {\rm bridges \ to} \ \
\lal \ \ {\rm in} \ a.$$

\noindent
Then $a \beta^{-1} \alpha^{-1}$ bridges to
$\lal$ in $a \alpha^{-1}$ so 
is in $\uu$.
As before consider $a \alpha^{-1} \beta^{-1}$.
Here $a \alpha^{-1}$ is either in $\lab$ or
bridges to $\lab$ in $b \tau^{-1}$ (the top
intersection of $\lab$ with $\lat$).
If $a \alpha^{-1}$ in $\lab$ then 
$a \alpha^{-1} \beta^{-1}$ is in $\lab$ so in $\uu$,
as above contradiction.
If it bridges to 
$\lab$ in $b \tau^{-1}$ then
$a \alpha^{-1} \beta^{-1}$ bridges to $\lab$ in
$b \tau^{-1} \beta^{-1} = b \alpha \tau^{-1}$.
Since in this case

$$b \alpha \ > \ a \ \ {\rm in} \ \ \lat,
\ \ \ {\rm then} \ \ b \alpha \tau^{-1}  \ > \  
a \tau^{-1} \ \ {\rm in} \ \lat \ \ \ 
\Rightarrow \ \ \ 
a \alpha^{-1} \beta^{-1} \ \in \ \uu,$$

\noindent
again
a contradiction.

\vskip .05in
The second possibility is that $a \alpha^{-1} = b$.
Here we have to split further into two options:

Recall that $a \tau \geq a \alpha^{-1}$ in $\lat$.
First consider the case that $a \tau = a \alpha^{-1}$,
see fig. \ref{10}, b.
We have the equalities
$a \beta^{-1} = a \tau \alpha \tau^{-1} = a \tau^{-1}$.
Use

$$(a \alpha^m ) \tau^{-1} \alpha^{-1} \tau
\ = \ a \alpha^m \alpha^{-m} \beta^{-1} \gamma^{-1}
\ = a \beta^{-1} \gamma^{-1}
\ = \ a \tau^{-1} \gamma^{-1} \not \in \uu$$

\noindent
Hence $a \alpha^m \tau^{-1} \alpha^{-1}$ is not in $\uu$
and bridges to $\lal$ in $a$,
$a \alpha^m \tau^{-1}$ bridges to $\lal$ in $a \alpha$.
But 

$$a \alpha^m \ \in \ \lal \ \ \Rightarrow \ \ 
a \alpha^m \tau^{-1} \ \in \ \lab \ \ \ 
\Rightarrow \ \ \ 
\lal \cap \lab \ = \ [a , a \alpha],$$

\noindent
see fig. \ref{10}, b.
Now evaluate $\gamma^{-1} = \beta \alpha \beta^{-1} \alpha^{-1}$
on $a \tau^{-1}$:

$$(a \tau^{-1}) \gamma^{-1} \ = \
(a \beta^{-1}) \beta \alpha \beta^{-1} \alpha^{-1}
\ = \ a \alpha \beta^{-1}  \alpha^{-1}.$$

\noindent
Notice that $a \alpha$ is in $\lab$ so $a \alpha \beta^{-1}$
is in $\lab$. Either $a \alpha \beta^{-1}$ is in
$\lal$ and then $a \alpha \beta \alpha^{-1}$  is in
$\lal \subset \uu$ (contradiction) $-$
or 

$$a \alpha \beta ^{-1} \ \not \in \ \lal
\ \ \ {\rm so \
bridges \ to} \ \  \lal \ \ \ {\rm in} \ \  a \ \ \
{\rm and} \ \
a \alpha \beta \alpha^{-1} \ \ \ {\rm bridges \ to} \ \ 
\lal \ \ {\rm in} \ \ a \alpha^{-1}$$

\noindent
 and again
this point is in $\uu$.
In either case $\uu \gamma \cap \uu \not = \emptyset$,
contradiction.

The last option of the second possibility
$a \alpha^{-1} = b$ 
is that $a \tau > b = a \alpha^{-1}$ in
$\lat$.
Then

$$b \tau^{-1} \ = \ a \tau^{-1} \beta \ < \ a
\ \ {\rm in} \ \lat 
\ \ \ \Rightarrow \ \ \ 
\lal \ \cap \ \lab \ = \ \emptyset,$$

\noindent
see fig. \ref{10}, c.
Here use $\alpha \tau = \tau \alpha \beta \alpha^{m-1}$
applied to $z$: The point $z \alpha$ bridges to
$a$ in $\lat$ and $z \alpha \tau$ bridges to $a \tau$
in $\lat$. Since $a \tau > b$, then $z \alpha \tau$
bridges to $b = a \alpha^{-1}$ in $\lal$.

On the other hand $z \alpha$ bridges to $b \tau^{-1}$ in
$\lab$ hence $z \alpha \beta$ bridges to 
$b \tau^{-1} \beta$ in $\lab$, hence to $a$ in
$\lal$. Finally
$z \alpha \beta \alpha^{m-1}$ bridges  to
$a \alpha^{m-1}$ in $\lal$.
Since $m$ is not $0$ this is a contradiction.

We conclude that situation I cannot happen.

\vskip .2in
\noindent
{\bf {Situation II}} $-$ \ $a \alpha^{-1} \ \pa \ a$ \ in $\lal$.

\vskip .2in
\noindent
{\bf {Situation II.1}} $-$ $a \alpha^{-m}$ is not in $\lat$.
Here use

$$z \alpha^{-1} \tau \ = \
z \tau^{-1} \alpha^{-1} \tau \ = \
z \alpha^{-m} \tau \alpha \tau^{-1} \gamma^{-1}$$

\noindent
is in $\uu$,
so $z \alpha^{-m} \tau \alpha$ is not in $\uu$.
It bridges to $\lal$ in $a$, hence $z \alpha^{-m} \tau$ bridges
to $\lal$ in $a \alpha^{-1}$ and hence bridges
to $\lat$ in $a$. On the other hand $z \alpha^{-m}$
bridges to $\lal$ in $a \alpha^{-m}$, so bridges
to $\lat$ in $b$.
It follows that $z \alpha^{-m} \tau$ bridges to
$\lat$ in $b \tau$ which then must be $a$.
So $a < a \tau^{-1}$ in $\lat$.


\begin{figure}
\centeredepsfbox{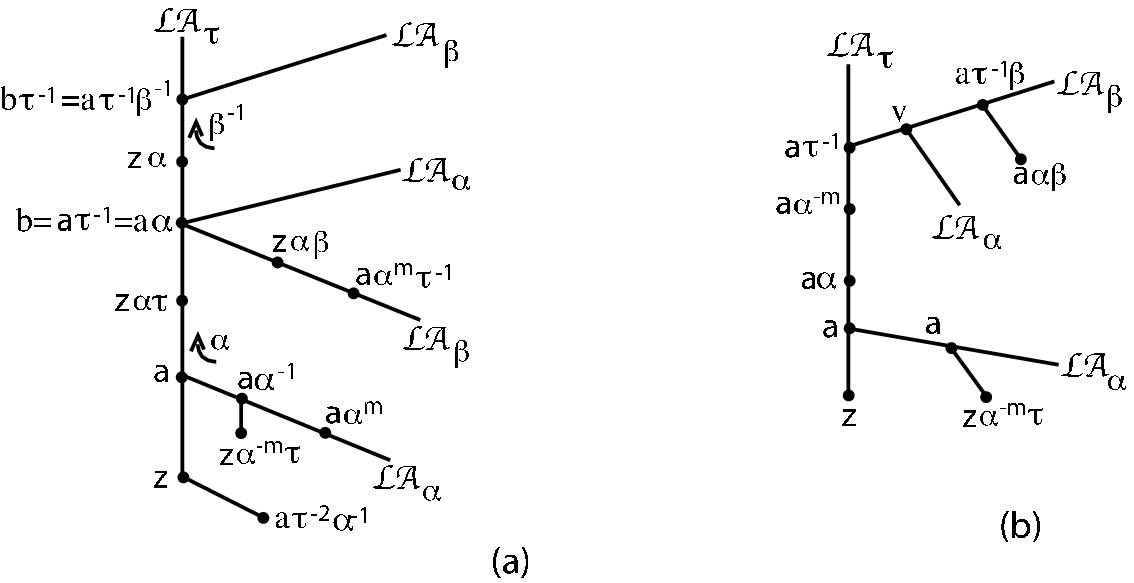}
\caption{
Case $a \alpha^{-1} \ \pa \ a$ in
$\lal$:
a. Picture when
$a \alpha^{-m} \not \in \lat$,
$a \alpha = a \tau^{-1}$,
b. Picture when $a \alpha^{-m} \in \lat, \ 
a \tau^{-1} \beta \not \in \lal$.}
\label{11}
\end{figure}

Notice $\lab \cap \lat$ is equal to $[a \tau^{-1}, b \tau^{-1}]$
and this intersects $\lal$ in $a \tau^{-1} = b$.

Suppose first that $a \alpha$ is not $a \tau^{-1} = b$.
Here

$$a \beta^{-1} \ \ \ {\rm bridges \ to}
\ \ \ \lab \ \ {\rm  in} \  a \tau^{-1} \beta^{-1},
\ \ \ {\rm so \ bridges \ to} \ \ \lal \ \ {\rm in} \ a \tau^{-1}.$$

\noindent
Then $a \beta^{-1} \alpha^{-1}$ bridges to 
$\lal$ in $a \tau^{-1} \alpha^{-1} \not = a$.
It follows that $a \beta^{-1} \alpha^{-1}$ is in $\uu$.

On the other hand $a \alpha^{-1}$ bridges to $\lab$ in $a \tau^{-1} = b$,
so $a \alpha^{-1} \beta^{-1}$ bridges to $\lab$ in $b \beta^{-1}$ which
is not $b$ and it follows that $a \alpha^{-1} \beta^{-1}$ is also
in $\uu$. As seen before  this implies $\uu \gamma \cap \uu$ is
not emptyset, contradiction.

\vskip .05in
The second option in situation II.1 is that 
$a \alpha = a \tau^{-1}$, see fig. \ref{11}, a.

Apply \ $\alpha^{-m}  \beta^{-1} \gamma^{-1}
\ = \ \tau^{-1} \alpha^{-1} \tau$ \  to $a \alpha^m$.
The right side becomes
$a \beta^{-1} \gamma^{-1}$. Here

$$a \beta^{-1} \ \in \ \uu \ \ \ \ \Rightarrow \ \ \
a \beta^{-1} \gamma^{-1} \ \not \in \ \uu \ \ \ \Rightarrow \ \ \ 
a \alpha^m \tau^{-1} \alpha^{-1} \ \not \in \ \uu$$

\noindent
and bridges to $\lal$ in $a$. \
It follows that $a \alpha^m \tau^{-1}$ bridges to 
$\lal$ in $a \alpha = a \tau^{-1} = b$.
But $a \alpha^m$ is in $\lal$,
so  $a \alpha^m \tau^{-1}$ is in $\lab$.
Consequently 
$\lal \cap \lab = a \tau^{-1} = b$, \  see fig. \ref{11}, a.

%
%
%
%

The point
$a \beta^{-1}$ is in $\uu$, hence 

$$a \beta^{-1} \gamma^{-1} \ = \ a \alpha \beta^{-1} \alpha^{-1}
\ = \ a \tau^{-1} \beta^{-1} \alpha^{-1} \ = \ 
a \tau^{-2} \alpha^{-1}$$

\noindent
is not in $\uu$.
Not only that, but also $a \beta^{-1} \gamma^{-1}$ is
not equal to $z$ $-$ else some point near $a \beta^{-1}$
in $\uu$ will have image under $\gamma$ in
$\uu$, which is disallowed.
Then 

$$z \ \in \ (a, a \tau^{-2} \alpha^{-1})
\ \ \Rightarrow \ \ 
z \alpha \ \in \ (a \alpha, a \tau^{-2}) \ = \ 
(a \tau^{-1}, a \tau^{-2}) \ \ \Rightarrow \ \ 
z \alpha \tau \ \in \ (a, a \tau^{-1}).$$

\noindent
In particular $z \alpha \tau$ is in $\lal$ and
$z \alpha \tau \alpha^{1-m}$ is in $\lal$ as well.
This point is equal to $z \alpha \beta$.

On the other hand 

$$z \alpha \ \in \ (a \tau^{-1}, a 
\tau^{-2}) \ = \ (a \tau^{-1}, a \tau^{-1} \beta^{-1})
\ \ \ \Rightarrow \ \ \ 
z \alpha \beta\ \in \ 
(a \tau^{-1}, a \tau^{-1} \beta).$$

\noindent
But then
$z \alpha \beta$ is not in $\lal$, contradiction.

This finishes the analysis of situation II.1,
$a \alpha^{-m}$ is not in $\lat$.

\vskip .2in
\noindent
{\bf {Situation II.2}} $-$ $a \alpha^{-m}$ is in $\lat$.

In particular $a \alpha$ is in $(a,b]$.
Here $z \alpha^{-m} \beta^{-1} \gamma^{-1} =
z \tau^{-1} \alpha^{-1} \tau$ is in $\uu$.
As usual this implies 
$z \alpha^{-m} \tau \alpha$ is not in $\uu$
and bridges to $\lal$  
in $a$ and $z \alpha^{-m} \tau$ bridges to $\lal$
in $a \alpha^{-1}$, see fig. \ref{11}, b; so
$z \alpha^{-m} \tau$
bridges to $\lat$ in $a$.
So

$$z \alpha^{-m} \ \ \ {\rm bridges \ to} \ \lat \ \ {\rm in}
\ a \tau^{-1} \ \ \ \Rightarrow \ \ \ 
a \tau^{-1} \ > \ a \ \ {\rm in} \ \ \lat.$$

Notice $z \alpha^{-m}$ bridges to $\lal$ in $a \alpha^{-m}$.
If \ $a \alpha^{-m} \pa \ b$ \ in $\lal$, then
$z \alpha^{-m}$ also bridges to $\lat$ in $a \alpha^{-m}$
and $a \alpha^{-m} = a \tau^{-1}$.
If 

$$a \alpha^{-m} \ = \ b \ \ \ {\rm then} \ \ 
z \alpha^{-m} \ \ \ {\rm bridges \ to} \ \ 
\lat \ \ \ {\rm in \ a  \ point } \ 
\geq a \alpha^{-m},$$

\noindent
that is, $a \tau^{-1} \geq a \alpha^{-m}$ in
$\lat$.
In any case 
$a \alpha^{-m} \leq a \tau^{-1}$ in $\lat$ and 
$a \alpha < a \tau^{-1}$ in $\lat$.

Now compute $a \gamma = a \alpha \beta \alpha^{-1} \beta^{-1}$.
Here $a \alpha$ is in $[a, a \tau^{-1}]$
and bridges to $\lab$ in $a \tau^{-1}$.
Hence $a \alpha \beta$ bridges to $\lab$
in $a \tau^{-1} \beta$.
There are two options:
First if $a \tau^{-1} \beta$ is not in $\lal$, then
$a \alpha \beta$ bridges to a point
$v$ in $\lal$ and  $v \in (a, a \tau^{-1} \beta)$ $-$
see fig. \ref{12}, b.
Here 
$v$ could be in $\lat$. 
Then

$$a \alpha \beta \alpha^{-1} \ \ \ {\rm bridges \ to \ a \ point} \ 
\ v \alpha^{-1} \ \ {\rm in} \ \ \lal
\ \ \Rightarrow \ \ {\rm it \ bridges \ a \ point} \ \ 
c \ \ {\rm in}  \ \lab, \ \ c \ \in \
(b \tau^{-1}, a \tau^{-1} \beta).$$

\noindent
It follows that 
$a  \gamma = a \alpha \beta \alpha^{-1} \beta^{-1}$ bridges
to a point in $\lab$ which is not $a \tau^{-1}$, hence
$a \gamma$ is in $\uu$, contradiction.

The second option here is
that $a \tau^{-1} \beta$ is in $\lal$.
Here $a \tau^{-1}$ is in $\lal$.
Then consider $a \tau^{-1} \alpha^{-1}$ which is in $\lal$ and hence
in $\uu$. Then

$$(a \tau^{-1} \alpha^{-1}) \alpha \beta \alpha^{-1}
\ = \ a \tau^{-1} \beta \alpha^{-1}$$

\noindent
is in $\lal$ and 
\ $a \tau^{-1} \beta \alpha^{-1}  \ \pa \ 
a \tau^{-1} \beta$  \ in $\lal$.
Therefore

$$a \tau^{-1} \beta \alpha^{-1} \ \ \ 
{\rm bridges \ to \ a \ point \ in} \ 
\lab \ \ \ {\rm contained \ in} \ \ 
(b \tau^{-1}, a \tau^{-1} \beta).$$

\noindent
Apply $\beta^{-1}$ $-$ the resulting point
bridges to a point in
$\lab$ which is not $a \tau^{-1}$, hence
$(a \tau^{-1} \alpha^{-1}) \gamma$ is in $\uu$, again a contradiction.

This finishes the analysis of situation $II$.
Hence this finishes the analysis of case B.1.3,
$\uu \gamma$ is not equal to $\uu$.

%
%
%
%
%

\vskip .2in
\noindent
{\bf {Case B.1.4}} $-$ Suppose $\uu \gamma = \uu$.

Since the boundary $\partial \uu$ in $T$ is the point $z$ this
implies that $z \gamma = z$.
 Here $(\lat) \gamma \cap \lat \not = \emptyset$,
choose
$c \gamma$ in this intersection.
So $c, c \gamma$ are disjoint and in $\lat$. If follows
that $z, c, c \gamma$ are aligned (the particular
order is not important) and $c$ is in a local
axis of $\gamma$. But $c \gamma^{-q} = c \tau^p $ is
also in $\lat$ and it follows easily that the local
axis is contained in and therefore equal to
the local axis $\lat$ of $\tau$
so $\gamma, \tau$ and hence $\kappa$ leaves
$\lat$ invariant. This sort of argument
will be used from time to time from now on.

Here the ideal would be to apply the proof of
case A, where $\tau$ acted freely and $\at$
was invariant by $\gamma$ and $\tau$. We already
have $\lat$ invariant under $\gamma$ and $\tau$,
however $\lat$ is not properly embedded in
$T$ - at least in the $z$ direction.
In order to apply the proof of case A, we analyse
the relative positions of $(\lat) \alpha, (\lat) \alpha \tau$
and so on. In particular for that analysis to work
we must have $(\lat) \alpha$ contained in $\uu$ and so
on. So first we 
do preparation work, showing all images
of the local axis are in $\uu$ and then we can
apply the proof of case A.

For simplicity of notation in case B.1.4 we do the
following: 
$\kk$ will denote the local axis
$\lat$ which is contained in $\uu$ and has an
ideal point $z$.
Again as we want to use this in section C as well,
we will consider a local axis $\lal$ for $\alpha$.
The key result is the following:

%
%
%


\begin{lemma}{}{}
We have $\kk \alpha \subset \uu$,
$\kk \alpha^{-1} \subset \uu$ and $\kk \alpha \tau \alpha^{-1}
\subset \uu$.
\label{probl}
\end{lemma}

\begin{proof}{}
We treat each case in turn:

\vskip .1in
\noindent
{\bf {Problem 1}} $-$ Is $\kk \alpha \subset \uu$?

Suppose not. Then as $a \alpha$ is in $\lal$ contained
in $\uu$ there is $t$ in $\kk$ with $t \alpha = z$ or
$z \alpha^{-1}$ is in $\kk$, see fig. \ref{13}, a.
Here $z$ bridges to $a$ in $\lal$ so $z \alpha^{-1}$
bridges to $a \alpha^{-1}$ in $\lal$. So $z \alpha^{-1}$
can only be in $\kk$ if
$b$ is in $(z, z \alpha^{-1})$
and $a \alpha^{-1} = b$. In particular 
\ $a \alpha \ \pa \ a$ \ in $\lal$.


\begin{figure}
\centeredepsfbox{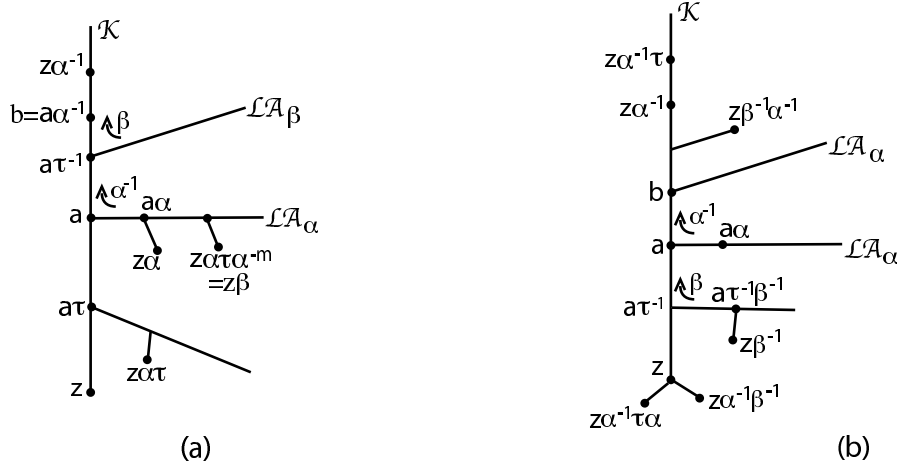}
\caption{
Analysing $z \alpha^{-1} \in \kk$:
a. Picture when $a \tau \in [z,a)$,
\ b. Picture when $a \tau^{-1} \in [z,a)$.}
\label{12}
\end{figure}

There are two possibilities depending on whether
$\tau$ is expanding away from $z$ or not:

First suppose $a \tau$ is in $[z, a)$, see fig.
\ref{12}, a.
As $z \alpha$ bridges to $a$ in $\kk$ then $z \alpha \tau$
bridges to $a \tau$ in $\kk$ so bridges to
$a$ in $\lal$. Then $z \alpha \tau \alpha^{-m}$  bridges
to $a \alpha^{-m}$ in $\lal$.
The point $z \alpha \tau \alpha^{-m}$
is equal to $z \beta$ and
bridges to $a$ in $\kk$ so bridges to $a \tau^{-1}$ in
$\lab$. But $z$ also bridges 
to $a \tau^{-1}$ in $\lab$, contradiction.

The second option is $a \tau > a$ in $\kk$, see fig. \ref{12}, b.
Here $z \beta^{-1} $ bridges to $a \tau^{-1} \beta^{-1}$ 
in $\lab$ and
so to $a$ in $\lal$.
Hence 

$$z \beta^{-1} \alpha^{-1} \ \ \ {\rm bridges \ to} \ 
a  \alpha^{-1} \ \ {\rm in} \  \lal
\ \ \ \Rightarrow \ \ \ 
z \beta^{-1} \alpha^{-1} \ \in \ \uu.$$

\noindent
On the other hand $z \alpha^{-1} \beta^{-1} =
z \alpha^{-1} \tau \alpha \tau^{-1}$.
Here 

$$z \alpha^{-1} \tau \ \in \ \kk \ \ \Rightarrow \ \ 
z \alpha^{-1} \ \in \ (z, z \alpha^{-1} \tau) \ \ \Rightarrow \ \ 
z \alpha^{-1} \tau \alpha \ \not \in \ \uu \ \ \Rightarrow \ \ 
z \alpha^{-1} \beta^{-1} \ \not \in \ \uu.$$

\noindent
But $z \beta^{-1} \alpha^{-1} \gamma
= z \alpha^{-1} \beta^{-1}$, leading
to $\uu \gamma \not = \uu$,
contradiction
to case B.1.3.

\vskip .05in
So we obtain $z \alpha^{-1} \in \uu$ is impossible.
Hence $\kk \alpha \subset \uu$.
If 
$\kk \alpha$ intersects $\kk$ in at most one point we can
use the analysis of Case B.1.2
(or of case A) and disallow it.
If

$$\kk \ \cap \ \kk \alpha \ = \ (z, t),
\ \ \ {\rm then} \ \ 
\kk, \ \kk \alpha \ \ \ {\rm share \ a \ ray}.$$

\noindent
The orientations in
$\kk$ and $\kk \alpha$ may agree or not. In the first
case $z \alpha = z$ and in the second case there is
a fixed point $r$ of $\alpha$ in $\lat = \kk$.
If $z \alpha = z$, then $z$ is a global fixed point,
impossible by non trivial action. In the second option 
let $\uu_1$ (respectively $\uu_2$) be the component
of $T - \{ r \}$ containing $r \tau$ (respectively
$r \tau^{-1}$). The condition
$\kk \alpha \cap \kk = (z,t)$ implies that $\uu_1 \alpha 
= \uu_2$. This is now disallowed by lemma \ref{turn}.

Now consider the situation that $\kk$ has another ideal
point $v$.
Then $v \kappa = v$.
Suppose first that $v$ is in $\lal$. Here we split
into cases: if $\alpha$ acts freely then $v$ is
a fixed point of $\tau$ in the axis of $\alpha$ and
this falls under case B.2. 
Consider then the case that
$\alpha$ does not act freely. Then $(w,v)$ has no
fixed point of  $\alpha$ (as $v$ is in $\lal$) and
also no fixed point of $\tau$ or $\gamma$.
Also $T_w(v)$ is invariant under $\alpha$ and
$T_v(w)$ is invariant under $\tau$. 
Then $v$ in $\lal$ is 
disallowed by lemma \ref{noal}.

\vskip .05in
It follows that 
$v$ has the same properties as $z$. 
In any case one obtains that

$$\kk \alpha \ \cap \ \kk \  
 = \ [t, r], \ \ t  \ \not = \ r,
\ \ t \ \ \ {\rm closest \ to} \ \ z$$

\noindent
and if $\kk$ is not properly embedded in the
other direction then $r$ is an actual point in $\kk$.
Then $\kk \alpha \tau \cap \kk = [t \tau, r \tau]$.
So the intersections are the same as occurred
in Case A so far.


\vskip .2in
\noindent
{\bf {Problem 2}} $-$ Is $\kk \alpha^{-1} \subset \uu$?

This is similar to problem 1.
As before if $\kk \alpha^{-1}$ not contained in $\uu$, then
$z \in \kk \alpha^{-1}$ and $z \alpha \in \kk$.
This can only happen if $b \in (z, z \alpha)$,
$a \alpha = b$ and $a \alpha^{-1} \ \pa \ a$ in $\lal$.

First suppose that $a \tau^{-1} \in [z,a]$.
Then 

$$a \tau^{-1} \alpha \ \in \ [z \alpha, a \alpha] \ = \  [b, 
z \alpha] \ \ \Rightarrow \ \ a \tau^{-1} \alpha \ \in \ \kk
\ \ \Rightarrow \ \ a \tau^{-1} \alpha \tau \ \in \ \kk$$

\noindent
and this last point 
bridges
to $b$ in $\lal$. 
Then $a \tau^{-1} \alpha \tau \alpha^{-m} = a \gamma \beta$
bridges to $b \alpha^{-m}$ in $\lal$.
But

$$b \alpha^{-m} \ \pa \ b \ \ {\rm in} \ \lal \ \ \ \Rightarrow \ \ \ 
a \gamma \beta \ \ {\rm bridges \ to} \ \   
b \tau^{-1} \ \ {\rm in} \ \lab.$$

\noindent
On the other hand $a \gamma \in [z, a \tau^{-1}]$ and
bridges to $a$ in $\lab$, so $a \gamma \beta$ bridges to 
$a \tau^{-1} \beta$. 
Since 
$a \tau^{-1} \beta$ is a point
in $\lab - \kk$ it is not equal to $b \tau^{-1}$, leading
to a contradiction.

\vskip .05in
The second option is $a \tau^{-1} > a$ in $\kk$.
Here use 

$$z \beta^{-1} \ = \ z \alpha \tau^{-1} \ \in \ \kk, \ \ \ \
z \alpha \ \in \ [z, z \beta^{-1}) \ \ \ \Rightarrow \ \ \ 
z \beta^{-1} \alpha^{-1} \ \not \in \ \uu.$$

\noindent
On the other hand $z \alpha^{-1}$ bridges to
$a \alpha^{-1}$ in $\lal$ so bridges to $a \tau^{-1}$
in $\lab$. So $z \alpha^{-1} \beta^{-1}$ bridges
to $a \tau^{-1} \beta^{-1}$ in $\lab$ and is in
$\uu$. As above this is a contradiction.

We conclude that problem 2 does not occur.

As in problem 1, this implies that

$$\kk \alpha^{-1} \ \cap \ \kk \ = \ 
[t', r'], \ \ \ {\rm with} \ \ t' \ \not = \ r',
t' \ \not = \ z$$

\noindent
and if
$\kk$ not properly embedded on the other side then
$r'$ has to be finite in $\kk$.

Then clearly $\kk \alpha^{-1} \tau^{-1} \subset \uu$ and
intersects $\kk$ in a segment. 

The last problem is the following:

\vskip .15in
\noindent
{\bf {Problem 3}} $-$ Does $\kk \alpha \tau \alpha^{-1}
\subset \uu$?

Suppose not, that is, $\kk \alpha \tau \alpha^{-1}
\not \subset \uu$.
We have to be careful here.
First a preliminary claim:

\vskip .08in
\noindent
{\underline {Claim}} $-$ $z \in \kk \alpha \tau \alpha^{-1}$.

If this is not true then
$\kk \alpha \tau \alpha^{-1} \cap \uu = 
\emptyset$.
Notice that 

$$\kk \alpha \tau \ \cap \ \lal  \ \not = \ \emptyset
\ \ \ \Rightarrow \ \ \ 
\kk  \alpha \tau \alpha^{-1} \ \cap \ \lal \
\not = \ \emptyset \ \ \ {\rm and} \ \ \ 
\kk \alpha \tau \alpha^{-1} \
\cap \ \uu \ \not = \ \emptyset,$$

\noindent 
contrary to assumption here.

So consider $\kk \alpha \tau \cap \lal = \emptyset$.
Also here 
$\kk \alpha \tau \cap \kk$ is a non trivial
segment.
If $\kk \alpha \tau$ bridges to $a$ in $\lal$ then
$\kk \alpha \tau \alpha^{-1}$ is contained in
$\uu$ and we are done.
If follows that $\kk \alpha \tau$ has to bridge
to $b$ in $\lal$ and hence $z \alpha$ has to be in
the this bridge. But then $z \alpha$ is in $\kk$,
which was disallowed in problem 2.
This proves the claim.

\vskip .05in
We now analyse what happens when 

$$z \ \in \ \kk \alpha \tau \alpha^{-1} \ \ \ {\rm so} \ \ 
z \tau^{-1} \ = \ z \in \kk \alpha \beta
\ \ \ {\rm and} \ \ \ z \beta^{-1} \alpha^{-1} \ \in \ \kk.$$

\noindent
Also $z \beta^{-1} \alpha^{-1} \gamma = z \alpha^{-1} \beta^{-1}$
is in $\kk$ as well.


\begin{figure}
\centeredepsfbox{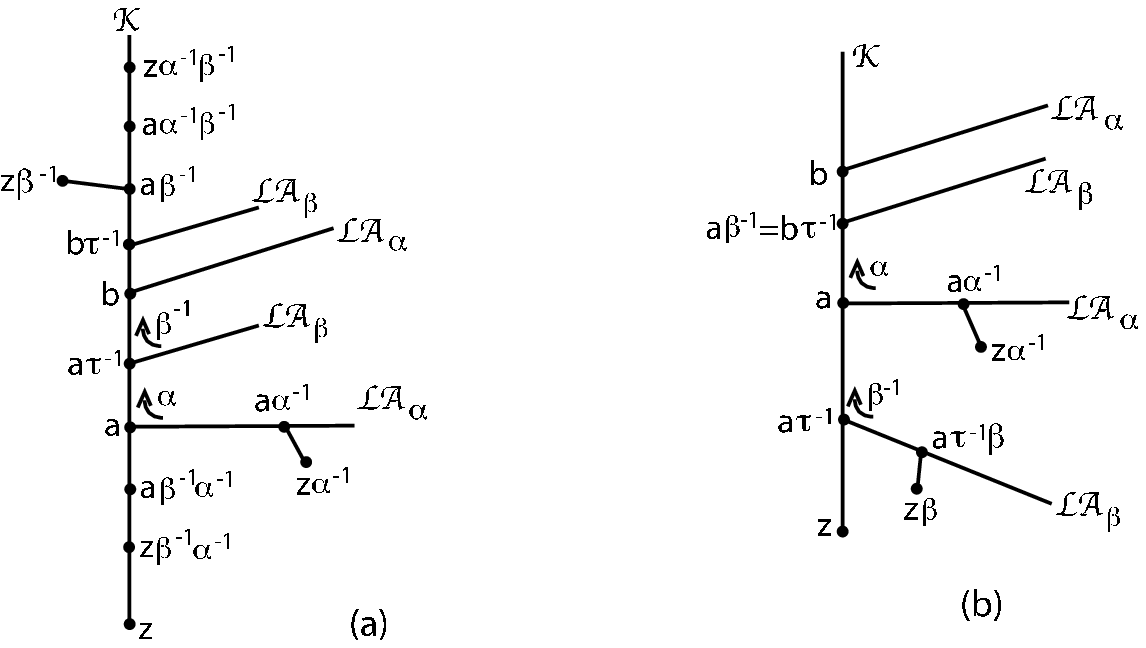}
\caption{
Situation \ $a \alpha^{-1} \ \pa \ a$ \ in $\lal$:
a. Picture when $a \tau < a$ in $\kk$,
b. Picture when $a \tau^{-1} < a$ in $\kk$.}
\label{13}
\end{figure}

\vskip .2in
\noindent
{\bf {Situation I}} $-$ \ $a \alpha^{-1} \ \pa \ a$ in $\lal$.

\vskip .2in
\noindent
{\bf {Situation I.1}} $-$ $a \tau  < a$ in $\kk$.

Here $z \alpha^{-1}$ bridges to $a \alpha^{-1}$ in
$\lal$, so it bridges to $a \tau^{-1}$ in $\lab$. 
Also

$$z \alpha^{-1} \beta^{-1} \ \in \ \kk \ \ \  {\rm and} \ \ \ 
a \tau^{-1} \ \prec \ a \ \prec \ 
a \alpha^{-1} \ \prec \ z \alpha^{-1}.$$

\noindent
As
$\beta^{-1}$ moves 
points up along $\kk$, it follows that
$z \alpha^{-1} \beta^{-1} > b$ in $\kk$ and
$a \tau^{-1} \beta^{-1} = b \tau^{-1}$.
Here $a \alpha^{-1} \in
[a \tau^{-1}, z \alpha^{-1}]$, see fig. \ref{13}, a.
Then 

$$a \tau^{-1} \beta^{-1} \ = \ b \tau^{-1} \ \prec \
a \beta^{-1} \ \prec \ a \alpha^{-1} \beta^{-1} \ = \ 
v_1 \ \prec \ z \alpha^{-1} \beta^{-1} \ = \ v_2$$

\noindent
and all are in $\kk$.
Also $a \beta^{-1} \in (b, \ a \alpha^{-1} \beta^{-1})
\subset \kk$ and $z \beta^{-1}$ bridges to $\kk$ in
$a \beta^{-1}$ so bridges to $\lal$ in $b$.
Then $z \beta^{-1} \alpha^{-1} 
= v_2 \gamma^{-1} \in \kk$ bridges
to $a$ in $\lal$ and $a \beta^{-1} \alpha^{-1}
= v_1 \gamma^{-1}$ is in $(z \beta^{-1} \alpha^{-1}, \ a)$,  \
see fig. \ref{13}, a.
Then

$$z \beta^{-1} \alpha^{-1} \ \prec \ a \beta^{-1} \alpha^{-1}
\ \prec \ a \alpha^{-1} \beta^{-1} \ \prec \ z \alpha^{-1} \beta^{-1}$$

\noindent
all points in $\kk$.
This contradicts the fact that $\gamma$ acts as
a translation in $\kk$.

\vskip .2in
\noindent
{\bf {Situation I.2}} $-$
$a \tau > a$ in $\kk$.

Here $z \alpha^{-1}$ bridges to $a$ in
$\kk$, see fig. \ref{13}, b.
If $a \geq b \tau^{-1}$ in $\kk$ then
$z \alpha^{-1}$ bridges to 
a point $t \geq_{\beta} b \tau^{-1}$ in $\lab$, so

$$z \alpha^{-1} \beta^{-1} \ \ \ {\rm bridges \ to} \ \ 
\lab \ \ {\rm in \ a \ point} \ 
\geq_{\beta} \ b \tau^{-1} \beta \ \ \ {\rm and} \ \ \
z \alpha^{-1} \beta^{-1} \ \not \in \ \kk,$$

\noindent
contradiction.
Hence $a < b \tau^{-1}$ in $\kk$ and 
$z \alpha^{-1}$ bridges to $a$ in $\lab$ so
$z \alpha^{-1} \beta^{-1}$ bridges to
$a \beta^{-1}$ in $\lab$ and as $z \alpha^{-1} \beta^{-1}$
is in $\kk$ then 

$$z \alpha^{-1} \beta^{-1} \ > \ b \tau^{-1} \ \ {\rm in} \ \
\kk  \ \ \ {\rm and} \ \ \  a \beta^{-1} \ = \ b \tau^{-1}
\ \ \ {\rm or} \ \ 
a \tau \alpha \ = \  b.$$

\noindent
Now

$$a \beta^{-1} \ = \ b \tau^{-1} \ \ \ {\rm so} \ \ \ 
a \alpha \ = \ a \tau^{-1} \beta^{-1} \tau \ < \ a \beta^{-1} \tau
\ = \ b$$

\noindent
so in particular 
$a \alpha$ is in $\kk$.
Also $z \beta$ bridges to $a$ in $\lal$ and so does $z$.
Hence $z \beta \alpha = z \alpha \beta$ and $z \alpha$
bridge to $a \alpha$ in $\lal$. 
Since $a \alpha < b$ then $z \alpha, z \alpha \beta$ bridge
to $a \alpha$ in $\lat$ as well.

If $a \alpha < b \tau^{-1}$ in $\kk$ then
$z \alpha, z \alpha \beta$ bridge 
to $a \alpha$ in $\lab$, impossible $-$ they have
to bridge to distinct points in $\lab$.
If 

$$b \tau^{-1} \ \in \ (a, a \alpha) \ \ \ \Rightarrow \ \ 
z \alpha, \ z \alpha \beta \ \ \ {\rm bridge \ to}  \ \ 
b \tau^{-1} \ \ {\rm in}  \ \lab,$$

\noindent
also contradiction. Therefore $a \alpha = b \tau^{-1}$
or $a \alpha \tau = b$. 
Now 

$$a \alpha \tau \alpha^{-1} \tau^{-1}
\ = \ b \alpha^{-1} \tau^{-1} \ = \ a
\ \ \ \ {\rm so} \ \ \ a \gamma = a \alpha^{-1} \beta^{-1}.$$

\noindent
Notice $a \gamma \in [z, a \tau^{-1}]$.
But $a \alpha^{-1}$ bridges to $a$ in $\lab$ so
$a \alpha^{-1} \beta^{-1}$ bridges to $a \beta^{-1}
= b \tau^{-1}$ in $\lab$ and
$a \alpha^{-1} \beta^{-1}$ cannot be
$a \gamma$, contradiction.

This finishes the analysis of situation I.

The remaining options are extremely similar and
have shortened proofs.


\vskip .2in
\noindent
{\bf {Situation II}} $-$  $a \alpha \ \pa \ a$ \ in $\lal$.

\vskip .2in
\noindent
{\bf {Situation II.1}} $-$ $a \tau^{-1} < a$ in $\kk$.

This is as situation I.1 above.
Here $z \beta^{-1}$ bridges to $a$ in $\lal$,
so $z \beta^{-1} \alpha^{-1}$ bridges to $a \alpha^{-1}$
in $\lal$ and $a \alpha^{-1} = b$.
It follows that

$$b \ \prec \ a \tau^{-1} \alpha^{-1} \ \prec \ 
a \tau^{-1} \beta^{-1} \alpha^{-1} \ \prec \
z \beta^{-1} \alpha^{-1},$$

\noindent
all points in $\kk$.

On the other hand $a \tau^{-1} \alpha^{-1} \in 
(b, \ (a \tau^{-1}) \beta^{-1} \alpha^{-1}) \subset \kk$.
The point $z \alpha^{-1}$ bridges to $(a \tau^{-1}) \alpha^{-1}$
in $\kk$. 
It follows that

$$z \alpha^{-1} \beta^{-1} \ \prec \ 
(a \tau^{-1}) \alpha^{-1} \beta^{-1} \ \prec \ 
(a \tau^{-1}) \beta^{-1} \alpha^{-1} \ \prec \ z \beta^{-1} \alpha^{-1},$$

\noindent 
all points in $\kk$.
As before this contradicts the fact that $\gamma$ acts as
a translation in $\kk$.

\vskip .2in
\noindent
{\bf {Situation II.2}} $-$ $a \tau < a$ in $\kk$.

This is very much like situation I.2.
Here $z \beta^{-1}$ bridges to $a \tau^{-1}$ in $\kk$.
If $a \tau^{-1} \geq b$ in $\kk$, 
then 

$$z \beta^{-1} \alpha^{-1} \ \ \ {\rm bridges \ to \ a \ point}
\ \ 
>_{\alpha} \  b \ \ {\rm in} \ \ \lal \ \ \ \Rightarrow \ \ \ 
z \beta^{-1} \alpha^{-1} \ \not \in \ \kk,$$

\noindent
 contradiction.
Hence 

$$a \tau^{-1} \ <  \ b \ \ {\rm in} \ \kk, \ \ \ 
z \beta^{-1} \alpha^{-1}  \ > \ b \ \ {\rm in } \  \kk
\ \ \ {\rm and} \ \ \ 
a \tau^{-1} \alpha^{-1} \ = \ b \ \ {\rm or} \ \ 
a \ = \ b \alpha \tau.$$

\noindent
In addition 

$$z \alpha, z \ \  {\rm bridge \ to} \ \lab \ \ {\rm in} \ 
a \tau^{-1}  \  \ \ \Rightarrow \  \ \
z \beta \alpha  =  
z \alpha \beta,  \ z \beta \ \ {\rm bridge \ to } \ 
\lab \ \ {\rm in} \ a \tau^{-1} \beta$$

\noindent
and similarly to situation I.2, this implies
$a \tau^{-1} \beta = b$ or
$ a  = b \tau \alpha$. Then $b \alpha \beta = b$
and $b \gamma =  b \alpha^{-1} \beta^{-1}$.
But $b \gamma \geq b \tau^{-1}$ in $\kk$ and
$b \alpha^{-1}$ bridges to $b$ in $\lab$, so
$b \alpha^{-1} \beta^{-1}$ bridges to $b \beta^{-1}
= a \tau^{-1}$ and cannot be equal to $b \tau^{-1}$.

This contradiction shows that problem 3 cannot occur.
This finishes the proof of lemma \ref{probl}.
\end{proof}

It follows from lemma \ref{probl} that 
$\kk \alpha \tau \alpha^{-1} \subset \uu$,
so $\kk \alpha \beta \subset \uu$ as
is $\kk \gamma \beta \alpha$. So all of the
sets $\kk, \ \kk \alpha, \ \kk \alpha \tau$, 
\ $\kk \alpha \tau \alpha^{-1}, \
\kk \alpha \beta, \ \kk \alpha^{-1}, \  \kk \alpha^{-1} \tau^{-1}$
\ and $\kk \alpha^{-1} \tau^{-1} \alpha$ are contained in
$\uu$ and none has $z$ as an ideal point.
If $\kk$ has another ideal point $v$, then
$v$ has the same properties as $z$ and
the same situation
occurs with respect to this other ideal point.

Given these facts, an analysis exactly as in
case A.2 can be applied here.
That analysis then shows that case B.1.3 is not possible.


Hence case B.1.4 is disallowed.
This also finishes the proof of case B.1.

\vskip .05in
For case B.2 we return to the study of $\alpha$
acting freely using the axis $\ap$.

\vskip .3in
\noindent
{\bf  {Case B.2}} $-$ $Fix(\tau) \cap \ap \not = \emptyset$.

This is the key case of the proof for essential laminations.
In this case the topology will be important, in particular, the
exact condition $|p - 2q| = 1$ will be used in a crucial manner.
Let $z \in Fix(\tau) \cap \ap$.
Let $\uu_1$ (respectively $\uu_2$) 
be the component of $T - \{ z \}$ containing
$z \alpha$ (respectively $z \alpha^{-1}$).
A priori we do not know whether $z$ is also a fixed point
of $\gamma$.
In some subcases, the tricky part will be in fact to
show that $z \gamma = z$.

\vskip .1in
\noindent
{\bf {Case B.2.1}} $-$ $\uu_1 \tau = \uu_1$.

Notice that $\uu_1 \alpha$ is  contained in $\uu_1$.
Here use 
$z \alpha \tau = z \tau \gamma \beta \alpha^m =
z \gamma \beta \alpha^m$.

$$z \alpha \in \uu_1 \ \Rightarrow \ z \alpha \tau \in \uu_1
\ \Rightarrow \ z \alpha \tau \alpha^{-m} \in \uu_1 \alpha^{-m} \subset \uu_1
\ \Rightarrow \ z \gamma \beta \in \uu_1.$$

\begin{figure}
\centeredepsfbox{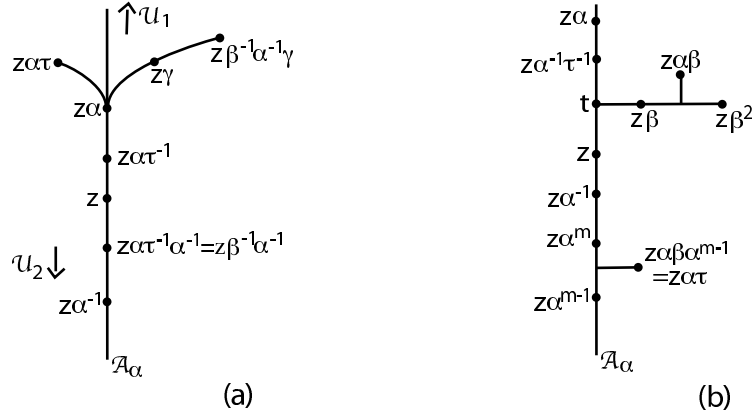}
\caption{
Case B: a. Picture when
$\uu_1 \tau = \uu_1$,
b. Picture when $\uu_1 \tau^{-1} = \uu_2$
and
$[z, z \beta] \cap [z, z \alpha] =
[z,t]$
.}
\label{14}
\end{figure}

\noindent
So $z \gamma \tau \alpha^{-1} \tau^{-1}$ is in $\uu_1$ and
then  $z \gamma \alpha^{-1} $ is in $\uu_1$ or
$z \gamma$ is in $\uu_1 \alpha$. In particular
$z \prec z \alpha \prec z \gamma$, see fig. \ref{14}, a.
We stress that 
in this case $z \gamma$ is not equal to $z$!

Use now 
$z \alpha \tau = z \alpha \beta \alpha^{m-1} =
z \alpha \tau \alpha^{-1} \tau^{-1} \alpha^{m-1}$.

$$z \alpha \tau \alpha^{1-m} \in \uu_1 \ \Rightarrow \
z \alpha \tau \alpha^{-1} \tau^{-1} \in \uu_1 \ \Rightarrow \
z \alpha \tau \alpha^{-1} \in \uu_1
\ \Rightarrow \ z \alpha \tau \in \uu_1 \alpha.$$

\noindent
In particular $z \prec z \alpha  \prec z \alpha \tau$
and $z \prec z \alpha \tau^{-1} \prec z \alpha$
and so $z \alpha \tau^{-1} \alpha^{-1} \in (z \alpha^{-1}, z)$.
In other words

$$z \alpha \tau^{-1} \alpha^{-1} \ = \ 
z \tau \alpha \tau^{-1} \alpha^{-1} \ = \ 
z \beta^{-1} \alpha^{-1} \ \ \ \in \ (z \alpha^{-1}, z).$$

\noindent
Then $z \beta^{-1} \alpha^{-1}$ is in $\uu_2$ so
$z \beta^{-1} \alpha^{-1} \gamma$ is in $\uu_2  \gamma$.
Notice $z \beta^{-1} = z \tau \alpha \tau^{-1} = z \alpha \tau^{-1}$
with $z \alpha \in \uu_1$, $z \alpha \tau^{-1}$ also in $\uu_1$.


Recall that $z \gamma \not = z$.
If  $\uu_1 \gamma \subset \uu_1$ this implies 
that $z$ is in a local axis for $\gamma$ contradicting
$z \gamma^q = z \tau^{-p} = z$. 
Therefore $\uu_1 \gamma$ is not contained in $\uu_1$ and
consequently $\uu_2 \gamma$ is contained in $\uu_1$
and so $z  \gamma$ separates $\uu_2 \gamma$  from $z$.
Hence

$$z  \alpha \ \ \ {\rm separates} \ \  \uu_2 \gamma
\ \ {\rm from} \   z \ \ \ \ {\rm and} \ \ \ 
z \beta^{-1} \alpha^{-1} \gamma \ \in \ \uu_2 \gamma.$$

\noindent
But $z \beta^{-1}  \alpha^{-1} \gamma = z \alpha^{-1} \beta^{-1}
= z \alpha^{-1} \tau \alpha \tau^{-1}$.
Now $z \alpha$  separates $z$ from 
$z \alpha^{-1} \tau \alpha \tau^{-1}$
which is in $\uu_2 \gamma$. Apply $\tau$:
$z \alpha \tau$ separates  $z$ from $z \alpha^{-1} \tau \alpha$.
Then

$$z \alpha \tau \ \in \ \uu_1 \alpha \ \ \Rightarrow \ \ 
z \alpha^{-1} \tau \alpha \ \in \ \uu_1 \alpha \ \ \Rightarrow \ \ 
z \alpha^{-1} \tau \ \in \ \uu_1 \ \ {\rm and} \ \ 
z \alpha^{-1} \ \in \ \uu_1 \tau^{-1} = \uu_1.$$

\noindent
But this contradicts $z \alpha^{-1}$ is in $\uu_2$.
This is an impossible case.

We conclude that $\uu_1 \tau \not = \uu_1$.

\vskip .2in
\noindent
{\bf {Case B.2.2}} $-$ $\uu_1 \tau \not = \uu_2$.

Then $z \alpha \tau$ is not in $\uu_2$, which 
implies $z \alpha \tau \alpha^{1-m}$ is in $\uu_1$,
or $z \alpha \beta \in \uu_1$
and $z \alpha \tau \alpha^{-1} \tau^{-1}$ is
in $\uu_1$.
By assumption $z \alpha \tau \not \in \uu_1$, hence
$z \alpha \tau \alpha^{-1} \in \uu_2$ and
$z \alpha \tau \alpha^{-1} \tau^{-1} \in \uu_2 \tau^{-1}$.
This would imply $\uu_2 \tau^{-1} = \uu_1$
or $\uu_1 \tau = \uu_2$, so the assumption is incompatible.

We conclude that $\uu_1 \tau = \uu_2$.

\vskip .2in
\noindent
{\bf {Case B.2.3}} $-$ $\uu_1 \tau^{-1} = \uu_2$.

This is a very interesting case. Here we only use
the fact that $p$ is odd.

First consider $z \beta = z \tau \alpha^{-1} \tau^{-1}
= z \alpha^{-1} \tau^{-1}$ which is in
$\uu_2 \tau^{-1} = \uu_1$.
Then $z \alpha, z \beta$ are in the component $\uu_1$,
hence $[z, z \alpha]$, $[z,  z \beta]$ share
a subprong.
Suppose first that

$$[z, z \beta] \cap [z, z \alpha] =
[z,t],  \ \ t \not = z \alpha, z \beta,
\ \ \ {\rm that \ is} \ \ \ z \alpha \not \in [z, z \beta],
 \  \ z \beta \not \in [z, z \alpha]$$

\noindent
see fig. \ref{14}, b.
Then $z \alpha \beta$ bridges to $t$ in $\ap$
and  $z \alpha \beta \alpha^{m-1}$  bridges
to $\ap$ in $t \alpha^{m-1}$ which is
a point in $(z \alpha^{m}, z \alpha^{m-1})$.
But 

$$z \alpha \beta \alpha^{m-1} = z \alpha \tau \ \ \ \Rightarrow \ \ \ 
z \alpha^{-1} \tau^{-1} \ \in \ 
[z, z \alpha) \ \ \ \Rightarrow \ \ \ 
z \beta = z \alpha^{-1} \tau^{-1} \ \in \ [z, z \alpha),$$

\noindent
contradiction.

So either $z \beta \in [z, z \alpha]$ or
$z \alpha \in [z, z \beta]$.

\vskip .2in
\noindent
{\bf {Situation I}} $-$ $z \alpha$ is in $[z, z \beta]$.

Use $z \beta \tau = z \tau \alpha^{-1} = z \alpha^{-1}$.
As $z \alpha$ is in $[z, z \beta]$, then
$z \alpha \tau \in [z, z \beta \tau] = [z, z \alpha^{-1}]$
and $z \alpha \tau \alpha^{1-m} \in [z \alpha^{-m}, z \alpha^{1-m}]$.
But 

$$z \alpha \tau \alpha^{1-m} \ = \ 
z \tau^{-1} \alpha \tau \alpha^{1-m} \ = \ z \alpha \beta,
\ \ \ \ \ {\rm so} \ \ \ z \alpha \beta \in [z \alpha^{-m}, z \alpha^{1-m}]
\subset \ap.$$

\noindent
We stress that $z \alpha \beta \in \ap$.
Here $z \beta^{-1} \prec z \prec z  \alpha$,
hence $z \prec z \beta \prec z \alpha \beta$.
It follows that 

$$z \beta \ \in \ \ap \ \ \ {\rm and} \ \ \
z \beta \ \in \ [z, z \alpha \beta] \ \ \ \Rightarrow \ \ \ 
z \alpha \beta \alpha^{-1} \ \in \
[z \alpha^{m-1}, z \alpha^{-m}].$$

We want $z \gamma = z$ or $z \alpha \beta = z \beta \alpha$.
We first analyse the other two possibilities.

\vskip .2in
\noindent
{\bf {Situation I.1}} $-$ $z \alpha \beta \alpha^{-1} > z \beta$ 
in $\ap$.


Then $z \beta \prec z \alpha \beta \alpha^{-1} \prec z \alpha \beta$,
so $z \prec z \gamma \prec z \alpha$, or $z \gamma \in (z, z \alpha)$,
so $z \gamma \in \uu_1$.
Clearly $z \beta \alpha \in \ap$.
Here $z \alpha \beta > z \beta \alpha$ in $\ap$.
Then

$$z \prec z \beta \alpha \prec z \alpha \beta
\ \ \ {\rm all \ in } \ \ \ap \ \ 
\Rightarrow \ z \beta^{-1} \prec z \beta \alpha \beta^{-1}
\prec z \alpha \ \ \ {\rm and} \ \ 
z \beta^{-1} \alpha^{-1} \prec z \gamma^{-1} \prec z.$$

\noindent
But $z \beta^{-1} = z \alpha \tau^{-1} \in \uu_2$,
hence $z \beta^{-1} \alpha^{-1}$ is in $\uu_2$.
%
Now $z \gamma \in \uu_1, \ z \gamma^{-1} \in \uu_2$,
therefore $z$ is in a local axis for $\gamma$, hence
$z \gamma^q \not = z$, contradiction.

\vskip .2in
\noindent
{\bf {Situation I.2}} $-$ Suppose  \ $z \alpha \beta \ 
\pa  \ z \beta \alpha$.

Then 

$$z \ \prec \ z \alpha \beta \alpha^{-1} \ \prec \
z \beta \ \ \ \ \Rightarrow \ \ \ \ 
z \beta^{-1} \ \prec \ z \gamma \ \prec \ z.$$

\noindent
As $z \beta^{-1} = z \alpha \tau^{-1}$ is in $\uu_2$,
then $z \gamma$ is in $\uu_2$.


Now \ $z \alpha \beta \ \pa  \ z \beta \alpha$. \
If $\ab$ contains elements in $\ap$ above 
$z \alpha \beta$, that is,
$\ab \cap \ap \supset [z,t)$ with
$t >_{\alpha}  z \alpha \beta$.
Then 

$$z \ \prec \ z \alpha \ \prec \ t \beta^{-1} \ \prec \ 
z \beta \alpha \beta^{-1}, \ \ \ {\rm with} \ \
t \beta^{-1} \ \in \ \ap \ \ \ 
\Rightarrow \ \ \ 
z \ \prec \ t \beta^{-1} \alpha^{-1}
\ \prec \ z \gamma^{-1}$$

\noindent
 with $t \beta^{-1} \alpha^{-1}$ in
$\ap$ so $z \gamma^{-1}$ is in $\uu_1$ and not in $\uu_2$.


On the other hand if
$\ab$ escapes $\ap$ in $z \alpha \beta$, then
$z \beta \alpha \beta^{-1}$ bridges to
$\ab$ in $z \alpha$,
hence bridges to $\ap$ in $z \alpha$ as $z \alpha \in (z, z \beta)$. 
Hence $z \beta \alpha \beta^{-1} \not \in \uu_2 \alpha$
and $z \beta \alpha \beta^{-1} \alpha^{-1} = z \gamma^{-1}$
bridges to $\ap$ in $z$ and $z \gamma^{-1}$ is not
in $\uu_2$.
In any case $z \gamma^{-1}$ is not in $\uu_2$
and $z \gamma$ is in $\uu_2$ so
$z$ separates $z \gamma$ from $z \gamma^{-1}$ and
$z$ is in a local axis for $\gamma$, impossible.

We conclude that $z \alpha \beta = z \beta \alpha$ 
or that $z \gamma = z$.

\vskip .2in
\noindent
{\bf {Situation I.3}} $-$ $z \gamma = z$.

Then $\gamma$ leaves invariant the set of components
of $T - \{ z \}$.
Recall that $\uu_1 \tau^{-1} = \uu_2$
and $\uu_1 \tau = \uu_2$ in situation I.
Use 
$z \beta^{-1} \alpha^{-1} \gamma \ = 
\ z \alpha^{-1} \beta^{-1}$.
The left side is
$z \tau \alpha \tau^{-1} \alpha^{-1} \gamma =
z \alpha \tau^{-1} \alpha^{-1} \gamma$.

$$z \alpha \in \uu_1 \ \Rightarrow \ 
z \alpha \tau^{-1} \in \uu_1 \tau^{-1} \not = \uu_1,
\ \ \ {\rm so} \ \ \ z \alpha \tau^{-1} \alpha^{-1}
\in \uu_2 \ \ \ {\rm and} \ \ \
z \alpha \tau^{-1} \alpha^{-1} \gamma \in \uu_2 \gamma.$$

\noindent
On the other hand the right side is
$z \alpha^{-1} \tau \alpha \tau^{-1}$:

$$z \alpha^{-1} \in \uu_2 
\ \Rightarrow \ z \alpha^{-1} \tau \in \uu_2 \tau = \uu_1,
\ \ z \alpha^{-1} \tau \alpha \in \uu_1
\ \ \ {\rm and} \ \ \ z \alpha^{-1} \tau \alpha \tau^{-1}
\in \uu_1 \tau^{-1} = \uu_2.$$

\noindent
So $\uu_2 \gamma \cap \uu_2 \not = \emptyset$.
Since $\gamma$ now preserves the set of components
of $T - \{ z \}$ it follows that
$\uu_2 \gamma = \uu_2$
and 
$ \uu_1 \gamma \ = \ \uu_2 \tau \gamma \ = \
\uu_2 \gamma \tau \ = \ \uu_2 \tau \ = \ \uu_1$.
Now we use $p$ odd and $\tau^p \gamma^q = id$:

$$\uu_1 \ = \ \uu_1 \gamma^q \tau^p \ = \
\uu_1 \tau^p \ = \ \uu_1 \tau^{p(mod 2)}
\ = \ \uu_1 \tau.$$

\noindent
This contradicts $\uu_1 \tau \not = \uu_1$ and
finishes the analysis of situation $I$.

\vskip .2in
\noindent
{\bf {Situation II}} $-$ $z \beta \in [z, z \alpha]$.

This is very similar to the previous case if we think
of it in the appropriate way. The trick here is to
switch the roles of $\alpha$ and $\beta$, which
can be done.
Notice first that $z \beta \in \uu_1$ and
$z \beta^{-1} = z \tau \alpha \tau^{-1}
= z \alpha \tau^{-1}$ is in $\uu_2$.
So the component of $T - \{ z \}$ containing
$z \beta$ (respectively $z \beta^{-1}$)
is the $\uu_1$ (respectively $\uu_2$).
First rewrite the relations as

$$ \tau \alpha \tau^{-1} \ = \ \beta^{-1}
\ \ \ \ \ \ \ \tau \beta \tau^{-1} \ = \
\gamma^{-1} \alpha \beta^{m} \ = \ 
\beta \alpha \beta^{m-1}$$

\noindent
As $z \beta$ is in $[z, z \alpha]$
then $z \beta \tau^{-1}$ is in
$[z \tau^{-1}, z \alpha \tau^{-1}] = 
[z, z \beta^{-1}]$.
So 

$$z \tau \beta \tau^{-1} \beta^{1-m} \ = \
z \beta \tau^{-1} \beta^{1-m} \ = \ z \beta \alpha
\ \ \in \ \ [z \beta^{-m}, z \beta^{1-m}] \ \subset
\ \ab.$$

\noindent
As $z \beta \in [z, z \alpha]$, then $z \beta \alpha$ is
in $[z \alpha, z \alpha^2]$
and 

$$z \alpha \in [z, z \beta \alpha] \ \subset \ 
[z, z \beta^{1-m}] \ \subset \ \ab.$$

\noindent
Therefore $z \alpha$ is in $\ab$ and similarly
$z \alpha \beta$, $z \beta \alpha$ are in $\ab$. 

From this point on the proof is entirely similar
to the analysis in situation I:
consider whether $z \alpha \beta \ \pb \ z \beta \alpha$, \ \
$z \alpha \beta \ >_{\beta}  \ z \beta \alpha$, \
or $z \alpha \beta = z \beta \alpha$,
with completely analogous proofs. 

Therefore this case is disallowed. This finishes the
analysis of case B.2.3, \ $\uu_2 \tau = \uu_1$.

\vskip .2in
\noindent
{\bf {Case B.2.4}} $-$ $\uu_1 \tau = \uu_2$,
\ $\uu_1 \tau^{-1} \not = \uu_2$.

This is the most interesting case which relates
to the topology in a crucial way.


Use $z \beta^{-1} \alpha^{-1} \gamma = z \alpha^{-1}
\beta^{-1}$.
The right side is $z \tau \alpha \tau^{-1} \alpha^{-1} \gamma
= z \alpha \tau^{-1} \alpha^{-1} \gamma$.

$$ z \alpha \in \uu_1 \ \ \Rightarrow \   z \alpha \tau^{-1}
\in \uu_1 \tau^{-1} \not = \uu_1 \ \Rightarrow
\ z \alpha \tau^{-1} \alpha^{-1} \in \uu_2.$$

\noindent
Hence $z \beta^{-1} \alpha^{-1} \gamma$ is
in $\uu_2 \gamma$.
On the other hand $z \alpha^{-1} \beta^{-1} =
z \alpha^{-1} \tau \alpha \tau^{-1}$:

$$z \alpha^{-1} \tau \ \in \uu_2 \tau \not = \uu_2
\ \Rightarrow z \alpha^{-1} \tau \alpha \in \uu_1 \ \
\Rightarrow \ 
z \alpha^{-1} \tau \alpha \tau^{-1} \in \uu_1 \tau^{-1}
\ \not = \ \uu_2.$$

\noindent
We conclude that

$$\uu_2 \gamma \ \cap \ \uu_1 \tau^{-1} \ \not = \emptyset,
\ \ \ \ {\rm or} \ \ \ \ 
\uu_1 \tau \gamma \cap \uu_1 \tau^{-1} \not = \emptyset
\ \ \ \ \ \ \ (*).$$

\noindent
What we actually want is that these two sets
are equal. A priori we have to be careful because
$\gamma$ may not preserve the set of components
of $T - \{ z \}$, or equivalently we may have
$z \gamma \not = z$. So we first deal with this case.
We will need the following useful lemma:

\begin{lemma}{}{}
Let $\eta$ be a homeomorphism of a tree $V$ so that
$\eta^l$ has a fixed point $a$, where $l$ is not $0$.
Then there is a fixed point of $\eta$ in $[a, a \eta]$.
\end{lemma}

\begin{proof}{}
Consider $a \eta^2$.
If $a \eta^2$ is in $[a, a \eta]$ and
not equal to $a \eta$, then $\eta$ sends
$[a, a \eta]$ into itself and has a fixed
point there, done.
If $a \eta$ is in $(a, a \eta^2)$ then $a$ is in 
a local axis of $\eta$ and $a \eta^l$ is not $a$,
impossible.
If $a$ is in $(a \eta, a \eta^2)$, then $\eta^{-1}$ sends
$[a \eta, a \eta^2]$ into itself (into $[a, a \eta]$)
producing a fixed point there, done.


\begin{figure}
\centeredepsfbox{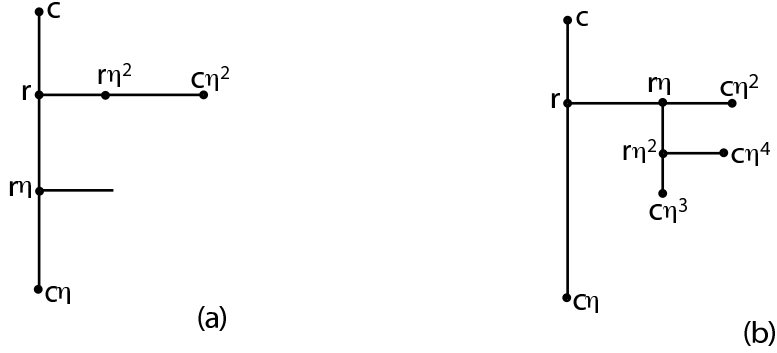}
\caption{
a.  $r \eta \in [r, a \eta]$,
b. $r \eta \in (r, a \eta^2]$.}
\label{15}
\end{figure}

We can now assume $a \eta^2$ bridges to $[a, a \eta]$
in a point
$r$ which is in $(a, a \eta)$, see fig. \ref{15}, a.
If $r \eta = r$ we are done.
Assume $r \eta \not = r$.
Then $r \eta$ is in $[a \eta, a \eta^2]$.

Suppose first that 
$r \eta$ is in $[r, a \eta]$, see fig. \ref{15}, a.
Then $r \eta^2$ is in $[r \eta, a \eta^{2}]$ so either
$[r \eta, r]$ is contained in its image 
under $\eta$ or vice versa.
In any case there is a fixed point of 
$\eta$ in $[r, r \eta]$.

Suppose now that $r \eta$ is in $(r, a \eta^2]$ see
fig. \ref{15}, b.
Hence $a \prec r \prec r \eta$ and $a \eta \prec r \eta \prec
r \eta^2$.
Then $r \in (a \eta, r \eta)$ and 
$r \eta \in (r, r \eta^2)$, so $r$ is in a local axis for
$\eta$. 
This implies
that $a \eta^t \not = a$ for any nonzero $t$ in ${\bf Z}$,
contradiction.
This finishes the proof.
\end{proof}

We are back to case B.2.4.

\vskip .1in
\noindent
{\bf {Situation I}} $-$
$z \gamma \not = z$.

Suppose first that $z \gamma \in \uu_2$.
Notice $\uu_2 \tau \not = \uu_1$ and also $\not = \uu_2$.
Then there is $c$ in $[z, z \gamma]$ fixed
by $\gamma$ so $c$ is in $\uu_2$. 
This implies 

$$\uu_2 \tau \gamma \ \subset \ 
\uu_2 \ \ \ \Rightarrow \ \ \ 
\uu_1 \tau^2 \gamma \ \subset \ \uu_2, \ \ {\rm or} \ \ 
\uu_1 \tau \gamma \ \subset \ \uu_1.$$

\noindent
But by $(*)$
$\uu_1 \tau \gamma \cap \uu_1 \tau^{-1} \not = \emptyset$,
which now implies 
$\uu_1 \tau^{-1} \cap \uu_1 \not = \emptyset$.
This is impossible and rules out this case.

The second possibility is that $z \gamma \in \uu_1$.
Here $\uu_2 \gamma \subset \uu_1$ so
$\uu_1 \tau \gamma \subset \uu_1$.
As $\uu_1 \tau \gamma \cap \uu_1 \tau^{-1} \not
= \emptyset$ then $\uu_1 \tau^{-1} \cap \uu_1 \not
= \emptyset$, also impossible.

The final option is
$z \gamma \not \in \uu_1 \cap \uu_2$, 
$z \gamma \in \uu_3$ (which may be $\uu_2 \tau$ or not).
Here there is $y$ fixed by $\gamma$ with
$y \in \uu_3$. Here first use

\noindent
$$\uu_2 \gamma \ \subset \ \uu_3, \ \ {\rm or} \ \ 
\uu_1 \tau \gamma \ \subset \ \uu_3 \ \ \ \Rightarrow \ \ \ 
\uu_1 \tau^{-1}  \ \cap \ \uu_3 \ \not = \ 
\emptyset \ \ {\rm and} \ \ 
\uu_1 \tau^{-1} \ = \ \uu_3.$$

\noindent
Use $\uu_1 \gamma \subset \uu_3$, so

$$\uu_1 \tau \gamma \subset \uu_3 \tau
\ \ \ {\rm and} \ \ \ 
\uu_1 \tau^{-1} \cap \uu_3 \tau \not = \emptyset
\ \ \ {\rm or} \ \ \ \uu_1 \tau^{-1} = \uu_3 \tau.$$

\noindent
Then $\uu_3 = \uu_3 \tau$ or $\uu_1 \tau^{-1} = \uu_1 \tau^{-2}$,
so $\uu_1 \tau = \uu_1$, impossible.
This rules out this final option.

We conclude that:

\vskip .1in
\noindent
{\bf {Situation II}} $-$ $z \gamma = z$.

This is a crucial case. In fact there is an essential
lamination in $M_{p/q}$ whenever $|p - 2q| \geq 2$ and this
essential lamination may satisfy these properties:
$\tau$ has a fixed point, $\alpha$ has an axis (or at
least a local axis) which contains the fixed point of
$\tau$. See more below.
So here is a part of the proof where the
specific condition $|p - 2q| = 1$ needs to be used. See remark
below on the topological significance of this condition.

Here is the proof.
Since $z \gamma = z$ , then $\gamma$ permutes components
of $T - \{ z \}$. Since $\uu_1 \tau \gamma \cap
\uu_1 \tau^{-1} \not = \emptyset$, it now follows
that 

$$\uu_1 \tau \gamma  \ = \ \uu_1 \tau^{-1}
\ \ \ \ \ {\rm or} \ \ \ \ 
\uu_1 \gamma \tau^2 \ = \ \uu_1.$$

\noindent
We now compute

$$\uu_1 \ = \ \uu_1 \tau^p \gamma^q \ = \ 
\uu_1 \tau^{p-2q} \tau^{2q} \gamma^q \ = \
\uu_1 (\gamma \tau^2)^q \tau^{p-2q} \ = \
\uu_1 \tau^{p-2q}.$$

\noindent
When $|p - 2q| = 1$ then either
$\uu_1 = \uu_1 \tau$ or $\uu_1 = \uu_1 \tau^{-1}$ $-$
so in either case $\uu_1 = \uu_1 \tau$!
But this contradicts that we proved before
that in case B, $\uu_1 \tau$ is not equal to $\uu_1$.
This is a contradiction showing that case B.2.4 cannot happen.
This is quite straightforward, but it needed all
the previous steps.

This finishes the proof of case B: $Fix(\tau) \not = \emptyset$,
$Fix(\alpha) = \emptyset$.

\vskip .2in
\noindent
{\bf {Remark}} $-$
We now analyse the topology of this situation.
Consider the original stable foliation in the
torus bundle over the circle (the manifold $M$). This
produces a lamination $\lambda_1$ in $M - N(\delta)$.
The solid torus complementary component of $\lambda_1$
have degeneracy locus $(1,2)$ that is $\gamma \tau^2$.
This means the $\gamma \tau^2$ is a curve in the boundary
leaf of the complementary component and it also
preserves the ``outer" side of this complementary
component. Now do $p/q$ Dehn filling on $M - N(\delta)$
and look at the tree $T$ produced. The leaf
through $\delta$ collapses to a fixed point $z$ of $\tau$
(and $\gamma$ too). Usually neither $\tau$ nor $\gamma$
preserves the complementary components of $z$, but
the above fact about the degeneracy locus means
that $\gamma \tau^2$ does preserve these components $-$
if $\uu_1$ is one such component of $T - \{ z \}$
then $\uu_1 \gamma \tau^2 = \uu_1$
After $(q,p)$ Dehn surgery,
the leaf space $T$ of the lamination has a singularity
at $z$ with exactly
$|p - 2q|$ prongs. 
The transformation $\tau$ rotates by one in the set of
prongs, hence
$\tau^{p - 2q}$ preserves each of
the prongs.
This is also detected by $\gamma \tau^2$ preserving the
set of prongs and $\tau^p \gamma^q$ being null homotopic.
All is well when $|p - 2q| \geq 2$, because we have 
2 or more prongs and the lamination is essential and
the action is very nice.
However when $|p - 2q| = 1$ there is only
one prong and the lamination is not essential.
It is amazing that this sort of difficulty 
can still be detected on
the level of group action on trees.
Notice that this is exactly what the proof shows that
$\uu_1 \tau = \uu_1$, which must happen if there is
only one prong.

\vskip .2in

\section{Case C $-$ $\alpha$ has a fixed point
and $\tau$ has a fixed point}

Let $s$ in $Fix(\kappa)$, $w$ in $Fix(\alpha)$ with
$(s,w] \cap Fix(\kappa) = \emptyset$ and
$[s,w) \cap Fix(\alpha) = \emptyset$.
The following notation will be very useful in
this section. Given $u \not = v$ in $T$ let

$$T_u(v) \ = \ \{ {\rm component \ of} \ \  T - \{ u \}
\ \ \ {\rm containing} \ \ v \ \}.$$

\noindent
Let $\we = T_s(w)$, $\vv = T_w(s)$.
First in this section we will try to prove that
$\we$ is invariant under $\tau$ and $\vv$ is
invariant under $\alpha$. This will produce local
axes for $\alpha$ and (eventually) for $\tau$ and we will
see how the 2 axes interact.

\vskip .2in
\noindent
{\bf {Case C.1}} $-$ Suppose $\we \tau \not = \we$.

\vskip .2in
\noindent
{\bf {Case C.1.1}} $-$ Suppose $w \in [s, s \alpha]$.

This is equivalent to $\vv \alpha \not = \vv$.
Notice $s \alpha \not = w$. We know $s \alpha \beta = s \beta \alpha$,
and $s \beta \alpha = s \alpha^{-1} \tau^{-1} \alpha$.
Then 

$$s \alpha^{-1} \not \in \vv \ \ \Rightarrow \ \ 
s \alpha^{-1} \in \we \ \ \Rightarrow \ \ 
s \alpha^{-1} \tau^{-1}  \in  
\we \tau^{-1} \subset \vv \ \ \Rightarrow \ \ 
s \alpha^{-1} \tau^{-1} \alpha \in \vv \alpha \subset \we.$$

On the other hand
$s \alpha \beta = s \alpha \tau \alpha^{-1} \tau^{-1}$.
Here 

$$s \alpha\ \in \ \vv \alpha \ \subset \ \we \ \ \ \Rightarrow \ \ \
s \alpha \tau \ \in \ \we \tau \ \subset \ \vv \ \ \ \Rightarrow \ \ \ 
s \alpha \tau^{-1} \alpha^{-1} \ \in \ 
\vv \alpha^{-1} \ \subset \ \we$$

\noindent
and $s \alpha \beta$ is
in $\we \tau$.
These two facts together imply $\we = \we \tau$,
contrary to assumption.

Conclusion: if $\we \tau \not = \we$, then $\vv \alpha = \vv$.

\vskip .2in
\noindent
{\bf {Case C.1.2}} $-$ $s \alpha^{-1} \not \in [s,w]$,
\ $s \alpha \not \in [s, w]$.

This implies 
$s \alpha, s \alpha^{-1}$ are in $\we$.
For otherwise if $s \alpha$ is not in $\we$, then
$s$ is in $(w,s \alpha]$ and
so $s \alpha^{-1}$ is in $[w,s]$.

In this case $s \alpha^{-1}$ bridges to $[s,w]$ in
a point $r$ with $r \in (s,w)$ $-$ the important
fact
is that
$r$ is not one of the endpoints which would occur
if $s \alpha^{-1}$ is not in $\we$ or $\vv$.
Then 

$$r \ \in \ [w, s] \ \cap \ [w, s \alpha^{-1}] \ \ \ \Rightarrow \ \ \ 
r \alpha^{-1} \ \in \
[s, w \alpha^{-1}].$$

\noindent
Notice $r \alpha^{-1}$ is not equal to $r$.
If $r \alpha^{-1}$ is in $(r, s \alpha^{-1})$,
then $s \alpha^{-2}$ bridges to $[r, s \alpha^{-1}]$
in $r \alpha^{-1}$, hence $s \alpha^{-2}$ bridges
to $[s,w]$ in $r$. The same happens for all
$s \alpha^n$ with $n$ negative.
If on the other hand $r \alpha^{-1}$ is in $(w,r)$
then $s \alpha^{-2}$ bridges to $r \alpha^{-1}$
in $[s,w]$ and $s \alpha^n$ bridges to $[s,w]$
in $r \alpha^{n+1}$ for all $n$ negative. Notice
then $r \alpha^n$ are all in $(w,r) \subset (w,s)$.
The important conclusion is that under the hypothesis
$s \alpha, s \alpha^{-1}$ both not in $[s,w]$ then
any $s \alpha^n$ bridges to $[s,w]$ in a point
in the interior of $[s,w]$, Hence all $s \alpha^n$
are in $\we$ and $\vv$.

Use $s \tau^{-1} \alpha \tau  = s \gamma \beta \alpha^m$.
Here $s \alpha$ is in $\we$, so $s \alpha \tau$ is in
$\we \tau$.
Also $s \beta = s \alpha^{-1} \tau^{-1}$ is in $\we \tau^{-1}$
and bridges to $s$ in $[s,w]$. Hence
$s \beta \alpha^m$ bridges to $s \alpha^m$ in
$[s \alpha^m, w]$. 
But $s \alpha^m$ is in
$\we$ and bridges to $[s,w]$ in
a point in the interior of $(s,w)$.
This implies $s \beta \alpha^m$ is in $\we$,
contradiction.

This case is impossible.

\vskip .2in
\noindent
{\bf {Case C.1.3}} $-$ Suppose $s \alpha \in [s,w]$.

This implies for instance that $\we \alpha \subset \we$ and \
$T_s(w \tau^{-1}) \beta^{-1} \ \subset \ T_s(w \tau^{-1})$.

\vskip .2in
\noindent
{\bf {Case C.1.3.1}} $-$ Suppose $s \alpha^{-1} \in \we \tau$.

Then $s \beta^{-1} = s \alpha \tau^{-1}$
is in $(s, w \tau^{-1}) \subset \we \tau^{-1}$.
Also $s \alpha^{-1} = s \beta \alpha^{-1} \beta^{-1}$.
Here $s \beta = s \alpha^{-1} \tau^{-1}$ is in $\we$.

In this case suppose first that 
$s \beta$ is not in $\vv$. 
Then

$$w \ \in \
[w \tau^{-1}, s \beta] \ \ \ {\rm and} \ \ \
w \beta^{-1} \ \in \ [w \tau^{-1}, s] \ \ \ \Rightarrow \ \ \
w \beta^{-1} \alpha^{-1} \ \in \ \we \tau,$$

\noindent
as $s \alpha^{-1}$ is in $\we \tau$.
Notice $w \beta^{-1} \alpha^{-1}$ is not $s$.
Then

$$w \beta^{-1} \alpha^{-1} \gamma \ = \
w \alpha^{-1} \beta^{-1} \ = \ 
w \beta^{-1} \ \ \ {\rm is \ in} \ \  \we \tau^{-1}.$$

\noindent
Notice if $w \beta^{-1} = s$, then

$$ w \beta^{-1} \alpha^{-1} \  = \ 
w \beta^{-1} \gamma^{-1} \ = \ s \gamma^{-1}  \ =  \ s
\ = \ w \beta^{-1},$$

\noindent
contradiction
because $s$ is not fixed by $\alpha$.

Collecting all of this together: $w \beta^{-1} \alpha^{-1}
\gamma$ is in $\we \tau \gamma$. This point is equal
to $w \beta^{-1}$ which is in $\we \tau^{-1}$.
Therefore 

$$\we \tau \gamma  \ = \ \we \tau^{-1} \ \ \ {\rm or} \ \ \ 
\we \tau^2 \gamma = \we, \ \ \ \ \ {\rm impossible \ when} \
|p - 2q| = 1,$$

\noindent
as in case B.2.4.

\begin{figure}
\centeredepsfbox{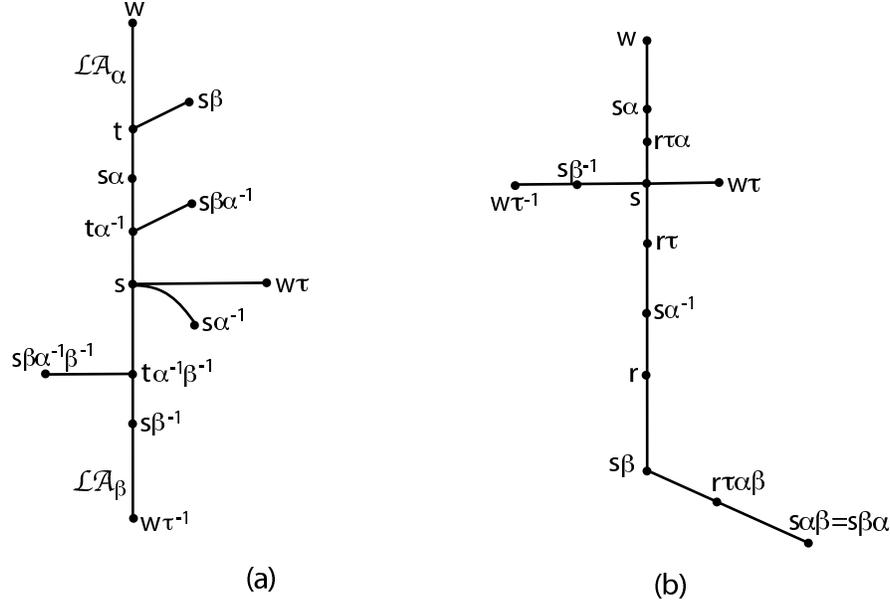}
\caption{
a. Case C.1.3.1, b. Case C.1.3.2.}
\label{16}
\end{figure}

The second option in case C.1.3.1 is that
$s \beta \in \vv$.
Recall that $s \alpha^{-1} \tau^{-1} = s \beta$ is in
$\we$.
Notice that 

$$\lab \ = \ (\lal) \tau^{-1} \ \ \ {\rm has \ a \ segment} \ \ 
[w \tau^{-1}, s] \ \subset \ 
\we \tau^{-1} \cup \{ s \}$$

\noindent
 and
then it goes into $\we$, as $s \beta$ is in $\we$.
Then either $s \beta = t \in (w,s)$ or
$s \beta$ bridges to $[w,s]$ in $t \in (w,s)$,
so bridges to $t$ in $\lal$.
In either case $s \beta \alpha^{-1}$ bridges 
to $t \alpha^{-1}$ in $\lal$ or is $t \alpha^{-1}$.
If $t \alpha^{-1}$ is in $[w,s)$, then
$s \beta \alpha^{-1}$ bridges to $t \alpha^{-1}$
in $\lab$, see fig. \ref{16}, a. Here $t \alpha^{-1}$ is
in $[w \tau^{-1}, s \beta)$.
If

$$s \ \in \ [t \alpha^{-1}, w] \ \ \ {\rm then} \ \ \
s \beta \alpha^{-1} \ \ \ {\rm bridges \ to} \ \ 
\lab \ \ {\rm in} \ \ r, \ \ {\rm with} \ \ 
r \ \in \ [s,w \tau^{-1}].$$

\noindent
This depends for instance on
whether $\we\tau = \we \tau^{-1}$ or not.
In any case 
$s \beta \alpha^{-1}$ bridges to $\lab$ in a point
in $[w \tau^{-1}, s \beta)$.
It follows that 
$s \beta \alpha^{-1} \beta^{-1}$ bridges
to a point $t$ in $\lab$ with $t$ in $[w \tau^{-1}, s)$,
that is, $s \beta \alpha^{-1} \beta^{-1}$ is in
$\we \tau^{-1}$.
Then

$$ s \alpha^{-1} \ \in \ \we \tau, \ \
s \beta \alpha^{-1} \beta^{-1} \ = \ 
s \alpha^{-1} \gamma \ \in \ \we \tau^{-1}
\ \ \  \Rightarrow \ \ 
\we \tau \gamma  \ = \  \we \tau^{-1},$$

\noindent
contradiction when $|p - 2q| = 1$.

This shows that case C.1.3.1 cannot occur.


\vskip .2in
\noindent
{\bf {Case C.1.3.2}} $-$ $s \alpha^{-1}$ is not in $\we \tau$.

Here $s \beta = s \alpha^{-1} \tau^{-1}$ is not in $\we$.
Also $s \beta^{-1} = s \alpha \tau^{-1}$ is not in $\we$.
It follows that 

$$\lab  \ \cap \ [w,s]  \ = \ \{ s \},$$

\noindent
so $s \alpha$ bridges to
$\lab$ in $s$ and 
$s \alpha \beta = s \beta \alpha$ bridges to $\lab$ in $s \beta$.
Hence $[s, s \beta] \subset (s \alpha, s \beta \alpha)$
and there is a fixed point $r$ of $\alpha$ 
in $(s, s \beta)$, see fig. \ref{16}, b.
It also implies that 

$$s \alpha^{-1} \in 
[s, r] \ \ \ {\rm and} \ \ \ 
s \alpha^{-1} \ \in  \ T_s (s \beta) \ = \ 
T_s(s \beta)   \tau,$$

\noindent
because $s \beta^{-1} = s \alpha \tau^{-1}$.
Now  apply
$\tau \alpha \beta = \alpha \tau \alpha^{1-m}$ to
$r$: \
$r \tau \alpha \beta = r \tau \alpha^{1-m}$.

As $s \beta \tau = s \alpha^{-1}$ and as
$r \in (s, s \beta)$, then

$$r \tau \ \in \ (s, s \alpha^{-1})
\ \  \Rightarrow \ \ 
r \tau \alpha \ \in \ (s, s \alpha) \ \ \Rightarrow \ \
r \tau \alpha \beta  \ \in  \ (s \beta, s \alpha \beta)
\ \subset \ T_r(s \beta).$$

\noindent
As $r \tau \alpha$ is in $(s, s \alpha^{-1}) \subset
T_r(s)$ this implies $r \tau \alpha^{1-m}$ is also
in $T_r(s)$.
Therefore 
$r$ separates $r \tau \alpha^{1-m}$ from
$r \tau \alpha \beta$,
contradiction.

This shows that case C.1.3, $s \alpha \in [s,w]$ cannot
occur.
Finally consider:

\vskip .2in
\noindent
{\bf {Case C.1.4}} $-$ Suppose $s \alpha^{-1} \in [s,w]$.

This implies that $\we \alpha^{-1} \subset \we$ and
\ $(\we \tau^{-1}) \beta \ \subset \ (\we \tau^{-1})$.

\vskip .2in
\noindent
{\bf {Case C.1.4.1}} $-$ Suppose $s \alpha \not \in \we \tau^{-1}$.

This case is very similar to case C.1.3.2.
Here $s \beta \in T_s(w \tau^{-1})$ which is not equal
to either $T_s(s \alpha)$ or $T_s(s \alpha^{-1})$.
Hence $s \beta$ bridges to $\lal$ in $s$ and $s \beta \alpha
= s \alpha \beta$ bridges to $\lal$ in $s \alpha$.
Hence 

$$s \beta \ \prec \ s \ \prec \ s \alpha \ \prec \ s \alpha \beta$$

\noindent
and there is a fixed point $r$ of $\beta$ in
$(s, s \alpha)$.
Then $s \beta^{-1} \in (s, r) \subset (s, s \alpha)$.
Now use $\beta \tau^{-1} \beta^{1-m} = \tau^{-1} \beta \alpha$
applied to $r$: \ $r \tau^{-1} \beta^{1-m} = r \tau^{-1} 
\beta \alpha$. As $s \alpha \tau^{-1} = s \beta^{-1}$ then

$$r \tau^{-1} \ \in  \ (s, s \beta^{-1})
\ \ \ \ \
{\rm so} \ \ \ \ \ r \tau^{-1} \beta^{1-m} \ \in \ 
(r, s \beta^{1-m}) \ \subset \ T_r(s).$$

\noindent
On the other hand $r \tau^{-1} \beta \alpha$ is in
$(s \alpha, s \beta \alpha) \subset T_r(s \alpha)$.
As $T_r(s \alpha) \not = T_r(s)$, this is a contradiction,
ruling out this case.


\vskip .2in
\noindent
{\bf {Case C.1.4.2}} $-$
$s \alpha$ is in $\we \tau^{-1}$.

This is similar to case C.1.3.1.
Suppose first that  $\we \tau^{-1} = \we \tau$.
Then $s \alpha \tau^{-1} = s \beta^{-1}$ is
in $\we$.
Also $\we \beta^{-1}$ is contained in $\we$.
It follows that

$$s \alpha^{-1} \beta^{-1} \ \in \ \we
\ \ \ {\rm and} \ \ \ 
s \alpha^{-1} \beta^{-1} \gamma^{-1} \ = \ s \beta^{-1} 
\alpha^{-1} \ \in \ \we.$$

\noindent
Hence $\we \gamma = \we$, $\we \tau^2 = \we$, leading
to contradiction when $p$ is odd.

Suppose now that 
$\we \tau^{-1} \not = \we \tau$.
Then $s \alpha \in \we \tau^{-1}$ and $s \alpha \tau^{-1}
= s \beta^{-1}$ is not in $\we$.
Also $s \beta^{-1}$ is in $\we \tau^{-2}$.
So $s \beta^{-1}$ bridges to $s$ in $\lal$ and
$s \beta^{-1} \alpha^{-1}$ bridges to
$s \alpha^{-1}$ in $\lal$
implying $s \beta^{-1} \alpha^{-1}$ is in $\we$.

Also $s \beta^{-1} \alpha^{-1} \gamma = s \alpha^{-1} \beta^{-1}$.
Here $s \alpha^{-1}$ bridges to $s$ in $\lab$,
$s \alpha^{-1} \beta^{-1}$ bridges to $s \beta^{-1}$ in 
$\lab$.
But 

$$s \beta^{-1} \ \in \ \we \tau^{-2} \ \ \ \Rightarrow \ \ \ 
s \alpha^{-1} \beta^{-1} \ \in \ \we \tau^{-2} \ \ \ \Rightarrow \ \ \ 
\we \gamma \ = \ \we \tau^{-2},$$

\noindent
again impossible
when $|p-2q| = 1$.

\vskip .1in
This finishes the analysis of case C.1.4,  \
$s \alpha^{-1} \in [s,w]$.

We conclude that case C.1, $\we \tau \not = \we$ is
impossible. This implies $\we \tau = \we$.
We stress that this does not yet produce
a local axis of $\tau$ in $\we$, because
we may have other fixed points of $\tau$ in
$(s,w)$. 

\vskip .2in
\noindent
{\bf {Case C.2}} $-$ Suppose that $\vv \alpha \not = \vv$.

Here we will use $s \alpha \tau = s \beta \alpha^m 
= s \alpha \beta \alpha^{m-1}$ many times.

\vskip .2in
\noindent
{\bf {Case C.2.1}} $-$ Suppose $w \tau, w \tau^{-1}$ are
not in $[s,w]$.

The bridge from $w \tau$ to $[s,w]$ is
$[w \tau, t]$, where $t$ is in $(s,w)$.
Since $s \alpha \not \in \vv$, then
$s \alpha \tau$ bridges to $t$ in $[s,w]$,
so $s \alpha \tau$ is in $\vv$.
Hence $s \alpha \tau \alpha^{-m}$ is 
in $\vv \alpha^{-m}$.
This point is equal to $s \beta = s \alpha^{-1} \tau^{-1}$.
In the same way $s \alpha^{-1}$ is not in $\vv$ and
bridges to $[s,w]$ in $w$. It follows that 
$s \alpha^{-1} \tau^{-1}$ bridges to a point $r$
in $[s,w]$, where $r$ is in
in $(s,w)$, hence $s \beta \in \vv$.
Therefore $\vv \alpha^m = \vv$.

On the other hand

$$s \alpha \tau \ = \ s \alpha \beta \alpha^{m-1} \ = \ 
s \alpha \tau \alpha^{-1} \tau^{-1} \alpha^{m-1}.$$

\noindent
The point $s \alpha \tau$ is in $\vv$ and
bridges to $t$ in $[s,w]$. So $s \alpha \tau \alpha^{-1}$
is in $\vv \alpha^{-1}$ and bridges to $w$ in 
$[s,w]$ so $s \alpha \tau \alpha^{-1} \tau^{-1}$ bridges
to $r$ in $[s,w]$ ($r$ as above) and as a result
this point is in $\vv$. Hence $s \alpha \beta \alpha^{m-1}$
is in $\vv \alpha^{m-1}$ and
$\vv \alpha^m = \vv \alpha^{m-1}$, contradicting
$\vv \alpha \not = \vv$.


\begin{figure}
\centeredepsfbox{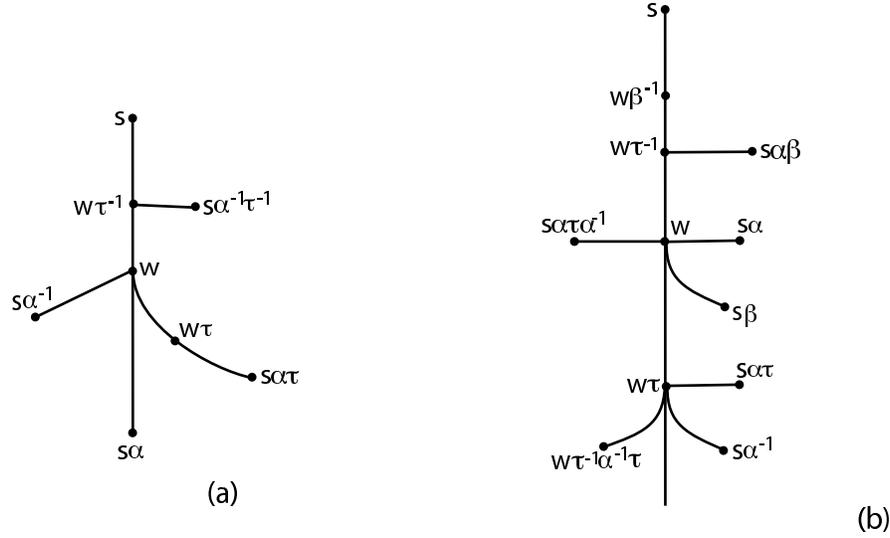}
\caption{
a. Case C.2.2.1, b. Case C.2.2.3.}
\label{17}
\end{figure}

\vskip .2in
\noindent
{\bf {Case C.2.2}} $-$ $w \tau^{-1} \in [s,w]$.

Here $\vv \tau^{-1}$ is contained in $\vv$.

The condition implies that 
$w$ is in a local axis $\lat$ of $\tau$
(this case will be ruled out, we only establish
the existence of a local axis of $\tau$ in $\we$ later).
Put an order $<$ in $\lat$ so 
$c < d$ in $\lat$ in $\lat$ if $s \prec c \prec d$ -
the order decreases
as points get closer to $s$.

\vskip .2in
\noindent
{\bf {Case C.2.2.1}} $-$
$w \tau \in \vv \alpha$, \ $w \tau \not \in \vv \alpha^{-1}$,
see fig. \ref{17}, a.

Here $\vv \alpha \tau \subset \vv \alpha$.

The conditions imply in particular that
$\vv \alpha \not = \vv \alpha^{-1}$.
Here $s \alpha \tau \in \vv \alpha$,
so $s \beta \alpha^m \in \vv \alpha$. 
Also $s \alpha^{-1}$ bridges to
$\lat$ in $w$
so $s \beta = s \alpha^{-1} \tau^{-1}$ bridges
to $\lat$ in  $w \tau^{-1}$.
It follows that $s \beta$ is in $\vv$ and
$s \beta \alpha^m$ is in $\vv \alpha^m$.
Hence $\vv \alpha^m = \vv \alpha$.

On the other hand
$s \alpha \tau = s \alpha \beta \alpha^{m-1}$.
Use 
$s \alpha \beta \ = \ s \alpha \tau \alpha^{-1} \tau^{-1}$.
Here

$$s \alpha \in \vv \alpha \ \ \Rightarrow \ \ 
s \alpha \tau \in \vv \alpha \ \ \Rightarrow \ \ 
s \alpha \tau \alpha^{-1} \in \vv \ \ \Rightarrow \ \ 
s \alpha \tau \alpha^{-1} \tau^{-1} \in \vv.$$

\noindent
Finally $s \alpha \beta \alpha^{m-1}$ is in
$\vv \alpha^{m-1}$. So $\vv \alpha^{m-1} = \vv \alpha$
and $\vv = \vv \alpha$, again contradicting the
assumption
in this case.

\vskip .2in
\noindent
{\bf {Case C.2.2.2}} $-$
Suppose $w \tau$ is not in $\vv \alpha$ and
$w \tau$ is not in $\vv \alpha^{-1}$.

Then $w \tau$ is in $\re$ another  component
of $T - \{ w \}$.
Then $s \alpha \tau$ is in $\re$.
Now $s \beta \alpha^m = s \alpha^{-1} \tau^{-1} \alpha^m$.
But

$$w \tau \ \not \in \ \vv \alpha^{-1} \ \ \ \Rightarrow \ \ \ 
s \alpha^{-1} \ \  \ {\rm bridges \ to} \ \ 
\lat \ \ {\rm in} \ w \ \ \ \Rightarrow \ \ \
s \alpha^{-1} \tau^{-1} \ \ {\rm bridges \ to}  \ \ \lat 
\ \ {\rm in } w \tau^{-1}$$

\noindent
and
$s \beta$ is in $\vv$. Therefore
$s \beta \alpha^m \in \vv \alpha^m = \re$.
Notice $\re \alpha^{-1}
\not = \re$ because $\re = \vv \alpha^m$ and
$\vv \alpha^{-1} \not = \vv$.
Use 

$$s \alpha \tau \ = \ s \alpha \beta \alpha^{m-1} \ = \
 s \alpha \tau \alpha^{-1} \tau^{-1} \alpha^{m-1}
\ \ \ \ {\rm and} \ \ \ \ 
s \alpha \tau \alpha^{-1}  \ \in 
\re \alpha^{-1}
\not = \re.$$

\noindent
Hence $s \alpha \tau \alpha^{-1}$ bridges
to $\lat$ in a point $\leq w$ in $\lat$ 
(it is in $[s,w]$) and $s \alpha \beta$
bridges to $\lat$ in
a point $\leq w \tau^{-1}$ in $\lat$.
Hence 

$$s \alpha \beta \in
\vv \ \ \ \Rightarrow \ \ \ s \alpha \beta \alpha^{m-1}
\in \vv \alpha^{m-1} \ \ \  \ \Rightarrow \ \  \ \ 
\vv \alpha^m = \vv \alpha^{m-1},$$ 

\noindent
contradiction.
Notice that here it doesn't matter whether
$\vv \alpha = \vv \alpha^{-1}$ or not.

\vskip .2in
\noindent
{\bf {Case C.2.2.3}} $-$ $w \tau$ is in $\vv \alpha^{-1}$,
see fig. \ref{17}, b.

This implies $\vv \alpha^{-1} \tau$ is a subset of $\vv \alpha^{-1}$.

Use $s \alpha \tau = s \beta \alpha^m = s \alpha^{-1}
\tau^{-1} \alpha^m$. 
Here

$$s \alpha \not \in \vv \ \ \ \Rightarrow \ \ \ 
s \alpha \tau \in T_w(w \tau) = \vv \alpha^{-1}
\ \ \ \Rightarrow \ \ \ 
s \alpha \tau \alpha^{-1} \ \in \ \vv \alpha^{-2} 
\not = \vv \alpha^{-1},$$

\noindent
so it bridges to a point $r$ in $\lat$ with $r \leq w$ in $\lat$.
Hence $s \alpha \beta$ is in $\vv$ and
$s \alpha \beta \alpha^{m-1}$ is in $\vv \alpha^{m-1}$.
Hence $\vv \alpha^{m-1} = \vv \alpha^{-1}$ or
$\vv \alpha^m = \vv$.


On the other hand 
$s \alpha \tau = s \beta \alpha^m$
is in $\vv \alpha^{-1}$,
so 

$$s \alpha^{-1} \tau^{-1} \ = \
 s \beta  
\ \ \ \ {\rm is \ in} \ \ \ \ 
\vv \alpha^{-1-m}  \ =  \
\vv \alpha^{-1}.$$


\noindent
Then $s \beta$
bridges to a point 
$> w$ in $\lat$. But $s \beta = s \alpha^{-1} \tau^{-1}$,
so $s \alpha^{-1}$ bridges to a point $> w \tau$ in
$\lat$, which implies $w \tau \in (w, s \alpha^{-1})$.
It follows that $w \tau \alpha \in (w, s)$ 
and
$w \beta^{-1} = w 
\tau \alpha \tau^{-1}$ is in $(w \tau^{-1}, s)$
is in $\we$ and in $\vv$.

The following arguments use the strategy of case R.2:

Now $w \beta^{-1} \alpha^{-1} \gamma = w \beta^{-1}$ is in
$\we$ and 
$w \beta^{-1} \alpha^{-1} = w \tau^{-1} \alpha^{-1} \tau$.
Use 

$$w \tau^{-1} \alpha^{-1} \ \in \ \vv \alpha^{-1} \
 = \ T_w(w \tau) \ \ \ {\rm so} \ \ \ 
w  \tau^{-1} \alpha^{-1} \tau \ \ \ {\rm is \ in} \ 
T_{w \tau}(w \tau^2) = \vv \alpha^{-1} \tau^{-1}
\subset \vv \alpha^{-1} \ \subset \ \we.$$

\noindent
Hence
$w \beta^{-1} \alpha^{-1} \in \we$.
From this it follows that $\we \gamma = \we$.
As usual this implies that $(\lat) \gamma = \lat$
so $\gamma, \tau$ have the common local axis $\lat$.
In addition $w \beta^{-1} \alpha^{-1} \gamma = w \beta^{-1}$
and as $w \beta^{-1}$ is in $\lat$, so does
$w \beta^{-1} \alpha^{-1}$.

If $w \tau \alpha \leq w \tau^{-1}$ in $\lat$ then
$w \beta^{-1} = w \tau \alpha \tau^{-1} \leq w \tau^{-2}$
in $\lat$. 
Also $w \tau, w \beta^{-1} \alpha^{-1}$ are in $\lat$
and $w \tau < w \beta^{-1} \alpha^{-1}$ in $\lat$.
Hence 

$$w \tau \gamma \ < \ w \beta^{-1} \alpha^{-1} \gamma \ = \ 
w \beta^{-1} \leq w \tau^{-2} 
\ \ \ {\rm in } \ \ \lat \ \ \ \ \Rightarrow \ \
p > 3q,$$

\noindent
contradiction to $|p-2q| = 1$.

If $w \tau \alpha > w \tau^{-1}$ in $\lat$ then 
$w \tau \alpha \tau^{-1} = w \beta^{-1} \in (w \tau^{-2}, 
w \tau^{-1})$. 
Here use 

$$(w \tau^2) \gamma \beta \alpha^m \ = \  w \tau \alpha \tau
\ 
\in T_w(w \tau) = \vv \alpha^{-1} \ \ \ \Rightarrow \ \ 
w \tau^2 \gamma \beta \ \in \vv \alpha^{-1},$$

\noindent
because 
$\vv \alpha^m = \vv$.
Therefore $w \tau^2 \gamma \beta$ bridges to $v$ in
$\lat$ with $v > w$ in $\lat$.
Hence $w \tau^2 \gamma < w \beta^{-1}$ in $\lat$ and
as $w \beta^{-1} < w \tau^{-1}$ we also obtain
$p > 3q$, contradiction.

This rules out the case C.2.2.3 and hence finishes
the analysis of case C.2.2, \ $w \tau^{-1} \in [s,w]$.
The next case is:

\vskip .2in
\noindent
{\bf {Case C.2.3}} $-$ $w \tau \in [s,w]$.

This implies that $\vv \tau \subset \vv$. The case is
similar to case C.2.2.

\vskip .2in
\noindent
{\bf {Case C.2.3.1}} $-$ 
$w \tau^{-1} \in \vv \alpha^{-1}, 
\ w \tau^{-1} \not \in \vv \alpha$.

This implies that 
$\vv \alpha^{-1} \tau^{-1} \subset \vv \alpha^{-1}$.

Here $w \tau^{-1} \alpha$ is in $\vv$,  \
$w \tau^{-1} \alpha \tau$ is in $\vv$ so
$w \alpha \beta \alpha^{m-1} = w \beta \alpha^{m-1}$
is in $\vv$.
Also 

$$w \tau \alpha^{-1} \ \in  \ \vv \alpha^{-1} \ \Rightarrow \
w \beta = w \tau \alpha^{-1} \tau^{-1} \in \vv \alpha^{-1}
\ \Rightarrow \
w \beta \alpha^{m-1} \in \vv \alpha^{m-2}$$ 

\noindent
which
must be equal to $\vv$.

On the other hand $s \alpha \tau = s \beta \alpha^m$.
Here $s \alpha \in \vv \alpha$ and bridges to $w$ in $\lat$,
so $s \alpha \tau$ bridges to $w \tau$ in $\lat$
and $s \alpha \tau \in \vv$.
Also

$$s \beta \ = \ s \alpha^{-1} \tau^{-1}  \ \in \ \vv \alpha^{-1}
\ \ \ {\rm and} \ \ \ 
s \beta \alpha^m \ \in \ \vv \alpha^{m-1}.$$

\noindent
It follows that $\vv \alpha^{m-1} = \vv \alpha^{m-2}$,
contradiction to $ \vv \not = \vv \alpha$.

\begin{figure}
\centeredepsfbox{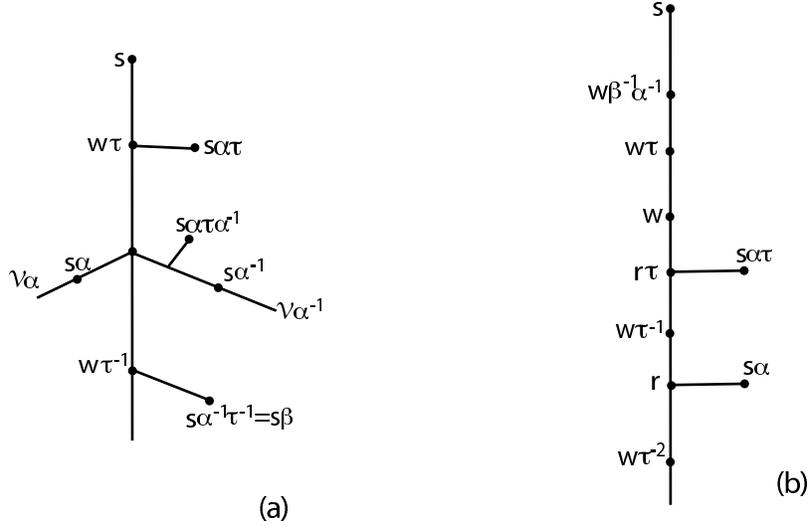}
\caption{
a. Case C.2.3, b. Case C.2.3.3, Situation III.}
\label{18}
\end{figure}

\vskip .2in
\noindent
{\bf {Case C.2.3.2}} $-$
$w \tau^{-1} \not \in \vv \alpha^{-1},
\ w \tau^{-1} \not \in \vv \alpha$,
see fig. \ref{18}, a.

Use $s \alpha \tau = s \beta \alpha^m 
= s \alpha \beta \alpha^{m-1}$.
In this case the point 
$s \alpha$ brides to $w$ in $\lat$
and $s \alpha \tau \in \vv$.
Also 
$s \alpha^{-1}$ bridges to $w$ in $\lat$
and $s \beta = s \alpha^{-1} \tau^{-1}$ bridges 
$w \tau^{-1}$ in $\lat$ so

$$s \beta \ \ \ {\rm is \ in} \ \
\re = T_w(w \tau^{-1}) \ \not = \ 
\vv \alpha, \vv \alpha^{-1} \ \ \Rightarrow \ \
s \beta \alpha^m \in \re \alpha^m = \vv.$$

\noindent
So in particular $\re \not = \re \alpha$.

On the other hand $s \alpha \tau \alpha^{-1} \in \vv \alpha^{-1}$
and bridges to $w$ in $\lat$ so $s \alpha \beta = 
s \alpha \tau \alpha^{-1} \tau^{-1}$ bridges
to $w \tau^{-1}$ in $\lat$ and is in $\re$.
Then $s \alpha \beta \alpha^{m-1} \in \re \alpha^{m-1}
= \vv \alpha^{-1}$. This would imply $\vv = \vv \alpha^{-1}$,
contradiction. 

The final case in C.2.3 is:

\vskip .2in
\noindent
{\bf {Case C.2.3.3}} $-$ $w \tau^{-1} \in \vv \alpha$.

Let $[s \alpha, r]$
be the bridge from
$s \alpha$ to $\lat$
with $r$ in $\lat$. Then $r > w$ in $\lat$.
Here we have to subdivide.

\vskip .2in
\noindent
{\bf {Situation I}} $-$ $r$ is in $(w, w \tau^{-1})$.

Then $s \alpha \tau$ bridges to $\lat$ in
$r \tau \in (w, w \tau)$ and $s \alpha \tau \in \vv$.
Hence 

$$s \alpha \tau \alpha^{-1} \not \in \vv \ \ \Rightarrow \ \ 
s \alpha \tau \alpha^{-1} \tau^{-1} \ = \ s \alpha \beta \
\in 
\vv \alpha \ \ \Rightarrow \ \ 
s \alpha \beta \alpha^{m-1} \ \in \ \vv \alpha^m
\ \ \Rightarrow \ \ 
\vv = \vv \alpha^m.$$

On the other hand $s \beta \alpha^m = s \alpha^{-1} \tau^{-1} 
\alpha^m$. Here $s \alpha^{-1} \tau^{-1}$ is in $\vv \alpha$
so $s \beta \alpha^m$ is in $\vv \alpha^{m+1}$, 
implying $\vv \alpha^m = \vv \alpha^{m+1}$ again a contradiction.

\vskip .2in
\noindent
{\bf {Situation II}} $-$ $r = w \tau^{-1}$.

Here $s \alpha \tau$ bridges to $\lat$ in $w$ hence
$s \alpha \tau \not \in \vv \alpha$ and $s \alpha \tau
\not \in \vv$.
So $s \alpha \tau$ is in $\re$, another component
of $T - \{ w \}$.
Also 

$$s \alpha^{-1} \ \not \in \vv \ \ \Rightarrow \ \ 
s \beta  \ = \ s \alpha^{-1} \tau^{-1} \ \in \  \vv \alpha
\ \ \Rightarrow \ \ 
s \beta \alpha^m \ \in \ \vv \alpha^{m+1} \ \ \Rightarrow \ \ 
\re = \vv \alpha^{m+1}.$$

On the other hand $s \alpha \beta \alpha^{m-1} \in \vv \alpha^{m+1}$,
so $s \alpha \beta \in \vv \alpha^2$.
Now $\vv \alpha^2 \not = \vv \alpha$ so $\vv \alpha^2 \tau$
is contained in $\vv$. Hence 
$s \alpha \tau \alpha^{-1} = s \alpha \beta \tau$ is in
$\vv$.
This would imply $s \alpha \tau$ is in $\vv \alpha$,
contradiction to the first conclusion in
this case.


\vskip .2in
\noindent
{\bf {Situation III}} $-$ $w  \tau^{-1} < r$ in $\lat$. 

This is a little more tricky.
Here $s \alpha \tau \in \vv \alpha$, see fig. \ref{18}, b.
Also 

$$w \beta^{-1}  \ = \  w \tau \alpha \tau^{-1} \ \in \ \vv \alpha
\ \subset \we.$$

\noindent
Now use 
$w \beta^{-1} \alpha^{-1}
= w \tau^{-1} \alpha^{-1} \tau$.
Here

$$w \tau^{-1} \in [w, s \alpha] \ \ \Rightarrow \ \ 
w \tau^{-1} \alpha^{-1} \in [w,s] \ \ \Rightarrow \ \ 
w \tau^{-1} \alpha^{-1} \tau \in
[s, w \tau] \subset \lat \subset \we.$$

\noindent
So $w \beta^{-1} \alpha^{-1}, w \beta^{-1}$ are both
in $\we$, with the usual implications
that $\we \gamma = \we$ and $\gamma$
leaves $\lat$ invariant.
As $w \beta^{-1} \alpha^{-1}$ is in $\lat$ then 
$w \beta^{-1}$ is in $\lat$ as well.
Also


$$w \beta^{-1} \alpha^{-1}  \ = \ 
w \tau^{-1} \alpha^{-1} \tau < w \tau \ \in \
\lat \ \ \ \Rightarrow \ \ \ 
w \tau \alpha \ \in \ \lat.$$

The proof is now analogous to previous arguments.
If

$$w \ \prec \ w \tau^{-1} \ \prec \ w \tau \alpha
\ \ \ \Rightarrow \ \ \ 
w \tau^{-1} \  \prec \ w \tau^{-2} \ \prec \ w \tau \alpha \tau^{-1}
\ = \ w \beta^{-1}.$$

\noindent
But 

$$w \beta^{-1} \alpha^{-1} \gamma \ = \ w \beta^{-1}
\ \ \ {\rm and} \ \ \ 
w \beta^{-1} \alpha^{-1} \ \in \ (s, w \tau)$$

\noindent 
implies
as before that
$p > 3q$, contradiction.

On the other hand if $w \prec w \tau \alpha \prec w \tau^{-1}$,
then $w \tau^{-1} \ \prec \ w \tau \alpha \tau^{-1} = w \beta^{-1}
\ \prec \ w \tau^{-2}$
all in $\lat$.
Here $s \alpha \tau \in \vv \alpha$. Now
$s \beta \alpha^m \ = \ s \alpha^{-1} \tau^{-1} \alpha^m$.
Also

$$s \alpha^{-1}  \ \not \in \ \vv \ \ \Rightarrow \ \ 
s \alpha^{-1} \tau^{-1} \ \in \ \vv \alpha \ \ \Rightarrow \ \ 
s \beta \alpha^m \ \in \ \vv \alpha^{m+1} \ \ \ \Rightarrow \ \ \ 
\vv \alpha = \vv \alpha^{m+1} \ \ {\rm or}  \
\vv \ = \ \vv \alpha^m.$$

Now use $w \tau^2 \gamma \beta \alpha^m = w \tau \alpha \tau$.
Here 

$$w \tau \ \prec \ w \tau \alpha \tau 
\ \prec  \ w \ \ \ {\rm in} \ \ \lat \ \ \ \Rightarrow \ \ \ 
w \tau \alpha \tau \ \in \ \vv, \ \ 
w \tau^2 \gamma \beta \ \in \ \vv \alpha^{-m} \ = \ \vv.$$

\noindent
So $w \tau^2 \gamma \prec w \beta^{-1} \prec w \tau^{-1}
\prec w$,
implying again $p > 3q$, contradiction.

This finishes the analysis of case C.2.3, $w \tau \in [s,w]$
and so proves that the case $\vv \alpha \not = \vv$ cannot
occur.
From now on in case C assume:

\vskip .2in
\noindent
{\bf {Case C.3}} $-$ 
$\we \tau = \we$ and $\vv \alpha = \vv$.

Since there is no other fixed point
of $\alpha$ in $(s,w)$,
this immediately implies there is a local axis $\lal$ of
$\alpha$ contained in $\vv$ 
with $w$ as an ideal point of $\lal$.
We stress that at this point we do not yet have an
axis for $\tau$, because there may be other fixed points
of $\tau$ in $(s,w)$.

\begin{lemma}{}{}
$s \alpha, s \alpha^{-1} \in \we$,
so $s \alpha, s \alpha^{-1}$ are not in $[s,w)$.
\label{all}
\end{lemma}

\begin{proof}{}
Suppose first that $s \alpha$ is not in $\we$.
Then

$$s \alpha^{-1} \ \in \ (s,w) \ \subset \ \we \ \ \ \Rightarrow \ \ \
s \alpha^{-1} \tau^{-1} \ \in \ \we \tau \ = \ \we.$$

\noindent
So $s \beta \in \we$ and bridges to $[s,w]$ in
a point $r$ which is in $(s,w]$.
Then $s \beta \alpha^m$ bridges to $[s,w]$ in
$r \alpha^m$ and $s \beta \alpha^m$ is in $\we$.
Therefore $s \alpha \tau$ is in $\we$ and
$s \alpha$ is in $\we \tau^{-1} = \we$, contradiction.

On the other hand suppose that $s \alpha^{-1} \not
\in \we$.
Then $s \alpha \in (s,w]$. Also $s \beta =
s \alpha^{-1} \tau^{-1} \not \in \we$, so bridges to $[s,w]$ in
$s$. Then $s \beta \alpha^m$ bridges to $[s \alpha^m, x]$
in $s \alpha^m$. Since $s \alpha^m \not \in \we$ this
implies $s \beta \alpha^m \not \in \we$, 
therefore $s \alpha \tau \not \in \we$. But then
$s \alpha$ is not in $\we$, contradiction.
This finishes the proof.
\end{proof}

We conclude that $s \alpha, s \alpha^{-1}$ are
in $\we \cap \vv$. Let $s \alpha$ bridge to $r$
in $[s,w]$, hence $r \in (s,w)$ and $s \alpha^{-1}$
bridges to $[s,w]$ in a point $t$ also in  $(s,w)$.

Let $z$ be the fixed point of $\tau$ in $[s,w]$ which is
closest to $w$. \ Then $z$ may be equal to $s$, but is not $w$.
Let $\uu = T_z(w)$. 
One important goal is to prove that $\uu \tau = \uu$.

\begin{lemma}{}{}
Let $\uu = T_z(w)$. Then $\uu \tau = \uu$.
If $z \not = s$  then
$z \gamma, w \gamma \not \in \we$, and
$z \alpha, \ z \alpha^{-1} \not \in (z,w)$.
\label{prel}
\end{lemma}

\begin{proof}{}
If $z = s$ then $\uu = \we$ and the result follows
from Case C.1.
For the rest of the proof of the lemma assume
that $s \not = z$.

We first analyse the possibility that $z \gamma \in \we$.
As $\kappa$ fixes $s$ then $z \gamma^{-1} \in \we$ also.
If $z \gamma = z$, then $z \kappa = z$, contradiction.

Suppose that $z \gamma$ or $z \gamma^{-1}$ is in
$[s,z)$.
Then as $s \gamma = s$, it follows that $z$ is in
a local axis for $\gamma$ and $z \gamma^q \not = z$,
contradiction to $z$ fixed by $\tau$.
Hence $z \gamma, z \gamma^{-1} \not \in [s,z]$.

Let $[z \gamma, r]$ be the bridge from
$z \gamma$ to $[s,z]$. Notice that
$r$ is in $(s,z)$, because $z \gamma, z \gamma^{-1}$
are not in $[s,w]$.
Then

$$r \ \in \ [s, z] \cap [s, z \gamma] \ \ \Rightarrow \ \ 
r \gamma^{-1} \in \ [s, z].$$

\noindent
If $r \gamma = r$, then $r \tau^p = r \gamma^{-q} = r$.
But $([s,z]) \tau = [s,z]$,
so this would imply $r \tau = r$.
Together these imply $r \kappa = r$, contradiction
to $s$ the fixed point of $\kappa$ in $[s,w]$
which is closest to $w$.

We conclude that $r \gamma \not = r$.
But as $s \gamma = s$, this implies that
$r$ is in a local axis $\lag$ of $\gamma$.
Compute $r^{nq}, n \in {\bf Z}$.
Assume without loss of generality 
that $r^{nq}$ moves away from $s$ as
$n \rightarrow +\infty$.
Then 

$$r \gamma^{nq} \ =  \ r \tau^{-np} \in 
[s,w], \ \forall n \ \ \ {\rm and} \ \ \
r \gamma^{nq} \rightarrow \ c  \in 
(s,z] \ {\rm as} \ n \rightarrow +\infty.$$

\noindent
Then $c \gamma = c$ and also $c \tau = c$, contradiction.

This contradiction shows that $z \gamma \in \we$
is impossible.
Notice that if $z \gamma$ is not in $\we$, then
$z \gamma$ separates $\we \gamma$ from $s$ and
hence from $\we$. It follows that
$\we \gamma \cap \we
= \emptyset$, so $w \gamma \not \in \we$.
This proves one assertion of lemma \ref{prel}.

We now consider where $z \alpha$ and $z \alpha^{-1}$ are.
Notice they are both in $\vv$.
Remember that for the rest of the proof $s \not= z$.

\vskip .1in
\noindent
{\bf {Situation I}} $-$
Suppose first that $z \alpha \in (z,w)$.

Use $\alpha \tau = \tau \gamma \beta \alpha^m$, applied
to $z$. Here $z \alpha$ is in
$\uu$ so $z \alpha \tau$ is in $\uu \tau$.
Suppose first that $\uu \tau \not = \uu \alpha^{-1}$.
Then $z \alpha \tau$ bridges to $\lal$ in a point
in $[z,w]$ and
hence $a = z \alpha \tau \alpha^{-m}$ bridges to
$\lal$ in a point in
$[z \alpha^{-m},w]$ and $a$ is in $\uu$.
Here

$$z \alpha \tau \alpha^{-m} \ = \ 
z \gamma \beta \ = \ z \gamma \alpha^{-1} \tau^{-1}
\ \ \ \Rightarrow \ \ \ 
z \gamma \alpha^{-1} \ \in \ \uu \tau \not = \uu.$$

\noindent
Again $z \gamma \alpha^{-1}$
bridges to $\lal$ in a point in $[z,w]$ and it follows that
$z \gamma$ is in $\uu$, hence 
$z \gamma \in \we$ contradicting $\we \gamma \cap \we = \emptyset$.

The remaining possibility is $\uu \tau = T_z( z \alpha^{-1})$,
so in particular $\uu \tau \not = \uu$, see fig. \ref{19}, a.
Consider $w \tau^{-1} \alpha^{-1} \tau$. The point
$w \tau^{-1}$ is not in $\uu$, hence it bridges
to $\lal$ in a point not in $(z,w]$. Therefore
$w \tau \alpha^{-1}$ bridges to $\lal$ in
a point not in $(z \alpha^{-1}, w]$, so $w \tau
\alpha^{-1}$ is in $T_z(z \alpha^{-1}) = \uu \tau$.
Hence 

$$w \tau^{-1} \alpha^{-1} \tau  \ = \ w \beta^{-1} \alpha^{-1}
\ \ \
{\rm is \ in} \ \ 
\uu \tau^2
\not = \ \uu \tau, \ T_z(s).$$

\noindent
Notice that

$$(T_z(s)) \tau \ = \ T_z(s), \ \ {\rm since } \ \ 
s \tau  \ = \ s, \ \ \ {\rm so} \ \ T_z(s) \ \not = \uu \tau^2.$$

\noindent
 In particular
$w \beta^{-1} \alpha^{-1}$ is in $\we$ and also bridges
to $\lal$ in a point which is in $[z,w]$. Then
$w \beta^{-1}$ bridges to $\lal$ in a point
which is in $[z \alpha, w]$ so in particular
$w \beta^{-1}$ is in $\uu \subset \we$.
But then $w \beta^{-1} \alpha^{-1}$ and $w \alpha^{-1}$ 
are both in $\uu$, contradicting 
$\we \gamma \cap \we = \emptyset$.


This finishes the analysis of possibility $z \alpha \in (z,w)$.


\begin{figure}
\centeredepsfbox{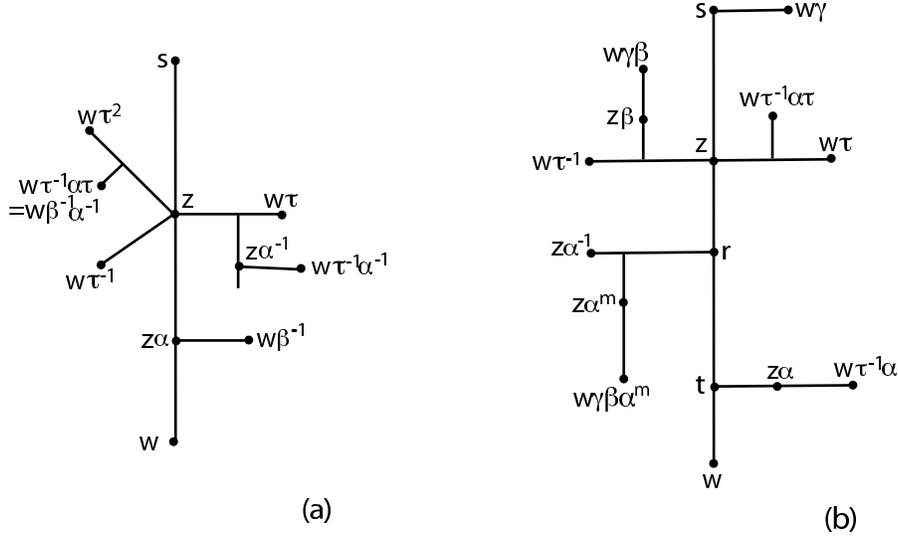}
\caption{
a. Situation I, Situation III.}
\label{19}
\end{figure}


\vskip .1in
\noindent
{\bf {Situation II}} $-$  Suppose $z \alpha^{-1} \in (z,w)$.

Consider first the case when $z \alpha \in \uu \tau^{-1}$,
that is $T_z( z \alpha) = T_z (w \tau^{-1})$. This is very
similar to Situation I, second part.
Since $z \alpha$ is not in $\uu$, this in particular
implies $\uu \tau \not = \uu$.
Here $w \tau \not \in \uu$, hence it bridges to $\lal$ in
a point which is not in $(z,w]$. It follows that
$w \tau \alpha$ bridges to $\lal$ in a point
which is not in $(z \alpha, w]$.
This implies  that $w \tau \alpha$ is in
$T_z(z \alpha) = T_z(w \tau^{-1})$.
Hence

$$w \beta^{-1} \ = \ w \tau \alpha \tau^{-1} 
\ \ \ {\rm is \ in} \ \
T_z(w \tau^{-2}) \not = \ 
T_z(s), \ T_z(w \tau^{-1}).$$

\noindent
The first
fact means that $w \beta^{-1}$ is in $\we$.
The second fact means that $w \beta^{-1}$ is not
in $T_z(z \alpha)$,  hence $w \beta^{-1}$
bridges to $\lal$ in
a point contained in $[z,w]$. Hence $w \beta^{-1}
\alpha^{-1}$ bridges to $\lal$ in a point contained
in $[z \alpha^{-1}, w]$ and is in $\we$. 
As $w \beta^{-1} \alpha^{-1} \gamma = w \beta^{-1}$,
this would imply $\we \gamma = \we$, again 
contradiction.
Hence this cannot occur.

Now we know $z \alpha$ is not in $T_z(w \tau^{-1})$.
The point $z \beta =
z \alpha^{-1} \tau^{-1}$ is in $T_z(w \tau^{-1})$,
hence it bridges to $\lal$ in a point
contained in $[z,w]$. It follows that $z \beta \alpha^m$
bridges to $\lal$ in a point contained in $[z \alpha^m, w]$.
But 

$$z \alpha^m  \in \ \uu \ \ \Rightarrow \ \  
z \beta \alpha^m \ \in \ \uu \ \ \Rightarrow \ \ 
z \gamma^{-1} \alpha \tau \ \in \ \uu \ \ {\rm or} \ \ 
z \gamma^{-1} \alpha \ \in  \ T_z(w \tau^{-1})$$

\noindent
and bridges to $\lal$ in a point in $[z,w]$.
It follows that $z \gamma^{-1}$ bridges to $\lal$
in a point in $[z \alpha^{-1}, w]$, hence
$z \gamma^{-1} \in \uu \subset \we$, impossible.

We conclude that situation II cannot occur.
This proves the last 2 assertions of the lemma \ref{prel}.
It also implies that the following situation must occur:


\vskip .2in
\noindent
{\bf {Situation III}} $-$ $z \alpha \not \in (z,w),
z \alpha^{-1} \not \in (z,w)$, see fig. \ref{19}, b.

What is left to prove of lemma \ref{prel} is that
$\uu \tau = \uu$. So suppose that $\uu \tau \not = \uu$.

Here $z \alpha^{-1}$ bridges to $[z,w]$ in a point
$r$ which is in $(z,w)$. Also $z \alpha$ bridges to $t$
in $[z,w]$ with $t$ also in $(z,w)$.

The point $w \gamma$ is not in $\we$, so it is in $T_z(s)$
and bridges to $[w \tau^{-1}, z]$ in $z$.
Hence $w \gamma \beta$ bridges to $[w \tau^{-1}, z \beta]$
in $z \beta$. But $z \beta = z \alpha^{-1} \tau^{-1}$
bridges to $[z, w \tau^{-1}]$ in $r \tau^{-1}$.
Then $w \gamma \beta$ bridges to $[z,w]$ in $z$
(this uses $\uu \tau \not = \uu$!).
Then 

$$w \gamma \beta \alpha^m \ \ \ {\rm bridges \ to} \ \ 
[z,w] \ \ {\rm in  \ a \ point \ in } \ \ 
(z,w) \ \ \ {\rm so} \ \ \ 
w \gamma \beta \alpha^m \ \in \uu.$$

On the other hand $w \tau^{-1}$ bridges to $[z,w]$ 
in $z$ so  $ w \tau^{-1} \alpha$ bridges to $[z,w]$
in a point in $(z,w)$
and $w \tau^{-1} \alpha$ is in $\uu$.
Then $w \tau^{-1} \alpha \tau$ is in $\uu \tau$.
Of course this implies $\uu \tau = \uu$, contrary to assumption.

So in any case we conclude that $\uu \tau = \uu$.
This finishes the proof of lemma \ref{prel}.
\end{proof}

This lemma is very useful.
Since there is no fixed point of $\tau$ in $(z,w)$ and
$T_z(w) \tau = T_z(w)$ it follows that there is
a local axis $\lat$ of $\tau$ contained in $\uu = T_z(w)$
with an ideal point $z$.

\begin{lemma}{}{}
$w$ is not in $\lat$.
\end{lemma}

\begin{proof}{}
Suppose not, that is, $w \in \lat$.
Here we will use lemma \ref{chax}:
Suppose that $\lat$ is a local axis for $\tau$ and
$w$ is a point in
$\lat$ with $w \alpha = w$. Then
at least one of the components of $T - \{ w \}$
containing
$w \tau, w \tau^{-1}$ is not invariant under $\alpha$.


\vskip .2in
\noindent
{\bf {Situation I}} $-$
 $w \tau^{-1} \in [z,w)$.

Here $\vv = T_w(z) = T_w(w \tau^{-1})$ is invariant under
$\alpha$. By lemma \ref{chax}, the set
${\cal R} = T_w(w \tau)$ is not invariant
under $\alpha$. Notice that $\ro \alpha$ is not equal
to $\vv$ either.

Use $w \alpha \tau = w \tau = w \tau \alpha \beta \alpha^{m-1}$.
Here 

$$w \tau \ \in \ \ro \ \ \Rightarrow \ \ 
w \tau \alpha  \ \in \ \ro \alpha \ \not = \ \vv 
\ \ \Rightarrow \ \ 
w \tau \alpha \tau \ \in \ \ro \tau \ \subset \ \ro
\ \ \Rightarrow \ \ 
c \ = \ w \tau \alpha \tau \alpha^{-1} \ \in \ \ro \alpha^{-1}
\ \not = \ \ro.$$

\noindent
So $c$ bridges to $w$ in $\lat$ and
then $w \tau \alpha \tau \alpha^{-1} \tau^{-1} = w \tau \alpha
\beta$ bridges to $w \tau^{-1}$ in $\lat$ and is then
in $\vv$. Finally $w \tau \alpha \beta \alpha^{m-1}$
is in $\vv \alpha^{m-1} = \vv$. This is not
$\ro$, contradiction.

\vskip .2in
\noindent
{\bf {Situation II}} $-$ $w \tau \in (z,w)$.

Here $\vv = T_w(w \tau) = T_w(z)$ is invariant under $\alpha$.
Let $\ro = T_w (w \tau^{-1})$, which is not
invariant under $\alpha$.
Use $w \tau^{-1} \alpha \tau = w \alpha \beta \alpha^{m-1}$.
Then $w \tau^{-1}$ is in $\ro$, so $w \tau^{-1} \alpha$
is not in $\ro$ or $\vv$ and bridges to $w$ in $\lat$.
Then $w \tau^{-1} \alpha \tau$ bridges to $w \tau$ in
$\lat$ and is in $\vv$. It follows that

$$w \tau^{-1} \alpha \tau \alpha^{1-m} 
\ = \  w \alpha \beta \ = \ 
w \beta = w \tau \alpha^{-1} \tau^{-1} \ \ \ {\rm is \ in}
\ \vv.$$

\noindent
Hence $w \tau \alpha^{-1}$ is in $\vv \tau$.
This implies

$$w \tau \alpha^{-1} \prec w \tau \prec w
\  \ \ \ \Rightarrow \ \  \ \ 
w \tau \prec w \tau \alpha \prec w \ \ \ \ \Rightarrow \ \  \ \ 
w \prec w \tau \alpha \tau^{-1} = 
w \beta^{-1} \prec w \tau^{-1}.$$

\noindent
In particular $w \beta^{-1}$ is in $\ro$ and $w \beta^{-1} \alpha^{-1}$
is in $\ro \alpha^{-1}$ which is not equal to $\vv$.
Also $w \beta^{-1} \alpha^{-1} = w \tau^{-1} \alpha^{-1}
\tau$. Here $w \tau^{-1} \alpha^{-1}$ is in $\ro \alpha^{-1}$ 
and bridges to $w$  in $\lat$ and so $w \tau^{-1} \alpha^{-1} \tau$
bridges to $w \tau$ in $\lat$ and so is in $\vv$.
As $\vv$ is not equal to $\ro \alpha^{-1}$, this
is a contradiction.

We conclude that situation II cannot happen either.
This finishes the proof of the lemma.
\end{proof}

Now we know that $w$ is not in $\lat$.

\begin{lemma}{}{}
$z$ is not in $\lal$.
\label{noal}
\end{lemma}

\begin{proof}{}
Suppose not, that is, $z \in \lal$.
This implies that either 
$z \alpha$ or $z \alpha^{-1}$ is in $(z,w)$.
Then lemma \ref{prel} implies that
$s   = z$.

Suppose first that $z \alpha \in (z,w]$.
So $z \alpha^{-1} \not \in T_z(w) = \uu$.
Use $z \alpha \tau = z \beta \alpha^m$
As $z \alpha \in \uu$, then $z \alpha \tau$ is
in $\uu$ also.
Then

$$z \alpha^{-1} \ \not \in \ \uu
\ \ \Rightarrow \ \  z \alpha^{-1} 
\tau^{-1} \ \not \in \ \uu \ \ \Rightarrow \ \ 
z \beta \ \ \ {\rm bridges \ to} \ \  [z,w] 
\ \ {\rm in} \ z$$

\noindent
 and
$z \beta \alpha^m$ bridges to $[z \alpha^m, w] \supset [z,w]$
in $z \alpha^m$. It follows that $z \beta \alpha^m$ is
not in $\uu$, contradiction.

Suppose now that $z \alpha^{-1}$ is in $[z,w]$. Then
$z \alpha^{-1} \tau^{-1} = z \beta$ is in $\uu$ and
bridges to $[z,w]$ in a point $t$ which is not
$z$. Then $z \beta \alpha^m$ bridges to $[z,w]$ in
$t \alpha^m$ and $z \beta \alpha^m$ is in $\uu$.
On the other hand $z \alpha$ is not in $\uu$ and
so $z \alpha \tau$ is not in $\uu$ either.
This is a contradiction.

This finishes the proof of the lemma.
\end{proof}

\vskip .1in
\noindent
{\bf {Summary in Case C.3}} $-$ So far we 
have proved: suppose
 that $w \alpha = w$, $s \kappa = s$,
no fixed points of $\kappa$ or $\alpha$  in $(s,w)$.
Let $z \in [s,w)$, the closest to $w$ with $z \tau = z$.
Then 

$$T_z(w) \tau \ = \ T_z(w), \ \ \ \ \ T_w(z) \alpha \ = \ T_w.$$

\noindent
If $\lat, \lal$ are the corresponding local axes 
of $\tau$ and $\alpha$ then
$z \not \in \lal$, $w \not \in \lat$.

\vskip .2in
\noindent
{\bf {Case C.3.0}} $-$ Suppose that $\lal \cap \lat$ has
at most one point.

This is very simple.
Let $[c,d]$ be the bridge from $\lat$ to $\lal$,
where $c = d$ if the intersection is one point.
We do the proof for $c \not = d$, the other is very 
similar.
Use $z \tau^{-1} \alpha \tau = z \alpha \beta \alpha^{m-1}$.
The right side is $z \alpha \tau$. Here $z \alpha$ bridges
to $\lal$ in $d \alpha$, hence bridges to $\lat$ in
$c$. So $z \alpha \tau$ bridges to $\lat$ in $c \tau$.

So $z \alpha \tau$ bridges to $\lal$ in $d$ so
$z \alpha \tau \alpha^{-1}$ bridges to $\lal$ in $d \alpha^{-1}$
and to $\lat$ in $c$. So $z \alpha \tau \alpha^{-1} \tau^{-1}
= z \alpha \beta$
bridges to $\lat$ in $c \tau^{-1}$ hence to 
$\lal$ in $c$. Finally $z \alpha \beta \alpha^{m-1}$ bridges
to $\lal$ in $d \alpha^{m-1}$ hence to $\lat$ in
$c$. Since $c \not = c \tau$ this is a contradiction.

\vskip .2in
\noindent
{\bf {Case C.3.1}} $-$
Now assume $\lal \cap \lat$ has more than one
point. 
We will use the analysis done in case B.

If $\uu \gamma$ is not equal $\uu$ then we use the
proof of case B.1.3 $-$ which was done also for
local axis of $\alpha$. This disallows this case.

The remaining case is that $\uu \gamma$ is equal to
$\uu$. As explained in case B.1.4 this implies
$\gamma$ leaves $\lat$ invariant.
Here we consider the intersection $\bb = \lal \cap \lat$.
First notice that $z$ is not in $\bb$. 
If $z$ were a limit point of $\bb$ then $\bb$ would
be $(z,r]$ (recall that $w$ is not in $\lat$).
Then as $\alpha$ leaves invariant $\lal$ we would
have $z \alpha = z$ also ruled out by non trivial
action of the group on $T$.
If $\lat$ is not properly embedded on the other
side let $v$ be the other ideal point of $\lat$.
Then 

$$v \kappa \ = \ v, \ \ \ (T_w(v)) \alpha \ = \ T_w(v),
\ \ \ (T_v(w)) \tau \ = \ T_v(w).$$

\noindent
Also $(w,v)$ has no fixed points of $\tau$.
Suppose that $v$ is in $\lal$. Then $(w,v)$ also
has no fixed points of $\alpha$. But then
$v$ has the same properties as $z$ and this
case is ruled out by lemma \ref{noal}.
It follows that $v$ is not in $\lal$. 
So if $\lat$ has another ideal point $v$, then
$\bb$ is $[r,t]$ with $t$ an actual point
in $\lat$.

Now we can apply the analysis of case B.1.4 which
was also done for $\alpha$ with a local axis.
The analysis rules out this situation.

This shows that case C.3.1 cannot happen either.

This finishes the proof of the main theorem.

\section{Remarks}

There are a lot of interesting questions still open.
First we discuss some internal questions about the proofs
in this article. The proof of the $\rrrr$-covered case
uses $p > 3q$ for $\alpha$ orientation reversing.
It would be useful to get a more general proof $-$
for instance showing that $p$ must be equal to $4$
or that $p$ has to be even. We obtained some preliminary
results, but not conclusive.
The same argument and condition $p > 3q$ are then used
in various places of the article so it would be
very good to discover a more general proof.

Also the best possible result for the manifolds $M_{p/q}$ 
described in this article would be the following:
If $p \geq q$, $p$ odd, $m \leq -4$ then the
only possible essential laminations are those 
coming from either stable or unstable lamination
in the original manifold $M$ $-$ these remain
essential whenever $|p - 2q| \geq 2$.
One way to interpret such a goal is a rigidity
result $-$ all laminations in this manifold have to
be of this type. Notice that Brittenham's results
for Seifert fibered spaces \cite{Br1} are of
this form. Also Hatcher and Thurston's results
for surgery on 2-bridge links \cite{Ha-Th} 
are along these lines.

Now on for more general goals: 
How far can the methods of this article be generalized?
Can they be used whenever $M$ is a punctured torus
bundle over $S^1$ with Anosov monodromy and
degeneracy locus $(1,2)$?
Probably a mixture of topological methods and group
action methods needs to be used.
How about surface bundles, where the surface has
higher genus? What about other degeneracy locus
as discovered by Gabai-Kazez \cite{Ga-Ka1}?

Since essential laminations do not exist in
every closed hyperbolic $3$-manifold, one
looks for useful generalizations. One possible 
idea was introduced by Gabai in \cite{Ga5}:
a lamination $\lambda$ in $M$, compact, orientable,
irreducible is {\em loosesse} if $\lambda$ satisfies:

0) \  $\lambda$ has no sphere leaves and

1) \ for any leaf $L$ of $\lambda$, the homomorphism
$\pi_1(L) \rightarrow \pi_1(M)$ induced by inclusion
is injective, and for any closed complementary region
$V$, the homomorphism $\pi_1(V) \rightarrow \pi_1(M)$
induced by inclusion is injective.

Gabai \cite{Ga5} conjectured that under
these conditions and $M$ closed then
$\widetilde \lambda$ is a product
lamination and $\mi$ is homeomorphic to $\rrrr^3$.
One test case is the class of manifolds $M_{p/q}$ studied
in this article.  When
$|p-2q| = 1$ the lamination coming from the stable
lamination has monogons. The leaves are either
planes or have ${\bf Z}$ fundamental group.
The complementary region is a solid torus.
So to check for loosesse one only needs to 
understand if leaves inject in the fundamental group
level.

{\footnotesize
{
\setlength{\baselineskip}{0.01cm}

\noindent
Florida State University

\noindent
Tallahassee, FL 32306-4510

}
}

\end{document}